%
%
%
\hsize=5in
\baselineskip=12pt
\vsize=20.6cm
\parindent=10pt
\pretolerance=40
\predisplaypenalty=0
\displaywidowpenalty=0
\finalhyphendemerits=0
\hfuzz=5pt
\frenchspacing
\footline={\ifnum\pageno=1\else\hfil\tenrm\number\pageno\hfil\fi}
%
%
\input amssym.def
\newsymbol\square 1003
\font\teneusm=eusm10
\font\seveneusm=eusm7
\font\fiveeusm=eusm5
\newfam\scriptfam
\textfont\scriptfam=\teneusm
\scriptfont\scriptfam=\seveneusm
\scriptscriptfont\scriptfam=\fiveeusm
\def\scr{\fam\scriptfam}
\font\tenib=cmmib10 
\skewchar\tenib='177
\def\boldfonts{\bf\textfont1=\tenib}
\font\ninerm=cmr9
\font\ninebf=cmbx9
\font\ninei=cmmi9
\skewchar\ninei='177
\font\ninesy=cmsy9
\skewchar\ninesy='60
\font\nineit=cmti9
\def\reffonts{\baselineskip=0.9\baselineskip
	\textfont0=\ninerm
	\def\rm{\fam0\ninerm}%
	\textfont1=\ninei
	\textfont2=\ninesy
	\textfont\bffam=\ninebf
	\def\bf{\fam\bffam\ninebf}%
	\def\it{\nineit}%
  \rm
}
%
%
\def\frontmatter{\vbox{}\vskip1cm\bgroup
	\leftskip=0pt plus1fil\rightskip=0pt plus1fil
	\parindent=0pt
	\parfillskip=0pt
	\pretolerance=10000
}
\def\endfrontmatter{\egroup\bigskip}
\def\title#1{{\baselineskip=1.44\baselineskip
  \font\titlef=cmbx12\titlef #1\par}}
\def\author#1{\bigskip#1\par}
\def\section#1\par{\ifdim\lastskip<\bigskipamount\removelastskip\fi
	\penalty-250\bigskip
	\vbox{\leftskip=0pt plus1fil\rightskip=0pt plus1fil
	\parindent=0pt
	\parfillskip=0pt
	\pretolerance=10000{\boldfonts#1}}\nobreak\medskip
}
\def\proclaim#1. {\medbreak\bgroup{\noindent\bf#1.}\ \it}
\def\endproclaim{\egroup
	\ifdim\lastskip<\medskipamount\removelastskip\medskip\fi}
\newdimen\itemsize
\def\setitemsize#1 {{\setbox0\hbox{#1\ }
	\global\itemsize=\wd0}}
\def\item#1 #2\par{\ifdim\lastskip<\smallskipamount\removelastskip\smallskip\fi
	{\leftskip=\itemsize
	\noindent\hskip-\leftskip
	\hbox to\leftskip{\hfil\rm#1\ }#2\par}\smallskip}
\def\Proof#1. {\ifdim\lastskip<\medskipamount\removelastskip\medskip\fi
	{\noindent\it Proof\if\space#1\space\else\ \fi#1.}\ }
\def\endproof{\hfill\hbox{}\quad\hbox{}\hfill\llap{$\square$}\medskip}
\def\Remark. {\ifdim\lastskip<\medskipamount\removelastskip\medskip\fi
        {\noindent\bf Remark. }}
\def\endremark{\medskip}
\def\emph#1{{\it #1}\/}
\def\text#1{\hbox{#1}}
\def\mathrm#1{{\rm #1}}
\def\refeq#1{\/$(#1)$}
%
%
\newcount\citation
\newtoks\citetoks
\def\citedef#1\endcitedef{\citetoks={#1\endcitedef}}
\def\endcitedef#1\endcitedef{}
\def\citenum#1{\citation=0\def\curcite{#1}%
	\expandafter\checkendcite\the\citetoks}
\def\checkendcite#1{\ifx\endcitedef#1?\else
	\expandafter\lookcite\expandafter#1\fi}
\def\lookcite#1 {\advance\citation by1\def\auxcite{#1}%
	\ifx\auxcite\curcite\the\citation\expandafter\endcitedef\else
	\expandafter\checkendcite\fi}
\def\cite#1{\makecite#1,\cite}
\def\makecite#1,#2{[\citenum{#1}\ifx\cite#2]\else\expandafter\clearcite\expandafter#2\fi}
\def\clearcite#1,\cite{, #1]}
%
%
\def\references{\section References\par
	\bgroup
	\parindent=0pt
	\reffonts
	\rm
	\frenchspacing
	\setbox0\hbox{99. }\leftskip=\wd0
	}
\def\endreferences{\egroup}
\newtoks\authtoks
\newif\iffirstauth
\def\checkendauth#1{\ifx\auth#1%
    \iffirstauth\the\authtoks
    \else{} and \the\authtoks\fi,%
  \else\iffirstauth\the\authtoks\firstauthfalse
    \else, \the\authtoks\fi
    \expandafter\nextauth\expandafter#1\fi
	}
\def\nextauth#1,#2;{\authtoks={#1 #2}\checkendauth}
\def\auth#1{\nextauth#1;\auth}
\newif\ifinbook
\newif\ifbookref
\def\nextref#1 {\par\hskip-\leftskip
	\hbox to\leftskip{\hfil\citenum{#1}.\ }%
	\initnextref}
\def\initnextref{\bookreffalse\inbookfalse\firstauthtrue\ignorespaces}
\def\paper#1{{\it#1},}
\def\InBook#1{\inbooktrue in ``#1",}
\def\book#1{\bookreftrue{\it#1},}
\def\journal#1{#1\ifinbook,\fi}
\def\BkSer#1{#1,}
\def\Vol#1{{\bf#1}}
\def\BkVol#1{Vol. #1,}
\def\publisher#1{#1,}
\def\Year#1{\ifbookref #1.\else\ifinbook #1,\else(#1)\fi\fi}
\def\Pages#1{\makepages#1.}
\long\def\makepages#1-#2.#3{\ifinbook pp. \fi#1--#2\ifx\par#3.\fi#3}
\def\inRus{{ \rm(in Russian)}}
\def\etransl#1{English translation in \journal{#1}}
%
%
\let\Rar\Rightarrow
\let\Lrar\Leftrightarrow
\let\hrar\hookrightarrow
\let\lhu\leftharpoonup
\let\rhu\rightharpoonup
\let\rhd\rightharpoondown
\let\ot\otimes
\let\sbs\subset
\let\<\langle
\let\>\rangle
\def\co#1{^{\mkern1mu\mathrm{co}\mkern2mu#1}}
\def\coev{\mathrm{coev}}
\def\cop{^{\mathrm{cop}}}
\def\ev{\mathrm{ev}}
\def\id{\mathrm{id}}
\def\lco#1{\vphantom{H}\co{#1}\mkern-1mu}
\def\op{^{\mathrm{op}}}
\def\triv{_{\mathrm{triv}}}
\def\lmapr#1#2{{}\mathrel{\smash{\mathop{\count0=#1
  \loop
    \ifnum\count0>0
    \advance\count0 by-1\smash{\mathord-}\mkern-4mu
  \repeat
  \mathord\rightarrow}\limits^{#2}}}{}}
\def\mapd#1#2{\llap{$\vcenter{\hbox{$\scriptstyle{#1}$}}$}\big\downarrow
  \rlap{$\vcenter{\hbox{$\scriptstyle{#2}$}}$}}
\def\lmapd#1#2#3{\llap{$\vcenter{\hbox{$\scriptstyle{#2}$}}$}
  \left\downarrow\vcenter to#1pt{}\right.\mskip-4mu
  \rlap{$\vcenter{\hbox{$\scriptstyle{#3}$}}$}}
\def\diagram#1{\vbox{\halign{&\hfil$##$\hfil\cr #1}}}
\def\cdiagram#1{\vcenter{\halign{&\hfil$##$\hfil\cr #1}}}
\let\al\alpha
\let\be\beta
\let\de\delta
\let\ep\varepsilon
\let\io\iota
\let\ka\kappa
\let\la\lambda
\let\ph\varphi
\let\si\sigma
\let\ze\zeta
\let\Ga\Gamma
\let\De\Delta
\let\Si\Sigma
\def\frm{{\frak m}}
\def\calA{{\cal A}}
\def\calE{{\cal E}}
\def\calF{{\cal F}}
\def\calM{{\cal M}}
\def\calR{{\cal R}}
\def\bbQ{{\Bbb Q}}
\def\bbR{{\Bbb R}}
\def\bbZ{{\Bbb Z}}
\def\ng{_{(-1)}}
\def\0{_{(0)}}
\def\1{_{(1)}}
\def\2{_{(2)}}
\def\3{_{(3)}}
\def\4{_{(4)}}
\def\AM{\vphantom{\calM}_A\mskip1mu\calM}
\def\BM{\vphantom{\calM}_B\mskip1mu\calM}
\def\Cd{C^\circ}
\def\CdQ{\vphantom{Q}^{\Cd}{\mskip-4mu}Q}
\def\CM{\vphantom{\calM}^C\!\calM}
\def\HAM{H\hbox{-}\AM}
\def\HAMA{\HAM_A}
\def\HAMB{\HAM_B}
\def\HBM{H\hbox{-}\BM}
\def\Hd{H^\circ}
\def\HdAM{\Hd\hbox{-}\mskip1mu\AM}
\def\HdAMA{\HdAM_A}
\def\HdAMC{\HdAM^\scrC}
\def\HdAMH{\HdAM_H}
\def\HdAMR{\HdAM_R}
\def\HdAMQ{\HdAM_Q}
\def\HdCM{\Hd\hbox{-}\mskip1mu\scrCM}
\def\HdCMH{\HdCM_H}
\def\HdCMR{\HdCM_R}
\def\HdHM{\Hd\hbox{-}\mskip1mu\HM}
\def\HdHMA{\HdHM_A}
\def\HdHMC{\HdHM^\scrC}
\def\HdHMQ{\HdHM_Q}
\def\HdM{\Hd\hbox{-}\mskip1mu\calM}
\def\HdMA{\HdM_A}
\def\HdMC{\HdM^\scrC}
\def\HdMQ{\HdM_Q}
\def\HdQM{\Hd\hbox{-}\mskip1mu\QM}
\def\HdQMH{\HdQM_H}
\def\HdQMQ{\HdQM_Q}
\def\HdRM{\Hd\hbox{-}\mskip1mu\RM}
\def\HdRMA{\HdRM_A}
\def\HdRMC{\HdRM^\scrC}
\def\HdRMQ{\HdRM_Q}
\def\Hdrat{\Hd_{\mathrm{rat}}\mskip1mu\hbox{-}\mskip1mu}
\def\HdratCM{\Hdrat\scrCM}
\def\HdratCMH{\HdratCM_H}
\def\HdratHMC{\Hdrat\HM^\scrC}
\def\HdratHMQ{\Hdrat\HM_Q}
\def\HdratMC{\Hdrat\MscrC}
\def\HdratMQ{\Hdrat\MQ}
\def\HdratQM{\Hdrat\mskip1mu\QM}
\def\HdratQMH{\HdratQM_H}
\def\HdratQMQ{\HdratQM_Q}
\def\HdxratAM{\Hd_{\mathrm{xrat}}\mskip1mu\hbox{-}\mskip1mu\AM}
\def\HdxratAMH{\HdxratAM_H}
\def\HdxratMA{\Hd_{\mathrm{xrat}}\mskip1mu\hbox{-}\mskip1mu\MA}
\def\HdxratHMA{\Hd_{\mathrm{xrat}}\mskip1mu\hbox{-}\mskip1mu\HM_A}
\def\HCM{\vphantom{\calM}_H^{\,C\!}\calM}
\def\HM{\vphantom{\calM}_H\calM}
\def\HMA{H\hbox{-}\mskip1mu\MA}
\def\HMB{H\hbox{-}\mskip1mu\MB}
\def\HMC{\HM^C}
\def\HMH{\HM^H}
\def\HMHH{\HMH_H}
\def\HMR{H\hbox{-}\mskip1mu\MR}
\def\HMscrC{H\hbox{-}\mskip2mu\MscrC}
\def\HQM{H\hbox{-}\QM}
\def\HRMR{H\hbox{-}\RM_R}
\def\LQ{\vphantom{Q}^LQ}
\def\LpQ{\vphantom{Q}^{L'}\!Q}
\def\MA{\calM_A}
\def\MAH{\MA^H}
\def\MB{\calM_B}
\def\MC{\calM^C}
\def\MH{\calM^H}
\def\MHH{\MH_H}
\def\Mk{\calM_k}
\def\MQ{\calM_Q}
\def\MR{\calM_R}
\def\MscrC{\calM^\scrC}
\def\QM{\vphantom{\calM}_Q\mskip1mu\calM}
\def\QMQ{\QM_Q}
\def\RM{\vphantom{\calM}_R\mskip1mu\calM}
\def\scrC{{\scr C}}
\def\scrCM{\vphantom{\calM}^\scrC\!\calM}
\def\scrI{{\scr I}}
\def\sqC{{}\mathbin{\square}_C}
\def\sqCH{\sqC\!H}
\def\Urd{\vphantom{U}^*{\mskip-1mu}U}
\def\defop#1#2{\def#1{\mathop{\rm #2}\nolimits}}
\defop\chr{char}
\defop\End{End}
\defop\gr{gr}
\defop\Hom{Hom}
\defop\Id{Id}
\defop\Img{Im}
\defop\Ker{Ker}
\defop\Max{Max}
\defop\rann{rann}
\defop\soc{soc}
\defop\Indinj{Ind.inj}
\defop\Indproj{Ind.proj}
\defop\Irr{Irr}
\defop\Lat{Lat}
\defop\Rat{Rat}
\defop\adc{adc}
\defop\lrk{lrk}
\defop\rrk{rrk}
\def\gr{_\mathrm{gr}}
\citedef
Ab-K59
Am-M05
An92
And-F
Az54
Bas60
Ber-Pl
Bia61
Bich23
Bj71
Br-W
Cae-MZ
Ch14
Ch24
Cl-PS77
DG
Di76
Doi83
Doi92
Er-Sk09
Fai
Fi33
Gab70
Gant
Kop93
Mes06
Moe-R86
Mo
Mue70
Nas-BO89
Nich78
Ob77
Scha94
Schn93
Scho
Sk06
Sk07
Sk10
Sk11
Sk21
Sk25a
Sk25b
Sk-Oy06
St
Sw
Sw75
Tak72
Tak79
Tak89
Tak94
Ulb90
Von96
\endcitedef

\frontmatter
\title{Hopf algebraic homogeneous spaces interpreted rationally:\break
the Abe-Kanno theorem}
\author{Serge Skryabin}
\endfrontmatter

\bigskip
{\reffonts
\leftskip=20pt\rightskip=20pt
{\bf Abstract.}
We present a Hopf algebraic generalization of the Abe-Kanno theorem on a 
correspondence between subgroups of an algebraic group and invariant subfields 
of the field of rational functions. It applies to residually finite-dimensional 
Hopf algebras admitting an artinian classical quotient ring and is used in the 
paper to derive some general properties of such Hopf algebras.\par}

\section
Introduction

Let $k[G]$ be the commutative Hopf algebra of regular functions on an affine 
algebraic group $G$. The right coideal subalgebras of $k[G]$ are precisely the 
subalgebras invariant under the action of $G$ on $k[G]$ induced by right 
translations of $G$. Each right coideal subalgebra $A$ determines a closed 
group subscheme of $G$, the stabilizer of the distinguished maximal ideal 
$A^+$ of $A$ consisting of functions vanishing at the identity element of $G$.  
This subscheme is represented by the factor Hopf algebra 
$\,k[G]\,/{\mskip1mu}k[G]{\mskip1mu}A^+\!$.

If $A$ is a right coideal subalgebra of a noncommutative Hopf algebra $H$, 
then $HA^+$ may not be a two-sided ideal, and therefore $H/HA^+$ is a factor 
coalgebra of $H$ with inherited left $H$-module structure, but not necessarily 
a Hopf algebra. This observation suggests that left $H$-module factor 
coalgebras of $H$ should be viewed as analogs of closed subgroups of an 
algebraic group.

The assignment $A\mapsto H/HA^+$ gives a map from the set of right coideal 
subalgebras of $H$ to the set of left $H$-module factor coalgebras of $H$. In 
the opposite direction Takeuchi \cite{Tak79} associated with a left $H$-module 
factor coalgebra $C$ the right coideal subalgebra $\lco CH$ consisting of 
all elements of $H$ invariant under the induced left coaction of $C$ on $H$. 
Several results on this correspondence were proved in \cite{Tak79} under the 
assumptions of faithful flatness or coflatness.

In general the assignments $A\mapsto H/HA^+$ and $C\mapsto\lco CH$ define a 
Galois connection between those two sets, suitably ordered. Chirvasitu 
\cite{Ch24} characterized closed elements of this connection in purely 
categorical terms as dominion subalgebras in the first set and codominion 
factor coalgebras in the second. It is very clear that the two maps are not 
bijective in general. If $K$ is the closed subgroup of an affine algebraic 
group $G$ corresponding to a right coideal subalgebra $A$ of $k[G]$, then the 
quotient $K\backslash G$ is a quasiaffine algebraic variety since its 
points are separated by functions in $A$. However, arbitrary homogeneous 
spaces for $G$ are quasiprojective, not necessarily quasiaffine, algebraic 
varieties.

Thus right coideal subalgebras of a Hopf algebra $H$ should be considered as 
objects representing Hopf algebraic analogs of quasiaffine homogeneous 
spaces. A more recent paper \cite{Sk10} adopted the view that Hopf algebraic 
quasiprojective homogeneous spaces should be interpreted by means of graded 
right $H$-comodule subalgebras of the Laurent polynomial algebra 
$H[t,t^{-1}]$. Several results on correspondence between left $H$-module 
factor coalgebras of $H$ and graded subalgebras of $H[t,t^{-1}]$ were proved 
in \cite{Sk10} for a Hopf algebra $H$ over a field $k$ satisfying the 
following two assumptions:

\setitemsize(A2)
\item(A1)
\enspace
$H$ is a residually finite-dimensional,

\item(A2)
\enspace
$H$ has an artinian classical quotient ring $Q(H)$.

An algebra over a field is called \emph{residually finite-dimensional} if its 
ideals of finite codimension have zero intersection. An overring $Q(R)$ of a 
ring $R$ is a \emph{classical right quotient ring} of $R$ if all 
nonzerodivisors of $R$ are invertible in $Q(R)$ and each element of $Q(R)$ 
can be written as $as^{-1}$ where $a,s\in R$ and $s$ is a nonzerodivisor. Such 
a ring exists if and only if the set $\Si$ of all nonzerodivisors of $R$ 
satisfies the right Ore condition, and then $Q(R)$ is isomorphic to the ring 
of fractions $R\Si^{-1}$. If $\Si$ satisfies both the right and the left Ore 
conditions, then $Q(R)$ is a \emph{classical (two-sided) quotient ring}.

Unlike what happens to be in the case of commutative Hopf algebras, it is not 
clear, and probably not true in general, that each left $H$-module factor 
coalgebra of $H$ corresponds to a model of a quasiprojective homogeneous space 
in the sense of \cite{Sk10}. This raises the question as to how the remaining 
left $H$-module factor coalgebras should be added to the general picture.

Rational functions provide richer information about an algebraic group and its 
quotients than just regular functions do. For example, rational functions 
separate cosets of any closed subgroup. This point of view allows us to remove 
previous restrictions in the Hopf algebraic version of the correspondence 
between subgroups and homogeneous spaces. Indeed, the quotient ring $Q(H)$ is 
the precise analog of the ring $k(G)$ of rational functions on an algebraic 
group $G$. There are natural left and right actions of the dual Hopf algebra 
$\Hd$ on $H$, both of which extend to an action on $Q(H)$. We may consider 
subalgebras of $Q(H)$ stable under the left action of $\Hd$. Now we state the 
main result of the present paper:

\proclaim
Theorem 0.1.
For a Hopf algebra $H$ satisfying assumptions {\rm(A1)} and {\rm(A2)} there 
is a canonical bijection between the set of left $H$-module factor coalgebras 
of $H$ and the set of left $\Hd$-invariant artinian subalgebras of\/ $Q(H)$.
\endproclaim

A theorem of Abe and Kanno \cite{Ab-K59, Th.~1} states that for a connected 
algebraic group $G$ over an algebraically closed field $k$ there is a 
bijection between the set of closed subgroups of $G$ and the set of right 
invariant subfields of $k(G)$ which contain the field $k$ and over which 
$k(G)$ is separably generated (see also \cite{Bia61}). Here right invariance 
means that a subfield is stable under the left action of $G$ on $k(G)$ induced 
by right translations of $G$. This action of $G$ is the restriction of the 
left action of the dual Hopf algebra $\Hd$ for $H=k[G]$, and we call such 
subfields left $\Hd$-invariant.

In the case when $\,\chr k>0\,$ the second set in the Abe-Kanno theorem does 
not include all right invariant subfields of $k(G)$ containing the field $k$. 
This is not surprising since all closed subgroups considered in \cite{Ab-K59} 
are reduced schemes, while there exist group schemes with nilpotents in their 
structure sheaves. Retaining the old perception of the notion of closed 
subgroups, Vonessen revisited that result and showed that each right invariant 
subfield of $k(G)$ containing the field $k$ is a purely inseparable extension 
of the subfield corresponding to some closed subgroup of $G$ \cite{Von96, 
Prop. 3.4}. The Hopf algebraic language makes no distinction between reduced 
and nonreduced algebras, nor is the algebraic closedness of the base field 
needed. However, Theorem 0.1 generalizes only the case of an affine group $G$ 
in the Abe-Kanno theorem.

In another paper \cite{Tak89} Takeuchi developed a Hopf algebraic approach to 
the Picard-Vessiot theory of field extensions. It was later extended by Amano 
and Masuoka \cite{Am-M05} to extensions of artinian $D$-simple $D$-module 
algebras for a smooth pointed cocommutative Hopf algebra $D$. Curiously, 
neither \cite{Tak89} nor \cite{Am-M05} mention the quotient ring $Q(H)$ of a 
finitely-generated commutative Hopf algebra $H$ as an example of a Hopf 
algebraic Picard-Vessiot extension of the base field $k$. Several results of 
these two papers apply to $Q(H)=k(G)$ in the case when $H=k[G]$ for a reduced 
affine algebraic group $G$, at least when $k$ is algebraically closed. In 
particular, a special case of \cite{Tak89, Th. 2.10} shows that for a connected 
$G$ the right invariant subfields of the field $k(G)$ are in a bijective 
correspondence with the closed group subschemes of $G$, thus strengthening the 
Abe-Kanno theorem.

All finitely-generated commutative (associative and unital) $k$-algebras are 
residually finite-dimensional \cite{Sw, Th. 6.1.3}. Since a finitely-generated 
commutative Hopf algebra is a Cohen-Macaulay ring \cite{Gab70, Prop. 1.1.1}, 
its total ring of fractions is artinian, i.e., such a Hopf algebra satisfies 
both (A1) and (A2).

However, conditions (A1) and (A2) hold for a much larger class of Hopf algebras, 
even in the commutative case. If $H$ is any commutative Hopf algebra over an 
algebraically closed field $k$ and $A$ its Hopf subalgebra, then each algebra 
homomorphism $A\to k$ extends to $H$ by \cite{DG, Ch. III, \S3, Cor. 7.6} or 
\cite{Tak72, Cor. 3.12}. If $H$ is moreover reduced, then it follows that its 
ideals of codimension 1 have zero intersection since this holds for all 
finitely-generated Hopf subalgebras of $H$ by Hilbert's Nullstellensatz. In 
particular, any commutative Hopf domain over an algebraically closed field 
satisfies (A1) and (A2).

Our main interest, however, is to generalize previously known results in the 
commutative case to noncommutative Hopf algebras. By \cite{Sk21, Th. 5.5} each 
residually finite-dimensional noetherian Hopf algebra satisfies also condition 
(A2), and so our results in the present paper apply to several important 
classes of noncommutative Hopf algebras.

The subfield of $k(G)$ corresponding to a closed subgroup $K$ of $G$ in the 
Abe-Kanno theorem is defined very easily as the set of elements of $k(G)$ 
invariant with respect to the action of $K$ induced by left translations of 
$G$. A similar construction is not available in the general case of Theorem 0.1. 
We associate a left $\Hd$-invariant artinian subalgebra of $Q(H)$ to a left 
$H$-module factor coalgebra of $H$ in two steps via an intermediate 
correspondence with a certain class of corings. This is very similar to the 
use of corings in the Hopf algebraic Picard-Vessiot theory, but we encounter 
more complications in the proofs. Section 1 of this paper explains our 
construction in detail and gives an overview of several parts of the proof 
contained in other sections.

One may wonder whether assumptions (A1) and (A2) imply other fundamental 
properties known for commutative Hopf algebras representing algebraic groups. 
A natural question concerns flatness over right coideal subalgebras. We are 
not able to answer it, yet Theorem 9.3 shows that flatness does hold for a 
right coideal subalgebra $A$ of a Hopf algebra $H$ satisfying (A1) and (A2) 
when the inclusion map $A\to\nobreak H$ admits a right or left $A$-linear 
retraction. By Theorem 10.1 all biideals of such a Hopf algebra $H$ are Hopf 
ideals. Verification of these facts makes use of Theorem 0.1 in an essential 
way. Concluding sections of the paper provide further information on the 
correspondence of Theorem 0.1.

$$
\vcenter{\baselineskip12pt
\halign{\hfil\ninerm#&\ninerm\enspace#\hfil\cr
\multispan2\hfil{\bf Contents}\hfil\cr
\noalign{\smallskip}
1.&Description of the correspondence\cr
2.&Simple module algebras and freeness of equivariant modules\cr
3.&First properties of the quotient ring and invariant subalgebras\cr
4.&Correspondence between corings and coalgebras\cr
5.&Coaction invariants and the second bijection\cr
6.&Twisting of comodules over $H$-module corings\cr
7.&A coring structure theorem\cr
8.&Equivalences between module and comodule categories\cr
9.&Exactness of induction and flatness over coideal subalgebras\cr
10.&Biideals are Hopf ideals\cr
11.&Conormal factor coalgebras\cr
12.&Comparison with the classical case\cr
13.&Quasiprojective homogeneous spaces\cr
\multispan2\ninerm
Appendix.\enspace Semiprimary rings arising naturally\hfill\cr
}}
$$

\section
1. Description of the correspondence

As pointed out in the introduction, the bijective correspondence of Theorem 0.1 
involves a certain class of corings as an intermediate link. A coring over a 
ring $R$ is an $R$-bimodule $\scrC$ equipped with a coassociative bimodule 
homomorphism
$$
\De_\scrC:\scrC\to\scrC\ot_R\scrC,
$$
called the comultiplication of $\scrC$, for which there exists a counit 
$\ep_\scrC:\scrC\to R$, also a bimodule homomorphism. In brief wording an 
$R$-coring is a coalgebra in the monoidal category of $R$-bimodules. For 
general facts concerning corings and comodules we refer to the book of 
Brzezi\'nski and Wisbauer \cite{Br-W}.

A coideal of an $R$-coring $\scrC$ is an $R$-subbimodule $\scrI$ such that 
$\ep_\scrC(\scrI)=0$ and $\De_\scrC(\scrI)$ is contained in the image of the 
canonical map
$$
(\scrI\ot_R\scrC)\oplus(\scrC\ot_R\scrI)\to\scrC\ot_R\scrC.
$$
For each coideal $\scrI$ the factor coring $\scrC/\,\scrI$ is defined in such 
way that the canonical surjection $\scrC\to\scrC/\,\scrI$ is a homomorphism of 
corings.

For each subring $A$ of a ring $R$ the $R$-bimodule $\scrC=R\ot_AR$ has a 
canonical coring structure. Assuming that $\scrC\ot_R\scrC$ is identified with 
$R\ot_AR\ot_AR$, the structure maps for this coring are defined by the 
formulas
$$
\De_\scrC(x\ot_Ay)=x\ot_A1\ot_Ay,\qquad\ep_\scrC(x\ot_Ay)=xy\eqno(1.1)
$$
for $x,y\in R$.

Further on we fix a base field $k$. The subscript $k$ in the notation for the 
tensor product functor $\ot_k$ will be omitted. Let $H$ be a Hopf algebra 
over $k$. We write $\De$, $\ep$, $S$ for the comultiplication, the counit, and 
the antipode of $H$.

A left $H$-module algebra $A$ is an algebra (associative and unital) in the 
monoidal category $\HM$ of left $H$-modules (see \cite{Mo} and \cite{Sw}). We 
denote by $\HMA$ and $\HAM$ the categories of \emph{$H$-equivariant} right and 
left $A$-modules. An object $M$ of $\HMA$ has a right $A$-module and a left 
$H$-module structures such that the structure map $M\ot A\to M$ is a morphism 
in the monoidal category $\HM$. In the case of the category 
$\HAM\,$ left $A$-module structures $A\ot M\to M$ are used instead.  One can 
interpret objects of $\HMA$ and $\HAM$ as modules over certain smash products. 
Given $M\in\HMA$ and $N\in\HAM$, the tensor product $M\ot_AN$ has an 
$H$-module structure which makes $M\ot_AN$ a factor module of the tensor 
product $M\ot N$ in the monoidal category $\HM$.

For two left $H$-module algebras $A$ and $B$ we denote by $\HAMB$ the category 
of \emph{$H$-equivariant} $(A,B)$-bimodules. An object $M$ of $\HAMB$ is a 
bimodule equipped with a left $H$-module structure with respect to which $M$ 
is an object of both $\HAM$ and $\HMB$. Tensoring with such a bimodule 
produces functors
$$
?{}\ot_AM:\,\HMA\to\HMB\qquad\text{and}\qquad M\ot_B{}?:\,\HBM\to\HAM.
$$
In particular, the category $\HAMA$ is monoidal with respect to the functor 
$\,\ot_A\,$.

\smallskip
Let $R$ be a left $H$-module algebra. By a \emph{left $H$-module $R$-coring} 
we understand an $R$-coring $\scrC$ equipped with a left $H$-module structure 
with respect to which $\scrC$ is an object of the category $\HRMR$ and the 
structure maps $\De_\scrC$, $\ep_\scrC$ are $H$-linear, i.e., $\De_\scrC$ and 
$\ep_\scrC$ are morphisms in the monoidal category $\HRMR$. In other words, a 
left $H$-module $R$-coring is a coalgebra in $\HRMR$. This generalizes the 
notion of a left $H$-module coalgebra which is a coalgebra in the monoidal 
category $\HM$ \cite{Mo}. If $A$ is an $H$-invariant subalgebra of $R$, then 
the canonical coring $\,R\ot_AR$ is a left $H$-module coring.

\smallskip
Suppose next that $H$ satisfies basic assumptions (A1) and (A2). Put $Q=Q(H)$. 
Already the assumption that $H$ has a right artinian classical 
right quotient ring $Q$ implies that the antipode $S$ of $H$ is bijective by 
\cite{Sk06, Th.~A}, and therefore, since $H\op$ is antiisomorphic to $H$, the 
ring $Q$ is automatically left artinian and also a left quotient ring of $H$. 
Thus the one-sided conditions on $Q$ are anyway equivalent to the two-sided 
conditions used in (A2).

Bijectivity of $S$ allows us to apply all results proved for $H$ to the Hopf 
algebras $H\op$ and $H\cop$ in which either the multiplication or the 
comultiplication of $H$ is changed to the opposite one.

The dual Hopf algebra $\Hd$ consists of all linear functions $H\to k$ vanishing 
on an ideal of finite codimension in $H$. The left and right actions of $\Hd$ 
on $H$ are given by the formulas
$$
f\rhu h=\sum\,f(h\2)h\1\,,\qquad h\lhu f=\sum\,f(h\1)h\2\eqno(1.2)
$$
where $f\in\Hd$ and $h\in H$. These actions make $H$ into a left $\Hd$-module 
algebra and a right one. By \cite{Sk-Oy06, Th. 2.2} the left action extends to 
$Q$, and the $H\cop$-variant of that result shows that so does the right action 
too. In this way $Q$ becomes an $\Hd$-bimodule algebra.

The two extended actions of $\Hd$ on $Q$ will be denoted by the same symbols 
$\rhu$ and $\lhu\,$. A subalgebra $A$ of $Q$ will be called \emph{left 
$\Hd$-invariant} (respectively, \emph{right $\Hd$-invariant}) if 
$f\rhu x\in A$ (respectively, $x\lhu f\in A$) for all $f\in\Hd$ and $x\in A$. 
We will use mostly the left action.

The base field $k$ is a subalgebra of $Q$ on which $\Hd$ acts trivially. 
Consider the canonical left $\Hd$-module $Q$-coring $Q\ot Q$ associated with 
the ring extension $k\sbs Q$. Its left $\Hd$-module factor corings 
$(Q\ot Q)/\,\scrI$ correspond to left $\Hd$-invariant coideals $\scrI$ of 
$Q\ot Q$. Note that $Q\ot Q$ is also a right $\Hd$-module coring with respect 
to the action $\lhu$ of $\Hd$ on each tensorand. This module structure on 
$Q\ot Q$ will be used only in section 11 of the paper where we encounter 
coideals of $Q\ot Q$ which are invariant under both actions of $\Hd$.

Left $H$-module factor coalgebras $H/I$ of $H$ correspond to coideals of $H$ 
which are also left ideals. Now we refine Theorem 0.1 in the following 
statement:

\setitemsize$\Lat_3$
\proclaim
Theorem 1.1.
For a Hopf algebra $H$ satisfying assumptions {\rm(A1)} and {\rm(A2)} there are 
canonical isomorphisms between the following three lattices:
$$
\eqalign{
\Lat_1\mskip1mu,&\quad\text{the lattice of coideal left ideals of $H,$}\cr
\Lat_2\mskip1mu,&\quad\text{the lattice of left $\Hd$-invariant coideals of 
the $Q$-coring $Q\ot Q,$}\cr
\Lat_3\mskip1mu,&\quad\text{the lattice of left $\Hd$-invariant artinian subalgebras 
of $Q$.}
}
$$
Thus the left $H$-module factor coalgebras of $H,$ the left $\Hd$-module 
factor corings of the left $\Hd$-module $Q$-coring $Q\ot Q,$ and the 
left $\Hd$-invariant artinian subalgebras of $Q$ are in a bijective 
correspondence with each other.
\endproclaim

The map $\Lat_1\to\Lat_2$ is defined by the assignment $I\mapsto\scrI$ 
described in the next lemma.

\proclaim
Lemma 1.2.
Define a $k$-linear map $\,\psi:H\to H\ot H\sbs Q\ot Q\,$ by the rule
$$
\psi(x)=\sum x\1\ot S(x\2),\qquad x\in H.
$$
For each coideal left ideal $I$ of $H$ the $Q$-subbimodule $\scrI$ of $Q\ot Q$ 
generated by the set $\psi(I)$ is a left $\Hd$-invariant coideal of the 
$Q$-coring $Q\ot Q$.
\endproclaim

\Proof.
Viewing $H$ as a right $H$-comodule with respect to the comultiplication of 
$H$, let $\,\rho:H\ot H\to(H\ot H)\ot H\,$ be the comodule structure on the 
tensor product of two copies of $H$ in the monoidal category $\MH$ of right 
$H$-comodules. Then
$$
\rho\bigl(\psi(x)\bigr)=\sum\,\bigl(x\1\ot S(x\4)\bigr)\ot x\2S(x\3)
=\sum\,\bigl(x\1\ot S(x\2)\bigr)\ot1=\psi(x)\ot1.
$$
This means that for each $x\in H$ the element $\psi(x)\in H\ot H$ is invariant 
under the coaction of $H$, and therefore under the corresponding action of the 
dual Hopf algebra $\Hd$. Hence $\scrI$ is stable under the left action of $\Hd$. 
For each $x\in I$ we have
$$
\ep_{Q\ot Q}\bigl(\psi(x)\bigr)=\sum x\1S(x\2)=\ep(x)1=0
$$
since $\ep(I)=0$, and, identifying $\,(Q\ot Q)\ot_Q(Q\ot Q)\,$ with 
$\,Q\ot Q\ot Q$,
$$
\eqalign{
\De_{Q\ot Q}\bigl(\psi(x)\bigr)
&=\sum\,x\1\ot1\ot S(x\2)\cr
&=\sum\,x\1\ot S(x\2)x\3\ot S(x\4)\cr
&=\sum\,\bigl(x\1\ot S(x\2)\bigr)\ot_Q\bigl(x\3\ot S(x\4)\bigr)\cr
&=\sum\,\psi(x\1)\ot_Q\psi(x\2)\in\scrI\ot Q+Q\ot\scrI
}
$$
since $\sum x\1\ot x\2=\De(x)\in I\ot H+H\ot I$. This shows that $\scrI$ is a 
coideal.
\endproof

Bijectivity of the map $\Lat_1\to\Lat_2$ is a special case 
of Corollary 4.8, taking into account Proposition 4.10. The inverse map 
$\Lat_2\to\Lat_1$ is given by the assignment
$$
\scrI\mapsto\psi^{-1}(\,\,\scrI^{\Hd})
$$
where $\,\scrI^{\Hd}$ stands for the vector subspace of $\Hd$-invariants in 
$\scrI$.

For each left $\Hd$-invariant subalgebra $A$ of $Q$ the canonical map 
$$
Q\ot Q\to Q\ot_AQ
$$
is a surjective homomorphism of left $\Hd$-module $Q$-corings whose kernel 
$\,\scrI_A$ is the $Q$-subbimodule of $Q\ot Q$ generated by the set 
$\{1\ot a-a\ot1\mid a\in A\}$. Hence $\,\scrI_A$ is a left $\Hd$-invariant 
coideal of the $Q$-coring $Q\ot Q$ and
$$
Q\ot_AQ\cong(Q\ot Q)/\,\scrI_A.
$$
We always identify $Q\ot_AQ$ with a left $\Hd$-module factor coring of $Q\ot Q$ 
in this way. The map $\,\Lat_3\to\Lat_2\,$ is defined by the assignment 
$\,A\mapsto\scrI_A$.

To explain the map in the opposite direction we recall the notion of 
dominions. Given a ring $R$, the \emph{dominion} of its subring $A$ is the 
larger subring consisting of all elements $x\in R$ such that
$$
x\ot_A1=1\ot_Ax\quad\text{in $R\ot_AR$}
$$
(see \cite{St, Ch. XI, Prop. 1.1}). We say that $A$ is a \emph{dominion subring} 
if $A$ coincides with its own dominion in $R$.

\proclaim
Lemma 1.3.
For each left $\Hd$-invariant coideal $\,\scrI$ of the $Q$-coring $Q\ot Q$ 
the set
$$
A=\{x\in Q\mid\,1\ot x-x\ot1\in\scrI\,\}
$$
is a left $\Hd$-invariant dominion subalgebra of $Q$.
\endproclaim

\Proof.
We have $A=\de^{-1}(\,\scrI\,)$ where $\de:Q\to Q\ot Q$ is the $k$-linear map 
defined by the rule $\,\de(x)=1\ot x-x\ot1\,$ for $x\in Q$. If $a,b\in A$, then
$$
\de(ab)=a\de(b)+\de(a)b\in\scrI
$$
since $\scrI$ is a subbimodule of $Q\ot Q$. Also $1\in A$ since $\de(1)=0$. 
This shows that $A$ is a subalgebra of $Q$. It is stable under the left action 
of $\Hd$ on $Q$ since $\de$ is actually a homomorphism of left $\Hd$-modules.

If $x\in Q$ is any element such that $\,1\ot_Ax=x\ot_A1\,$ in $Q\ot_AQ$, then 
$1\ot x-x\ot1$ lies in $\scrI_A$, but $\scrI_A\sbs\scrI$ by the definition of 
$A$ and $\scrI_A$, whence $x\in A$. Thus $A$ coincides with its dominion in $Q$.
\endproof

As an application of more general results proved in section 2 it will be 
asserted in Proposition 3.6 that the left $\Hd$-invariant dominion subalgebras 
of $Q$ are precisely the left $\Hd$-invariant artinian subalgebras of $Q$. 
Thus the algebra $A$ in Lemma 1.3 is artinian, and we define the map 
$\,\Lat_2\to\Lat_3\,$ by assigning this subalgebra to the left $\Hd$-invariant 
coideal $\,\scrI$.

If $A$ is any dominion subalgebra of $Q$, then the element $1\ot x-x\ot1$ 
lies in the kernel $\scrI_A$ of the canonical coring homomorphism 
$Q\ot Q\to Q\ot_AQ$ if and only if $x\in A$. Hence the map $\,\Lat_2\to\Lat_3\,$ 
is a left inverse of the map $\,\Lat_3\to\Lat_2\mskip1mu$.

For the bijectivity of the maps $\,\Lat_2\to\Lat_3\,$ and $\,\Lat_3\to\Lat_2\,$ 
we have to know that each left $\Hd$-invariant coideal of $Q\ot Q$ 
is generated as a $Q$-subbimodule by some subset of the set
$$
\{1\ot x-x\ot1\mid x\in Q\},
$$
or, in other words, each left $\Hd$-module factor coring of $Q\ot Q$ is 
$Q\ot_AQ$ for some subalgebra $A$ of $Q$. In the case when $Q$ is a division 
ring, i.e., a skew field, this conclusion follows immediately from Sweedler's 
fundamental theorem in \cite{Sw75} which determines all coideals of such 
corings. A short argument in section 5 generalizes this conclusion to the case 
where all simple factor rings of the artinian ring $Q$ are skew fields.

Unfortunately this argument does not work in general, and we need a 
considerable amount of extra work to complete the proof. This is done in a 
more general setup of section 7. The conclusion that $\,\scrC=Q\ot_AQ\,$ under 
suitable assumptions about a coring $\scrC$ is presented in Theorem 7.11. A 
special case of this theorem stated in Proposition 5.5 provides the final 
ingredient in the proof of Theorem 1.1.

We will refer to the bijection between $\Lat_1$ and $\Lat_2$ as the 
\emph{first bijection} of Theorem 1.1. The bijection between $\Lat_2$ and 
$\Lat_3$ is the \emph{second bijection}. Combining them we obtain the 
bijection of Theorem 0.1. The left $H$-module factor coalgebra $C=H/I$ of $H$ 
corresponding to a left $\Hd$-invariant artinian subalgebra $A$ of $Q$ can be 
described in this way:
$$
\displaylines{
C\cong(Q\ot_AQ)^{\Hd},\cr
\noalign{\smallskip}
I=\{\,x\in H\mid\ \sum\,x\1\ot_AS(x\2)=0\ \text{ in $\,Q\ot_AQ\,$}\}.
}
$$
As explained in section 4, the $Q$-coring structure on $\,Q\ot_AQ$ induces a 
coalgebra structure on the subspace of $\Hd$-invariants $\,(Q\ot_AQ)^{\Hd}\!$.

In two special cases considered in sections 12 and 13 of the paper we provide 
a more direct description of the left $\Hd$-invariant artinian 
subalgebra $A$ of $Q$ corresponding to a 
given left $H$-module factor coalgebra $C$ of $H$. However, intervention 
of corings seems to be unavoidable in the general case of our treatment.

The lattices $\,\Lat_1$, $\Lat_2$, $\Lat_3\,$ are complete. The supremum of 
a family of coideals of $H$ is the sum of coideals, and the infimum is the 
largest coideal contained in each coideal of the family. Supremums and 
infimums in $\,\Lat_2$ are described similarly. The supremum of a family of 
$\Hd$-invariant artinian subalgebras of $Q$ is the dominion of the subalgebra 
generated by all subalgebras in the family. The infimum is the intersection of 
subalgebras. The fact that this intersection is artinian follows again from 
Proposition 3.6 since the set of dominion subrings of a ring is closed under 
arbitrary intersections \cite{Scho, Lemma 7.20}.

\smallskip
Several categories of objects equipped with module or comodule structures 
satisfying some additional conditions will be used throughout the whole paper. 
A guide to sections where these categories are introduced is given below:
$$
\vcenter{\advance\baselineskip by3pt
\halign{\hfil#&\enspace#\hfil\cr
section 1:& $\HMA$, $\HAM$, $\HAMB$\cr
section 3:& $\HdratMQ$, $\HdratQM$\cr
section 4:& $\HdratQMQ$, $\HdMC$, $\HdratMC$, $\HdCM$, $\HdratCM$\cr
section 6:& $\HMscrC$\cr
section 8:& $\HdxratMA$, $\HdxratAM$, $\HdxratHMA$, $\HdxratAMH$,\cr
\omit& $\HdRMC$, $\HdCMR$, $\HdratHMQ$, $\HdratHMC$, $\HdratQMH$, $\HdratCMH$\cr
}}
$$
Standard notation $\HMC$ and $\HCM$ will be used for the categories of 
relative Hopf modules introduced by Takeuchi \cite{Tak79} and Doi 
\cite{Doi92}. Morphisms in each of those categories are maps respecting all 
relevant module and comodule structures.

The correspondence of Theorem 0.1 manifests itself also on the level of module 
and comodule categories. In section 8 it will be shown that the categories 
$\MC$ and $\CM$ of right and left $C$-comodules are equivalent to $\HdxratMA$ 
and $\HdxratAM$, respectively. The categories $\HMC$ and $\HCM$ are equivalent 
to $\HdxratHMA$ and $\HdxratAMH$. This will be important for the study of 
further properties in later sections of the paper.

With the exception of results in sections 2, 6, 7 developed in a more general 
context it is assumed elsewhere in this paper that $H$ is a Hopf algebra 
satisfying assumptions (A1) and (A2).

\section
2. Simple module algebras and freeness of equivariant modules

Let $H$ be an arbitrary Hopf algebra over the base field $k$ and $A$ a left 
$H$-module algebra. We say that an object $M$ of either $\HMA$ or $\HAM$ is 
\emph{$A$-finite} if $M$ is finitely generated as an $A$-module, and $M$ is 
\emph{locally $A$-finite} if $M$ is the union of its $A$-finite subobjects. An 
$H$-module algebra is \emph{$H$-simple} if it has no nonzero proper $H$-stable 
ideals. One important argument that we need is provided by \cite{Sk07, Th. 
7.6} which we repeat here:

\proclaim
Theorem 2.1.
Suppose that $A$ is a semilocal $H$-simple left $H$-module algebra. Then every 
locally $A$-finite object $M\in\HMA$ is projective in $\MA$. Moreover, $M$ is 
a free $A$-module if and only if $M/M\frm$ is a free $A/\frm$-module for at 
least one maximal ideal\/ $\frm$ of $A$.
\endproclaim

If the antipode of $H$ is bijective, then a similar result for objects of 
$\HAM$ is obtained by replacing $H$ with $H\cop$.

Any application of Theorem 2.1 requires the $H$-simplicity of $A$. Sometimes 
only weaker properties of the algebra are known initially. An $H$-module algebra 
$A$ is called \emph{$H$-semiprime} if $A$ has no nonzero nilpotent $H$-stable 
ideals. Furthermore, $A$ is \emph{$H$-prime} if $A\ne0$, and $IJ\ne0$ for all 
nonzero $H$-stable ideals $I$ and $J$ of $A$.

As follows from \cite{Sk-Oy06, Lemma 4.2}, every right artinian $H$-prime 
algebra $A$ is necessarily $H$-simple. Several ring-theoretic constructions 
generally produce semiprimary rather than artinian rings, and the artinianness 
of subalgebras in Theorem 0.1 will have to be derived. This will be based on 
Proposition 2.13 whose proof makes use of Theorem 2.1 applied to semiprimary 
module algebras.

Recall that a semilocal ring is \emph{semiprimary} if its Jacobson radical is 
nilpotent. I don't know whether every semiprimary $H$-prime algebra is 
$H$-simple. Still this will be asserted in the situation of Corollary 2.7, and 
here we will use the equivalence of Proposition 2.2 which is valid, more 
generally, for semiperfect module algebras. Recall that \emph{semiperfect} 
rings are semilocal rings with the property that all their finitely-generated 
modules admit projective covers.

To derive desired conclusions we will use the twisting operations introduced in 
\cite{Sk-Oy06}. For a right $A$-module $V$ and a right $H$-comodule $U$ consider 
$U\ot V$ as a right $A$-module with respect to the twisted action of $A$ defined 
by the formula
$$
(u\ot v)\,a=\sum u\0\ot v\bigl(S(u\1)a\bigr),\qquad u\in U,\ v\in V,\ a\in A.\eqno(2.1)
$$
In this way the category $\MA$ of right $A$-modules becomes a left module 
category over the monoidal category $\MH$ of right $H$-comodules. By 
\cite{Sk-Oy06, Lemma 1.1} $\,U\ot V\,$ is projective whenever so is $V$. If 
$\,\dim_kU<\infty$ and $V$ is finitely generated, then $U\ot V$ is finitely 
generated too.

Suppose that $\,\dim_kU<\infty$. Then the twisting functor $U\ot\,?$ on 
$\MA$ induces an endomorphism of the \emph{Grothendieck group} 
$K_0(A)$ of the category of finitely generated projective right $A$-modules. 
This abelian group is universal with respect to the property that it 
contains elements $[P]$ associated with the isomorphism classes of finitely 
generated projective right $A$-modules which satisfy the relations 
$[P]=[P']+[P'']$ for each triple of finitely generated projectives 
$P,P',P''\in\MA$ such that $P\cong P'\oplus P''$.

We get also an endomorphism of the dual group 
$K_0(A)^*=\Hom_\bbZ(K_0(A),\,\bbZ)$. It sends $f\in K_0(A)^*$ to another 
function $fU\in K_0(A)^*$ defined by the formula
$$
(fU)(x)=f(U\ot x),\qquad x\in K_0(A).\eqno(2.2)
$$
Each series of subcomodules $0=U_0\sbs U_1\sbs\cdots\sbs U_n=U$ gives rise to 
a series of submodules
$$
0=U_0\ot P\sbs U_1\ot P\sbs\cdots\sbs U_n\ot P=U\ot P.
$$
If $P$ is projective in $\MA$, then so are all factors of the above series, 
whence
$$
U\ot P\cong\bigoplus_{i=1}^n\,\bigl((U_i/U_{i-1})\ot P\bigr).
$$
This shows that the isomorphism class of the right $A$-module $U\ot P$ depends 
only on the set of composition factors of $U$. In other words, $K_0(A)$ and 
$K_0(A)^*$ are, respectively, left and right modules over the 
\emph{Grothendieck ring} $\,G_0(\MH)\,$ of the monoidal category of 
finite-dimensional right $H$-comodules.

If $A$ is semiperfect, then a projective right $A$-module $P$ is indecomposable, 
if and only if $P$ is the projective cover of a simple module, if and only if 
$P$ is isomorphic to some right ideal of $A$ generated by a primitive idempotent. 
In this case $K_0(A)$ is a free abelian group of finite rank with a standard 
basis $\,\Indproj(A)\,$ consisting of the isomorphism classes of 
indecomposable projective right $A$-modules. In $K_0(A)^*$ our considerations 
will refer to the dual basis.

We adopt the terminology used in the theory of nonnegative matrices. A matrix 
is \emph{nonnegative} (respectively, \emph{positive}) if all its entries are 
nonnegative (respectively, positive) real numbers. In a similar way we treat 
elements of the abelian groups $K_0(A)$ and $K_0(A)^*$ for a semiperfect algebra 
$A$. An element $x\in K_0(A)$ is \emph{nonnegative} if $x=[P]$ for some finitely 
generated projective module, and $[P]$ is \emph{positive} if, moreover, $P$ 
is a generator in $\MA$, i.e., each indecomposable projective is a direct 
summand of $P$. An additive function $f:K_0(A)\to\bbZ$ is \emph{positive} if 
$f([P])>0$ for each nonzero finitely generated projective module.

Denote by $\,\Max A\,$ the set of all maximal ideals of $A$. Since $A$ is 
semilocal, this set is finite and the factor rings $A/\frm$ are simple 
artinian for all $\frm\in\Max A$.

\setitemsize(a)
\proclaim
Proposition 2.2. 
For a semiperfect left $H$-module algebra $A$ the following three conditions 
are equivalent:

\item(a)
for each $\frm\in\Max A$ the largest $H$-stable ideal $\frm_H$ of $A$ 
contained in $\frm$ is a maximal $H$-stable ideal of $A,$ and so the factor 
algebra $A/\frm_H$ is $H$-simple,

\item(b)
there is an equivalence relation $\sim$ on the standard basis $\,\Indproj(A)$ 
of the group $K_0(A)$ such that for two indecomposable projective right 
$A$-modules $P$ and $P'$ one has $[P]\sim[P']$ if and only if $P'$ is 
isomorphic to a direct summand of $\,U\ot P$ for some right $H$-comodule $U,$

\item(c)
for each finite-dimensional right $H$-comodule $U$ whose composition factors 
include a sufficiently large finite set of simple comodules there exists a 
positive function $f\in K_0(A)^*$ such that $\,fU=(\dim_kU)f$.

\endproclaim

\Proof.
Denote by $\Ga$ the set of all pairs 
$\,([P],[P'])\in\Indproj(A)\times\Indproj(A)\,$ 
such that $P'$ is isomorphic to a direct summand of $U\ot P$ for some comodule 
$U\in\MH$. The binary relation on the set $\Indproj(A)$ defined by $\Ga$ is 
reflexive since $k\ot P\cong P$ where $k$ is the trivial one-dimensional right 
$H$-comodule. This relation is also transitive. Indeed, given indecomposable 
projectives $P,P',P''\in\MA$ such that $P'$ is isomorphic to a direct summand 
of $U\ot P$ and $P''$ is isomorphic to a direct summand of $U'\ot P'$ for some 
$U,U'\in\MH$, the $A$-module $P''$ is isomorphic to a direct summand of 
$U'\ot U\ot P$. Thus (b) will hold when the relation defined by $\Ga$ is also 
symmetric, which means that $([P'],[P])\in\Ga$ whenever $([P],[P'])\in\Ga$.

Since each comodule is the union of all its finite-dimensional subcomodules, the 
inclusion $([P],[P'])\in\Ga$ for two indecomposable projectives $P,P'\in\MA$ 
holds if and only if $P'$ is isomorphic to a direct summand of $U\ot P$ for 
some finite-dimensional comodule $U\in\MH$. Considering further a composition 
series of $U$, we deduce that $P'$ is then isomorphic to a direct summand of 
$S\ot P$ for some simple comodule $S\in\MH$.

Thus for each pair $\al=([P],[P'])\in\Ga$ we can pick a simple comodule 
$S_\al\in\MH$ such that $P'$ is isomorphic to a direct summand of $S_\al\ot P$. 
Let now $U\in\MH$ be any finite-dimensional comodule such that for each 
$\al\in\Ga$ there is a composition factor of $U$ isomorphic to $S_\al$. Then 
for each pair of indecomposable projectives $P,P'\in\MA$ we have 
$([P],[P'])\in\Ga$ if and only if $P'$ is isomorphic to a direct summand of 
$U\ot P$. In other words, $\Ga$ is the set of edges of the directed graph 
associated with the nonnegative matrix $M(\ph)$ of the endomorphism 
$$
\ph:K_0(A)\to K_0(A)
$$
induced by the functor $U\ot\,?$.

We may take the set $\{\,S_\al\mid\al\in\Ga\}$ or any larger finite set to be 
the set of simple comodules considered in (c). Let $U$ and $\ph$ be as in the 
previous paragraph. For each subset $I\sbs\Indproj(A)$ denote by $\<I\>$ the 
free abelian subgroup of $K_0(A)$ which has the set $I$ as its basis. Then 
$\<I\>$ is stable under $\ph$ if and only if $\Ga(P)\sbs I$ for each $[P]\in I$ 
where
$$
\Ga(P)=\{[P']\in\Indproj(A)\mid([P],[P'])\in\Ga\}.
$$
Since $\Ga$ is transitive, $\Ga(P)$ is the smallest subset of $\Indproj(A)$ 
containing $[P]$ which generates a $\ph$-stable subgroup of $K_0(A)$.

\smallskip
Suppose that (b) holds. Then $\Ga(P)$ is precisely the equivalence class 
containing $[P]$. In this case $K_0(A)$ is the direct sum of its $\ph$-stable 
subgroups $\<c\>$ corresponding to various $\Ga$-equivalence classes 
$c\sbs\Indproj(A)$. The restriction of $\ph$ to such a subgroup $\<c\>$ has a 
positive matrix $M(\ph)_c$. Assuming that the basis elements are suitably 
ordered, the whole matrix $M(\ph)$ is block diagonal with diagonal blocks 
$M(\ph)_c$ for various $c$.

Extend $\ph$ by linearity to the real vector space $K_0(A)\ot_\bbZ\bbR$, and 
denote by $\ph^*$ the dual linear operator acting on the vector space 
$K_0(A)^*\ot_\bbZ\bbR\cong\Hom_\bbZ(K_0(A),\bbR)$. This space is a direct sum 
of its $\ph^*$-stable subspaces canonically isomorphic to the spaces 
$\<c\>^*_\bbR=\Hom_\bbZ(\<c\>,\bbR)$ for various $\Ga$-equivalence classes 
$c$.

In each $\<c\>^*_\bbR$ consider the basis dual to the basis of the group 
$\<c\>$ formed by the elements of $c$. The matrix of the restriction of 
$\ph^*$ to $\<c\>^*_\bbR$ in this basis is positive since it coincides with 
the transpose of $M(\ph)_c$. By the Perron-Frobenius theorem there exists an 
eigenvector $f_c\in\<c\>^*_\bbR$ of the linear operator $\ph^*$ which is 
positive in the sense that $f_c([P])>0$ for all $[P]\in c$, and moreover the 
one-dimensional subspace spanned by $f_c$ is the whole eigenspace of the 
restriction of $\ph^*$ to $\<c\>^*_\bbR$ corresponding to the largest real 
eigenvalue $\la$. Note that $f_c([P])=0$ for all $[P]\in\Indproj(A)$ such 
that $[P]\notin c$.

We can show that $\la=d$ where $d=\dim_kU$, and so this value does not depend 
on $c$. The class $[A]$ of the rank 1 free right $A$-modules is clearly a 
positive element of $K_0(A)$. In particular, $f_c([A])>0$. Since $A$ may be 
viewed as an object of $\HMA$, it follows from \cite{Sk-Oy06, Lemma 1.2(i)} 
that $U\ot A$ is a free $A$-module of rank $d$. Hence
$$
\ph([A])=U\ot[A]=d\,[A]\quad\hbox{in $K_0(A)$},\eqno(2.3)
$$
and therefore $\ph^*(f_c)([A])=(f_c\circ\ph)([A])=d\cdot f_c([A])$. Since 
$\ph^*(f_c)=\la f_c$, we get $\la=d$, as claimed.

Since $\la\in\bbZ$ and the matrix $M(\ph)_c$ is also integer-valued, the 
one-dimensional $\la$-eigenspace of the operator $\ph^*|_{\<c\>^*_\bbR}$ 
contains a vector with all coordinates in $\bbQ$. Any such a vector is a 
scalar multiple of $f_c$. Replacing $f_c$ with a suitable scalar multiple we 
may therefore assume that $f_c([P])\in\bbZ$ for all $[P]\in c$.

Selecting such a function $f_c$ for each $\Ga$-equivalence class $c$, we put 
$f=\sum f_c$, the sum over all equivalence classes. Then $f$ is a positive 
element of the group $K_0(A)^*$ and $fU=\ph^*(f)=d\mskip1mu f$. These are 
exactly the conditions required by (c). Thus (b)$\,\Rar\,$(c).

\smallskip
Conversely, suppose that (c) holds. If the set $\Indproj(A)$ is written as a 
disjoint union of two subsets $I$ and $I'$, then $K_0(A)=\<I\>\oplus\<I'\>$. 
We claim that the subgroup $\<I\>$ is stable under $\ph$ whenever so is 
$\<I'\>$. This is checked as follows. Let us write $[A]=x+x'$ and 
$\ph(x)=y+y'$ where $x,y\in\<I\>$ and $x',y'\in\<I'\>$. Then
$$
\ph([A])=y+y'+\ph(x').
$$
On the other hand, $\ph([A])=dx+dx'$ by (2.3). Now compare the components of 
the two expressions lying in $\<I\>$. Since $\ph(x')\in\<I'\>$ by the 
assumptions, we get $y=dx$. However, $f\circ\ph=d\mskip1mu f$ according to 
(c). Hence
$$
f(y)+f(y')=f\bigl(\ph(x)\bigr)=d\mskip1mu f(x)=f(y),
$$
and we deduce that $f(y')=0$. Since $x$ is a nonnegative element of $K_0(A)$, 
so too is $\ph(x)$, and therefore also $y$ and $y'$. Hence $f(y')=0$ implies 
$y'=0$ by positivity of $f$. Thus $\ph(x)\in\<I\>$. Moreover, since $x$ is a 
linear combination of elements of $I$ with positive coefficients, we get 
$\ph([P])\in\<I\>$ for all $[P]\in I$. This does yield the desired inclusion 
$\ph(\<I\>)\sbs\<I\>$.

Suppose now that $[P],[P']\in\Indproj(A)$ are elements such that 
$([P],[P'])\in\Ga$. Take one subset $I'=\Ga(P')$, and let $I$ be the 
complementary subset. Then $[P]\notin I$ since $\<I\>$ is stable under $\ph$, 
while $[P']\notin I$. Hence $[P]\in\Ga(P')$, i.e., $([P'],[P])\in\Ga$. This 
verifies (b).

\smallskip
Thus we have proved that (b)$\,\Lrar\,$(c). This conclusion is essentially a 
consequence of the Perron-Frobenius theorem and another theorem of Frobenius 
according to which a nonnegative matrix $M$ is permutation similar to a block 
diagonal matrix with irreducible diagonal blocks of equal spectral radii if 
and only if both $M$ and the transpose of $M$ admit positive eigenvectors (see 
\cite{Ber-Pl, Th. 3.14} and \cite{Gant, Ch. XIII, Th.~7}).

To approach the other equivalence (a)$\,\Lrar\,$(b) we will interpret the 
binary relation $\Ga$ in terms of the maximal ideals of $A$. Let $P,P'\in\MA$ 
be two indecomposable projectives. These two modules are projective covers in 
the category $\MA$ of their simple factor modules, say $V$ and $V'$. Since $A$ 
is semilocal, each simple $A$-module is determined up to isomorphism by a 
maximal ideal of $A$ annihilating the module. Denote by $\frm,\frm'\in\Max A$ 
the annihilators of $V$ and $V'$, respectively. In order that $P'$ be 
isomorphic to a direct summand of $U\ot P$, it is necessary and sufficient 
that there exist a nonzero $A$-linear map $U\ot P\to V'$. However, by 
\cite{Sk-Oy06, Lemma 1.1(i)},
$$
\Hom_A(U\ot P,\,V')\cong\Hom_A(P,N)\eqno(2.4)
$$
where $N=\Hom_k(U,V')$ is endowed with an $A$-module structure by the rule
$$
(\eta a)(u)=\sum\eta(u\0)(u\1a),\qquad\eta\in N,\ a\in A,\ u\in U.\eqno(2.5)
$$
It follows that $P'$ is isomorphic to a direct summand of $U\ot P$ if and only 
if $V$ is isomorphic to a subfactor of the $A$-module $N$. It is immediately 
clear from (2.5) that the largest $H$-stable ideal $\frm'_H$ contained in 
$\frm'$ annihilates $N$. If $\,([P],[P'])\in\Ga$, then $\frm'_H$ must 
annihilate $V$, which means that $\frm'_H\sbs\frm$.

Conversely, suppose that $\frm'_H\sbs\frm$. Note that the isomorphism in (2.4) 
is valid for an arbitrary right $H$-comodule $U$. So we may take $U=H$ with 
the comodule structure given by the comultiplication of $H$. Taking linear 
maps $\eta\in N$ defined by the formula $\eta(u)=\ep(u)w$ for various elements 
$w\in V'$, we infer from (2.5) that the annihilator of $N$ in $A$ coincides 
with $\frm'_H=\{a\in A\mid Ha\sbs\frm'\}$. There is a primitive idempotent 
$e\in A$ such that $P\cong eA$ in $\MA$. Since $\frm$ annihilates $V$, we 
deduce that $e\mskip1mu\frm$ is contained in the maximal $A$-submodule of $eA$. 
Hence $e\notin e\mskip1mu\frm$. Since the right ideal $eA$ is a direct summand 
of $A$, we deduce that $e\notin\frm$. Then $e\notin\frm'_H$ as well, and 
therefore $Ne\ne0$. This implies that $\Hom_A(eA,N)\ne0$, i.e., 
$\Hom_A(P,N)\ne0$, and by (2.4) we do get $\,([P],[P'])\in\Ga$.

We conclude that condition (b) holds if and only if for each pair of maximal 
ideals $\frm,\frm'\in\Max A$ the two inclusions $\frm'_H\sbs\frm$ and 
$\frm_H\sbs\frm'$ are equivalent to each other, and then each inclusion is 
also equivalent to the equality of $H$-stable ideals $\frm_H=\frm'_H$. Since 
each ideal of $A$ is contained in a maximal ideal, that property amounts to 
condition (a).
\endproof                                                                                                                                              

\Remark.
It is not clear whether for an algebra $A$ satisfying all conditions in the 
statement of Proposition 2.2 there always exists a positive function 
$\,f\in K_0(A)^*\,$ which has the eigenvector property requested by condition 
(c) and does not depend on the comodule $U$. Clearly this question is answered 
in the affirmative when the center of the Grothendieck ring $G_0(\MH)$ 
contains arbitrarily large elements. For example, this happens to be true when 
$H$ is finite-dimensional in which case the class of the right $H$-comodule 
$H$ lies in the center of $G_0(\MH)$. If $H$ is a coquasitriangular Hopf 
algebra then the whole ring $G_0(\MH)$ is commutative.

Another case occurs when the subalgebra of $H$-invariants $A^H$ has an ideal 
$J_0$ contained in the Jacobson radical of $A$ such that $A/AJ_0$ is a right 
$A^H$-module of finite length. Under this assumption the desired function $f$ 
is defined explicitly by setting $f([P])$ for each finitely generated 
projective $P\in\MA$ to be the length of the right $A^H$-module $P/PJ_0$.  
\endremark

\proclaim
Corollary 2.3. 
For a semiperfect left $H$-module algebra $A$ the following conditions are 
equivalent:

\item(a)
each proper $H$-stable ideal of $A$ is contained in the Jacobson radical of $A,$

\item(b)
there exists a right $H$-comodule $U$ such that $U\ot P$ is a generator in 
$\MA$ for each nonzero projective right $A$-module $P$.

\endproclaim

\Proof.
Let $J$ be the Jacobson radical of $A$ and $J_H$ the largest $H$-stable ideal 
of $A$ contained in $J$. Condition (a) means that $J_H$ is a unique maximal 
$H$-stable ideal of $A$. Since $J\sbs\frm$ for each $\frm\in\Max A$, this can 
be rephrased by saying that $A$ has a unique maximal $H$-stable ideal which 
coincides with the ideal $\frm_H$ for each $\frm\in\Max A$, and so, in 
particular, $\frm_H=\frm'_H$ for all $\frm,\frm'\in\Max A$.

Condition (b) means that there is $U\in\MH$ such that for each pair of 
indecomposable projectives $P,P'\in\MA$ the projective $A$-module $U\ot P$ 
contains a direct summand isomorphic to $P'$, i.e., the set $\Ga$ defined in 
the proof of Proposition 2.2 is the whole $\,\Indproj(A)\times\Indproj(A)$. We 
conclude by referring to the proof of the equivalence (a)$\,\Lrar\,$(b) in 
Proposition 2.2.
\endproof

If $A$ and $B$ are two left $H$-module algebras then the twisting functors 
$U\ot\,?$ are defined on both $\MA$ and $\MB$. The next lemma shows that these 
endofunctors of module categories commute with the functor $\MA\to\MB$ 
obtained by tensoring with an $H$-equivariant $(A,B)$-bimodule $N$. Recall the 
definition of the twisted action in formula (2.1).

\proclaim
Lemma 2.4.
Let $A$ and $B$ be two left $H$-module algebras. For $\,U\in\MH\!,$ 
$V\in\MA\mskip1mu,$ and $N\in\HAMB$ there is a well-defined bijective 
$B$-linear map 
$$
\displaylines{
\ph:\ U\ot(V\ot_AN)\to(U\ot V)\ot_AN\cr
\text{such that}\quad\ph\bigl(u\ot(v\ot_An)\bigr)=\sum\,(u\0\ot v)\ot_Au\1n
}
$$
for $u\in U,$ $v\in V,$ and $n\in N$. Thus 
$\ U\ot(V\ot_AN)\cong(U\ot V)\ot_AN\ $ in $\MB$.
\endproclaim

\Proof.
This map is well-defined since for all $a\in A$ and $u,v,n$ as above we have
$$
\eqalign{
\sum\,(u\0\ot v)\ot_Au\1(an)
&=\sum\,(u\0\ot v)\ot_A(u\1a)(u\2n)\cr
&=\sum\,(u\0\ot v)(u\1a)\ot_Au\2n\cr
&=\sum\,\bigl(u\0\ot v(S(u\1)u\2a)\bigr)\ot_Au\3n\cr
&=\sum\,(u\0\ot va)\ot_Au\1n
}
$$
in $(U\ot V)\ot_AN$. Furthermore, $\ph$ is $B$-linear since
$$
\eqalign{
\ph\bigl((u\ot(v\ot_An))\,b\bigr)
&=\ph\bigl(\ \sum u\0\ot\bigl(v\ot_An(S(u\1)b)\bigr)\,\bigr)\cr
&=\sum\,(u\0\ot v)\ot_Au\1\bigl(n(S(u\2)b)\bigr)\cr
&=\sum\,(u\0\ot v)\ot_A(u\1n)\bigl(u\2S(u\3)b\bigr)\cr
&=\sum\,(u\0\ot v)\ot_A(u\1n)\,b
}
$$
for all $b\in B$. The inverse map $\ph^{-1}$ is given by the assignment
$$
(u\ot v)\ot_An\mapsto\sum\,u\0\ot\bigl(v\ot_AS(u\1)n\bigr)
$$
It is induced by a linear map $(U\ot V)\ot N\to U\ot(V\ot_AN)$ which sends 
$(u\ot v)a\ot n$ and $(u\ot v)\ot an$ for $a\in A$ to the same element
$$
\eqalign{
\sum\,u\0\ot\bigl(v(S(u\2)a)\ot_AS(u\1)n\bigr)
&=\sum\,u\0\ot\bigl(v\ot_A(S(u\2)a)(S(u\1)n)\bigr)\cr
&=\sum\,u\0\ot\bigl(v\ot_AS(u\1)(an)\bigr)
}
$$
in $\,U\ot(V\ot_AN)$, and so $\ph^{-1}$ is well-defined.
\endproof

As a special case of Lemma 2.4 the twisting endofunctors of module categories 
commute with the extension of scalars functors:

\proclaim
Corollary 2.5.
Given a homomorphism of left $H$-module algebras $A\to B,$ there are natural 
isomorphisms of right $B$-modules
$$
U\ot(V\ot_AB)\cong(U\ot V)\ot_AB
$$
for comodules $U\in\MH$ and modules $V\in\MA$.
\endproclaim

\Proof.
We apply Lemma 2.4, viewing $B$ as an object of $\HAMB$.
\endproof

\proclaim
Corollary 2.6.
Let $B$ be a semiperfect $H$-simple left $H$-module algebra. Then any 
semiprimary $H$-semiprime $H$-module subalgebra $A$ of $B$ is isomorphic to a 
direct product of finitely many $H$-simple $H$-module algebras.
\endproclaim

\Proof.
The $H$-simple algebra $B$ satisfies condition (a) of Proposition 2.2. Hence 
for each sufficiently large in the sense of (c) finite-dimensional right 
$H$-comodule $U$ there exists a positive function $g\in K_0(B)^*$ such that 
$\,gU=(\dim_kU)g$.

The extension of scalars functor $?\ot_AB:\MA\to\MB$ induces a homomorphism of 
abelian groups $K_0(A)\to K_0(B)$. If $e$ is any nonzero idempotent of $A$, 
then the $B$-module $eA\ot_AB\cong eB$ is nonzero. This means that nonzero 
nonnegative elements of $K_0(A)$ are taken to nonzero nonnegative elements of 
$K_0(B)$, and therefore the function $f\in K_0(A)^*$ defined by the rule
$$
f(x)=g(x\ot_A\!B),\qquad x\in K_0(A),
$$
is positive. It follows from Corollary 2.5 that $\,fU=(\dim_kU)f$. This shows 
that $A$ satisfies condition (c) of Proposition 2.2. Hence $\frm_H$ is a 
maximal $H$-stable ideal of $A$ for each $\frm\in\Max A$. Let $I_1,\ldots,I_r$ 
be all distinct maximal $H$-stable ideals of $A$. Their intersection is an 
$H$-stable ideal contained in all maximal ideals of $A$, and therefore in the 
Jacobson radical $J$ of $A$. Since $J$ is nilpotent, we get 
$I_1\cap\cdots\cap I_r=0$ by $H$-semiprimeness of $A$. Hence
$$
A\cong A/I_1\times\cdots\times A/I_r
$$
by the Chinese remainder theorem.
\endproof

\proclaim
Corollary 2.7.
Let $A$ be an $H$-module subalgebra of a semiperfect left $H$-module algebra $B$. 
If $A$ is semiprimary and $B$ contains no nonzero proper $H$-stable left ideals, 
then $A$ is $H$-simple.
\endproclaim

\Proof.
If $I$ is any nonzero $H$-stable ideal of $A$, then $BI$ is a nonzero $H$-stable 
left ideal of $B$, whence $BI=B$ by the hypothesis. It follows that $BIJ=B$, and 
therefore $IJ\ne0$ for any pair of nonzero $H$-stable ideals $I,J$ of $A$. Thus 
$A$ is $H$-prime. Applying Corollary 2.6, we conclude that $A$ must be $H$-simple.
\endproof

In Corollary 2.7 the algebra $B$ may be viewed as an object of $\HMA$, but it 
is not necessarily locally $A$-finite, and therefore we cannot say that $B$ is 
projective in $\MA$. Under suitable conditions projectivity and even freeness 
of an $H$-module algebra over an $H$-module subalgebra will be proved in 
Proposition 2.12. First we recall a well-known fact in the next lemma.

\proclaim
Lemma 2.8.
Let $M,N\in\HMA$. Then $\Hom_A(M,N)$ is a left $H$-module with respect to the 
action of $H$ defined by the formula
$$
(hf)(x)=\sum h\1f(S(h\2)x),\qquad h\in H,\ f\in\Hom_A(M,N),\ x\in M.
$$
\endproclaim

\proclaim
Corollary 2.9.
Let $M\in\HMA$. If $A$ is $H$-simple and $\Hom_A(M,A)\ne0,$ then $M$ is a 
generator in $\MA$.
\endproclaim

\Proof.
The evaluation map $\Hom_A(M,A)\ot M\to A$ is a morphism in the monoidal 
category of left $H$-modules. Hence its image $I$ is a nonzero $H$-stable 
ideal $A$. It follows that $I=A$ by the $H$-simplicity of $A$. In other words, 
the sum of images of all $A$-module homomorphisms $M\to A$ coincides 
with the whole $A$.
\endproof

The condition $\Hom_A(M,A)\ne0$ is satisfied when $M$ is a nonzero projective 
in $\MA$, and then the conclusion of Corollary 2.9 would tell us that $M$ is a 
projective generator in $\MA$. Under further restrictions we can prove even 
freeness of $M$ (this is not absolutely necessary for the main results of the 
present paper).

\proclaim
Lemma 2.10.
Suppose that $A$ is a semiperfect $H$-simple left $H$-module algebra. If 
$M\in\HMA$ is projective in $\MA$ and not $A$-finite, then $M$ is a free right 
$A$-module.
\endproclaim

\Proof.
Let $P_1,\ldots,P_r$ be a full set of pairwise nonisomorphic indecomposable 
projectives in $\MA$. Each projective right module over a semiperfect ring 
is a direct sum of indecomposable projectives 
(see \cite{Mue70, Th.~3} or \cite{And-F, Th. 27.11}). Hence
$$
M\cong P_1^{(X_1)}\oplus\cdots\oplus P_r^{(X_r)}\quad\hbox{in $\MA$}
$$
for some sets $X_1,\ldots,X_r$ where $P_i^{(X_i)}$ is the direct sum of the 
family of copies of $P_i$ indexed by the set $X_i$. Let $\ka$ be the largest 
among the cardinalities of these sets $X_1,\ldots,X_r$. Since $M$ is not 
$A$-finite, $\ka$ is an infinite cardinal.

By Corollary 2.3 there exists a right $H$-comodule $U$ such that $U\ot P_i$, 
regarded as a right $A$-module with respect to the twisted action of $A$, is a 
projective generator in $\MA$ for each $i$. Replacing $U$ with a suitable 
subcomodule, we may assume that $U$ is  finite-dimensional. Then $U\ot P_i$ is 
a finitely generated $A$-module, and therefore each indecomposable projective 
occurs as a direct summand of $U\ot P_i$ with a finite nonzero multiplicity. 
Now
$$
U\ot M\cong(U\ot P_1)^{(X_1)}\oplus\cdots\oplus(U\ot P_r)^{(X_r)}.
$$
Expressing each summand here as a direct sum of indecomposable projectives 
and counting multiplicities, we deduce that 
$\,U\ot M\cong P_1^{(\ka)}\oplus\cdots\oplus P_r^{(\ka)}\cong A^{(\ka)}\,$ is a 
free right $A$-module of rank $\ka$. On the other hand, $\,U\ot M\cong M^d\,$ 
where $d=\dim_kU$ by \cite{Sk-Oy06, Lemma 1.2(i)}. Hence $M\cong A^{(\ka)}$ as 
well.
\endproof

\proclaim
Lemma 2.11.
Suppose that $A$ is a semiperfect $H$-simple left $H$-module algebra and $M$ 
is an $A$-finite object of the category $\HMA$. If there exists a ring 
homomorphism $A\to\End_AM,$ then $M$ is a free right $A$-module.
\endproclaim

\Proof.
It follows from Theorem 2.1 that $M$ is a projective right $A$-module with the 
property that a suitable direct sum of several copies of $M$ is a free 
$A$-module. There is a smallest projective $E\in\MA$ having this property. 
In terms of indecomposable projectives $P_1,\ldots,P_r\in\MA$ we have a 
decomposition
$$
E\cong P_1^{a_1}\oplus\cdots\oplus P_r^{a_r}
$$
for some integers $a_1,\ldots,a_r>0$ such that $\gcd(a_1,\ldots,a_r)=1$ and
$$
E^d\cong P_1^{a_1d}\oplus\cdots\oplus P_r^{a_rd}\cong A\quad\hbox{in $\MA$}
$$
for some integer $d>0$. The ring homomorphism $A\to\End_AM$ allows us to view 
$M$ as an $A$-bimodule. Then
$$
M\cong A\ot_AM\cong E^d\ot_AM\cong L^d\quad\hbox{in $\MA$}
$$
where $L=E\ot_AM$. So $L$ is a projective right $A$-module with the property 
that a suitable direct sum of several copies of $L$ is a free $A$-module. 
Considering a decomposition of $L$ as a direct sum of indecomposables, we 
deduce that $L\cong E^n$ for some integer $n>0$. Hence $M\cong E^{nd}\cong A^n$ 
is a free right $A$-module.
\endproof

In the next proposition we have to assume that $A$ is not just semiperfect, 
but even right perfect. A ring is called \emph{right perfect} if each its 
right module admits a projective cover. Several characterizations of perfect 
rings were given by Bass \cite{Bas60}. Each right or left perfect ring is 
semiperfect, while semiprimary rings are right and left perfect.

A multiplicatively closed subset $T$ of a ring $R$ is a \emph{left Ore set} if 
for each $a\in R$ and each $s\in T$ there exists $t\in T$ such that $ta\in Rs$. 
Suppose that $T$ is a left Ore set consisting of nonzerodivisors of $R$. 
An overring $Q$ of a ring $R$ is called a \emph{left quotient ring} of $R$ with 
respect to $T$ if all elements of $T$ are invertible in $Q$ and 
$$
Q=\{s^{-1}a\mid a\in R,\ s\in T\}.
$$
Such a ring exists and is unique up to isomorphism.

In the case when $T$ is a left Ore set consisting of all nonzerodivisors of 
$R$ such an overring $Q$ of $R$ is a \emph{classical left quotient ring} of 
$R$. Right quotient rings are defined similarly.

\proclaim
Proposition 2.12.
Let $A$ and $B$ be two $H$-module subalgebras of a left $H$-module algebra $Q$ 
such that $A$ is right perfect and $H$-simple, while $B$ is the sum of its 
finite-dimensional $H$-submodules. Suppose that $Q$ is a left quotient ring 
of $B$ with respect to a left Ore set $T$ of nonzerodivisors of $B$.

Then $Q$ is a free $A$-module with respect to the action of $A$ by right 
multiplications.
\endproclaim

\Proof.
We may view $Q$ as an object of $\HAMA$. If $V\sbs Q$ is any $H$-submodule 
of finite dimension, then $VA$ is an $A$-finite $\HMA$-subobject. Hence $M=BA$ 
is a locally $A$-finite $\HMA$-subobject of $Q$. By Theorem 2.1 $M$ is 
projective in $\MA$. Then so too is $s^{-1}M\cong M$ for each $s\in T$. Since 
$T$ satisfies the left Ore condition, the set of right $A$-submodules 
$\{s^{-1}M\mid s\in T\}$ is directed by inclusion. Since $B$ is contained 
in $M$, it follows that
$$
\textstyle Q=\bigcup\limits_{s\in T}s^{-1}B
=\bigcup\limits_{s\in T}s^{-1}M
$$
is a directed colimit of projective right $A$-modules, whence $Q$ is right 
$A$-flat. Since all flat right modules over any right perfect ring are 
projective \cite{Bas60}, $Q$ is projective in $\MA$. Finally, freeness follows 
from either Lemma 2.10 or Lemma 2.11.
\endproof

The notions of dominions and dominion subrings have been already recalled in 
section 1. As an easy example, any subring $A$ of a ring $R$ which is a 
direct summand of $R$ as an $A$-module with respect to either left or right 
multiplications is a dominion subring of $R$.

\proclaim
Proposition 2.13.
Let $Q$ be a left $H$-module algebra, $B$ its $H$-module subalgebra. 
Suppose that

\item(1)
\ $B$ is the sum of its finite-dimensional $H$-submodules,

\item(2)
\ each nonzero $H$-stable left ideal of $B$ contains a nonzerodivisor of $B,$

\item(3)
\ $Q$ is a left artinian classical left quotient ring of $B$.

\noindent
Then for each $H$-module subalgebra $A\sbs Q$ the following conditions 
are equivalent:

\item(a)
$\,A$ is a dominion subalgebra of $Q,$

\item(b)
$\,A$ is semiprimary,

\item(c)
$\,A$ is left artinian.

Such a subalgebra $A$ has no nonzero proper $H$-stable left ideals.
\endproclaim

\Proof.
By a theorem of Schofield dominions of subrings in any semiprimary ring are 
semiprimary \cite{Scho, Th. 7.19} (an alternative proof is given in the 
Appendix). Hence (a)$\,\Rar\,$(b). The implication (c)$\,\Rar\,$(b) is a 
standard fact of ring theory.

Suppose that (b) holds. If $I$ is any nonzero $H$-stable left ideal of $Q$, 
then $I\cap B$ is a nonzero $H$-stable left ideal of $B$, whence $I$ contains 
an invertible element of $Q$ by conditions (2) and (3), i.e., $I=Q$. Thus the 
hypothesis of Corollary 2.7 is satisfied for the pair $A$ and $Q$, and we 
deduce that $A$ is $H$-simple. 

By Proposition 2.12 $Q$ is free as a right $A$-module, and so a generator 
in $\MA\mskip1mu$. This implies that $A$ is an $\MA$-direct summand of $Q$ 
\cite{Fai, 3.27}, and (a) follows. So does condition (c) too since the 
assignment $L\mapsto QL$ embeds the lattice of left ideals of $A$ into that of 
left ideals of $Q$. In fact $A\cap QL=L$ for each left ideal $L$ of $A$. If 
$L$ is stable under the action of $H$, then so is the left ideal $QL$ of $Q$, 
and we deduce that $L=0$ when $QL=0$, or $L=A$ when $QL=Q$.
\endproof

\Remark.
We note that in Proposition 2.13 the condition that \emph{$A$ is right 
perfect}, i.e., a weaker form of (b), can be added to the list of equivalent 
conditions (a), (b), (c). This follows from the fact that the conclusion of 
Corollary 2.7 still holds under this weaker assumption about $A$. If $A$ is 
right perfect, then its Jacobson radical $J$ is right $T$-nilpotent, whence 
$VJ\ne V$ for each nonzero right $A$-module. In Corollary 2.7 we get $BJ\ne 
B$, and this implies that $J$ cannot contain any nonzero $H$-stable ideals of 
$A$. With this observation the arguments in the proofs of Corollaries 2.6 and 
2.7 are still valid.
\endremark

At the end of this section we mention yet another freeness result.

\proclaim
Lemma 2.14.
Let $A$ be a semilocal $H$-simple left $H$-module algebra. Suppose that a ring 
homomorphism $A\to R$ is given where $R$ is a ring with the IBN (invariant 
basis number) property. If $M$ is a locally $A$-finite object of the category 
$\HMA$ such that $M\ot_AR$ is a free $R$-module, then $M$ is a free 
$A$-module.
\endproclaim

\Proof.
If $M$ is not $A$-finite, then its freeness in $\MA$ follows already from 
Theorem 2.1. So assume that $M$ is $A$-finite. Theorem 2.1 still shows that a 
direct sum of several copies of $M$ is $A$-free. Hence $M^n\cong A^m$ in $\MA$ 
for some integers $n,m\ge0$. Then $(M\ot_AR)^n\cong R^m$ in $\calM_R$. Since 
$M\ot_A\!R$ is a free $R$-module, the IBN implies that $n$ divides $m$. 

Let $F$ be the free right $A$-module of rank $m/n$, and let $J$ be the 
Jacobson radical of $A$. Since $M^n\cong F^n$, we also have 
$(M/MJ)^n\cong(F/FJ)^n$. It follows that each simple right $A$-module occurs 
as a direct summand in $M/MJ$ and in $F/FJ$ with equal multiplicities. Hence 
$M/MJ\cong F/FJ$ since the factor ring $A/J$ is semisimple artinian. Finally, 
$M\cong F$ in $\MA$ by projectivity of $M$.
\endproof

\section
3. First properties of the quotient ring and invariant subalgebras

We assume that $H$ satisfies basic assumptions (A1) and (A2). Assumption (A1) 
means that the finite dual $\Hd$ of $H$ is a dense subalgebra of the dual 
algebra $H^*$ (in the terminology of Sweedler \cite{Sw} residually 
finite-dimensional algebras are called \emph{proper}).

Recall that right $H$-comodules may be identified with rational left 
$H^*$-modules (see \cite{Mo} and \cite{Sw}). Given a comodule structure map 
$\rho:V\to V\ot H$ for a vector space $V$ each element $f\in H^*$ acts on $V$ 
by means of the linear operator obtained as the composite
$$
V\lmapr2\rho V\ot H\lmapr4{\id\ot f}V\ot k\cong V,
$$
i.e., $\,f\rhu v=\sum f(v\1)v\0$ for $v\in V$.

Each rational $H^*$-module is completely determined by its restriction to any 
dense subalgebra. Thus the category $\MH$ of right $H$-comodules may be 
identified with the category of rational left $\Hd$-modules. An arbitrary left 
$\Hd$-module $V$ contains a largest rational submodule which will be denoted 
$\,\Rat(V)$.

As recalled in section 1, the quotient ring $Q=Q(H)$ is an $\Hd$-bimodule 
algebra with respect to the two actions of $\Hd$ on $Q$ extending the well-known 
actions on $H$. In particular, $Q$ is a left $\Hd$-module algebra.

So we have the category of $\Hd$-equivariant right $Q$-modules $\HdMQ$. 
Its object $M$ will be called \emph{rationally generated} if $\,M=\Rat(M)Q$. 
We denote by $\HdratMQ$ the full subcategory of rationally generated objects 
of $\HdMQ$. 

\proclaim
Lemma 3.1.
The largest rational $\Hd$-submodule\/ $\Rat(M)$ of an equivariant module 
$M\in\HdMQ$ is also an $H$-submodule. The right $H$-comodule structure 
corresponding to the rational $\Hd$-module structure makes\/ $\Rat(M)$ an 
object of the category $\MHH$ of right Hopf modules.
\endproclaim

\Proof.
The $Q$-module structure map $M\ot Q\to M$ is $\Hd$-linear. Therefore the 
image of the rational $\Hd$-submodule $\Rat(M)\ot H\sbs M\ot Q$ is a rational 
$\Hd$-submodule of $M$, i.e., it is contained in $\Rat(M)$. By restriction we 
obtain a morphism of rational $\Hd$-modules $\Rat(M)\ot H\to\Rat(M)$ which 
defines an $H$-module structure on $\Rat(M)$ and which may be viewed as a 
morphism in the category $\MH$. This is precisely the compatibility 
condition required for objects of $\MHH$.
\endproof

The subspace of $\Hd$-invariants of a left $\Hd$-module $V$ is
$$
V^{\Hd}\!=\{v\in V\mid fv=f(1)v\,\hbox{ for all }f\in\Hd\}.
$$
Note that $\,V^{\Hd}\!\sbs\Rat(V)$. If $\,V$ is a rational left $\Hd$-module 
with the corresponding comodule structure $\,\rho:V\to V\ot H$, \ then
$$
V^{\Hd}\!=V\co H=\{v\in V\mid\rho(v)=v\ot1\}.
$$
For an arbitrary $\Hd$-module $V$ we have 
$\,V^{\Hd}\!=\Rat(V)^{\Hd}\!=\Rat(V)\co H$.

\proclaim
Proposition 3.2.
The functor of taking $\Hd$-invariants $\,M\mapsto M^{\Hd}\!$ gives an 
equivalence between $\HdratMQ$ and the category $\Mk$ of vector spaces.

The quasi-inverse functor $\Mk\to\HdratMQ$ is given by $V\mapsto V\ot Q$ with 
the two module structures on $V\ot Q$ arising from those on $Q$.
\endproclaim

\Proof.
Let $M\in\HdratMQ$. Since $\Rat(M)$ is stable under the right action of 
$H\sbs Q$ by Lemma 3.1, there is a canonical map $\Rat(M)\ot_HQ\to M$ which is 
injective by the properties of classical quotient rings, and also surjective 
since $M$ is rationally generated. Furthermore, $\Rat(M)\in\MHH$ by Lemma 3.1. 
The structural description of Hopf modules \cite{Sw, Th. 4.1.1} yields a 
canonical isomorphism
$$ 
\Rat(M)\cong\Rat(M)\co H\ot H=M^{\Hd}\!\ot H,
$$
whence
$$
M\cong\Rat(M)\ot_HQ\cong(M^{\Hd}\!\ot H)\ot_HQ\cong M^{\Hd}\!\ot Q.
$$
From this we also see that the equality $M=WQ$ for a vector subspace 
$W\sbs M^{\Hd}\!$ holds only when $W=M^{\Hd}\!$.

Applying the last observation to the $\Hd$-equivariant $Q$-module $M=V\ot Q$ 
where $V$ is an arbitrary vector space, we deduce that 
$\,(V\ot Q)^{\Hd}\!\cong V$.
\endproof

The full subcategory $\HdratQM\sbs\HdQM$ of rationally generated 
$\Hd$-equivariant left $Q$-modules is defined similarly. If $H$ is changed to 
$H\op$, then $Q$ is changed to $Q\op$ and $\Hd$ is changed to $(\Hd)\cop$. 
Since the category $(\Hd)\cop\hbox{-}{\mskip1mu}\calM_{Q\op}$ may be 
identified with $\HdQM$, the $H\op$-variant of Proposition 3.2 yields

\proclaim
Proposition 3.3.
The functor $\,\HdratQM\to\Mk{\mskip1mu},$ $\,M\mapsto M^{\Hd}\!,$ is 
an equivalence of categories.
\endproclaim

\proclaim
Corollary 3.4.
Up to isomorphism $Q$ is the only simple object of $\HdratMQ$ and the only 
simple object of $\HdratQM$. Thus $Q$ has no $\Hd$-invariant one-sided ideals 
other than the zero ideal and the whole $Q$.
Furthermore, $\,\Rat(Q)=H$ and $\,Q^{\Hd}\!=k$.  
\endproclaim

Moreover, each finite-dimensional $\Hd$-submodule of $Q$ is contained in $H$, 
and therefore is rational. This conclusion is a special case of Lemma 4.1.

\proclaim
Corollary 3.5.
Given an exact sequence\/ $0\to K\to M\to N\to0$ in $\HdMQ,$ the objects $K$ 
and $N$ are rationally generated provided that so is $M$.
\endproclaim

\Proof.
By Proposition 3.2 we have $M\cong M^{\Hd}\!\ot Q$. The image of $M^{\Hd}$ in 
$N$ is clearly contained in $N^{\Hd}$ and generates $N$. Hence $N$ is 
rationally generated, and therefore $N\cong N^{\Hd}\!\ot Q$. The morphism 
$M\to N$ may be identified with the map
$$
\ph\ot\id:\,M^{\Hd}\!\ot Q\to N^{\Hd}\!\ot Q
$$
for some $k$-linear map $\,\ph:M^{\Hd}\!\to N^{\Hd}\!$. Its kernel $K$ is the 
$Q$-module generated by its subspace $\,\Ker\ph\sbs K^{\Hd}\!$. Hence $K$ is 
rationally generated too.
\endproof

\proclaim
Proposition 3.6.
For a left $\Hd$-invariant subalgebra $A\sbs Q$ the following three conditions 
are equivalent:

\item(a)
$\,A$ is a dominion subalgebra of $Q$,

\item(b)
$\,A$ is semiprimary,

\item(c)
$\,A$ is left and right artinian.

\noindent
Such a subalgebra $A$ has no $\Hd$-invariant one-sided ideals other than the 
zero ideal and the whole $A$. The ring $Q$ is a free $A$-module on both sides.
\endproclaim

\Proof.
We meet all assumptions of Proposition 2.13 with $B=H$ regarded as a left 
$\Hd$-module subalgebra of $Q$. Condition (1) in the hypothesis of that 
proposition is satisfied since $H$ is a rational left $\Hd$-module. A left 
ideal $I$ of $H$ is left $\Hd$-invariant if and only if $I$ is a right 
coideal, and so an $\HMH$-subobject of $H$. Assuming that $I\ne0$, we get 
$I=H$, i.e., $1\in I$, since $H$ is a simple object of the category $\HMH$ by 
the fundamental theorem on Hopf modules applied to the Hopf algebra $H\op$. 
This verifies condition (2), and (3) is one of basic assumptions about $H$.

As a special case of Proposition 2.13 we thus obtain (c)$\,\Rar\,$(a)$\,\Rar\,$(b), 
and (b) implies that $A$ is left artinian. The $H\op$-variant of this 
conclusion shows that $A$ is also right artinian. Absence of nontrivial 
$\Hd$-invariant one-sided ideals is also stated in Proposition 2.13, while 
freeness of $Q$ over $A$ follows from Proposition 2.12 and its $H\op$-variant.
\endproof

\section
4. Correspondence between corings and coalgebras

We keep our basic assumptions about $H$. Since $Q$ is a left $\Hd$-module 
algebra, there is the category of $\Hd$-equivariant $Q$-bimodules $\HdQMQ$. 

An object $M\in\HdQMQ$ will be called \emph{rationally generated} if $M$ 
coincides with its $Q$-subbimodule generated by the largest rational 
$\Hd$-submodule $\Rat(M)$. Exactly as in the proof of Lemma 3.1 we deduce that 
$\Rat(M)$ is stable under the left and right actions of $H$, i.e., $\Rat(M)$ 
is an $H$-subbimodule of $M$.

If $M$ is rationally generated, then for each $x\in M$ there exist 
nonzerodivisors $s,t\in H$ such that $sxt\in\Rat(M)$.  Since nonzerodivisors 
of $H$ are invertible in $Q$, it follows that each $Q$-subbimodule of $M$ is 
generated by its intersection with $\Rat(M)$. Hence each $\HdQMQ$-subobject of 
$M$ is rationally generated, and so too is each factor object. In particular, 
the full subcategory
$$
\HdratQMQ\sbs\HdQMQ
$$
of rationally generated objects is abelian, and the inclusion functor is exact.

\proclaim
Lemma 4.1.
Let $M\in\HdQMQ\,,$ and let $N$ be any $H$-subbimodule of $M$ which is stable 
also under the action of $\Hd\!$ and satisfies $\,QNQ=M$. If $\,V\!$ is any 
finite-dimensional $\Hd$-submodule of $M$ then $V\sbs N$. Therefore 
$\,M^{\Hd}\!\sbs\Rat(M)\sbs N$.
\endproclaim

\Proof.
The set $I=\{a\in H\mid aV\sbs NQ\}$ is a nonzero left $\Hd$-invariant left 
ideal of $H$. Hence $I$ is also a right coideal. Since $H$ is a simple object 
of $\HMH\mskip1mu$, we get $I=H$, i.e., $1\in I$. Hence $V\sbs NQ$. Now the 
set $J=\{a\in H\mid Va\sbs N\}$ is a nonzero left $\Hd$-invariant right ideal 
of $H$. It is a right coideal as well. Since $H$ is a simple object of 
$\MHH\mskip1mu$, we deduce that $1\in J$, and so $V\sbs N$.

The final conclusion follows from the fact that $\Rat(M)$ is the sum of its 
finite-dimensional $\Hd$-submodules.
\endproof

\proclaim
Proposition 4.2.
For $M\in\HdQMQ$ the subspace of $\Hd$-invariants $\,M^{\Hd}\!$ is a left 
$H$-module with respect to the action of $H$ defined by the formula
$$
h\triangleright x=\sum h\1xS(h\2),\qquad h\in H,\ x\in M^{\Hd}\!.\eqno(4.1)
$$
The functor $\,M\mapsto M^{\Hd}\!$ gives an equivalence between $\,\HdratQMQ\,$ 
and the category $\HM$ of left $H$-modules.
\endproclaim

\Proof.
Since the action of $\Hd\!$ on $\Rat(M)$ arises from a right $H$-comodule 
structure
$$
\rho:\Rat(M)\to\Rat(M)\ot H,
$$
we may view $\Rat(M)$ as an object of the category $\HMHH$ of two-sided Hopf 
modules. Since $M^{\Hd}\!\!\sbs\Rat(M)$, we have $M^{\Hd}\!\!=\Rat(M)\co H$. 
If $y=h\triangleright x$ where $x\in M^{\Hd}\!$ and $h\in H$, then
$$
\rho(y)=\sum h\1xS(h\4)\ot h\2S(h\3)=h\1xS(h\2)\ot1=y\ot1
$$
since $\rho(x)=x\ot1$, and so $y\in M^{\Hd}\!$. Also, $1\triangleright x=x$ and 
$(ab)\triangleright x=a\triangleright(b\triangleright x)$ for all $a,b\in H$. 
Thus, $M^{\Hd}\!$ is indeed a left $H$-module.

If $M\in\HdratQMQ$, then
$$
M\cong Q\ot_H\Rat(M)\ot_HQ.\eqno(4.2)
$$
Conversely, given any object $N\in\HMHH$, the $Q$-bimodule $M=Q\ot_HN\ot_HQ$ is 
an object of $\HdratQMQ$ with respect to the $\Hd$-module structure which makes 
$M$ a factor module of the tensor product $Q\ot N\ot Q$ of three left 
$\Hd$-modules. Since $N$ is a free $H$-module, and therefore torsionfree, on 
both sides, it embeds in $M$ by the map sending $x\in N$ to 
$1\ot_Hx\ot_H1\in M$, and we get $\Rat(M)\cong N$ by Lemma 4.1. This 
shows that the functor $M\mapsto\Rat(M)$ gives an equivalence of categories
$$
\HdratQMQ\approx\HMHH\,.\eqno(4.3)
$$

Given $N\in\HMHH$, we have
$$
hx=\sum h\1xS(h\2)h\3=\sum(h\1\triangleright x)h\2\quad
\text{for $x\in N$ and $h\in H$}.\eqno(4.4)
$$
The Hopf module $N$ is freely generated by its subspace 
$\,V=N\co H\!=N^{\Hd}\!$ as a right $H$-module. Hence $N\cong V\ot H$ with 
the corresponding $H$-bimodule structure on $V\ot H$ defined by the formula
$$
h\cdot(x\ot y)\cdot h'=\sum\,(h\1\triangleright x)\ot h\2yh',\qquad 
x\in V,\quad h,y,h'\in H,\eqno(4.5)
$$
and the right $H$-comodule structure provided by the map $\id\ot\De$. Conversely, 
for any left $H$-module $V$ the bimodule and comodule structures just described 
make $V\ot H$ an object of $\HMHH$, and we have $\,(V\ot H)\co H\cong V$.

Thus the functor $N\mapsto N\co H$ gives an equivalence of categories 
$\HMHH\approx\HM$. This is one of equivalences described by Schauenburg 
\cite{Scha94, Th. 5.7} in the general context of Hopf algebras and Hopf 
modules in a symmetric monoidal category.

The composite of the two functors $\,\HdratQMQ\to\HMHH\to\HM\,$ considered 
above is the equivalence asserted in the statement of Proposition 4.2.
\endproof

\proclaim
Lemma 4.3.
Each object $M\in\HdratQMQ$ is a free $Q$-module on the right and on the left.
\endproclaim

\Proof.
Put $N=\Rat(M)Q$. Then $N$ is an object of $\HdratMQ$. By Proposition 3.2 $N$ 
is a free right $Q$-module. Now $M=QN$ coincides with the union of the 
directed family
$$
\{s^{-1}N\mid\hbox{$s$ is a nonzerodivisor of $H$}\}
$$
of free right $Q$-submodules. Hence $M$ is a flat right $Q$-module. Since $Q$ 
is artinian, this module is projective. We deduce that it is even free applying 
either Lemma 2.10 or Lemma 2.11.

Freeness on the left is obtained by symmetry. Formally speaking we apply the 
already proved conclusion with $\,H,\,Q,\,\Hd$ replaced by 
$\,H\op,\,Q\op,\,(\Hd)\cop$.
\endproof

The tensor product $\ot_Q$ makes $\HdQMQ$ into a monoidal category and $\HdMQ$ 
into a right module category over $\HdQMQ$. However, $M\ot_QN$ is not rationally 
generated when $M$ and $N$ are rationally generated objects of respective 
categories. Nevertheless, the rational submodules behave nicely under tensoring.
Note that the category $\HMHH$ is monoidal with respect to the functor $\ot_H$. 

\proclaim
Lemma 4.4.
For $M\in\HdratMQ$ and $N\in\HdQMQ$ there are canonical natural isomorphisms
$$
\Rat(M\ot_QN)\cong\Rat(M)\ot_H\Rat(N),\qquad
(M\ot_QN)^{\Hd}\!\cong M^{\Hd}\!\!\ot N^{\Hd}
$$
in $\MHH$ and $\Mk,$ respectively.
\endproclaim

\Proof.
By Proposition 3.2 $\,M\cong M^{\Hd}\!\ot Q$, whence 
$M\ot_QN\cong M^{\Hd}\!\ot N$. Since $\Hd$ acts trivially on $M^{\Hd}\!$, we get
$\,(M^{\Hd}\!\ot N)^{\Hd}\cong M^{\Hd}\!\ot N^{\Hd}$ and
$$
\Rat(M^{\Hd}\!\ot N)\cong M^{\Hd}\!\ot\Rat(N)
\cong(M^{\Hd}\!\ot H)\ot_H\Rat(N)\cong\Rat(M)\ot_H\Rat(N),
$$
as claimed.
\endproof

\proclaim
Lemma 4.5.
For $M\in\HdratQMQ$ and $N\in\HdQMQ$ the canonical maps
$$
\Rat(M)\ot_H\Rat(N)\to\Rat(M\ot_QN),\qquad
M^{\Hd}\!\ot N^{\Hd}\to(M\ot_QN)^{\Hd}
$$
are isomorphisms in $\HMHH$ and $\HM,$ respectively, provided that $N$ is flat 
as a left $Q$-module.
\endproclaim

\Proof.
Put $M'=\Rat(M)Q$. So $M'$ is an $\HdratMQ$-subobject of $M$, and $M=QM'\!$. 
By the left $Q$-flatness of $N$ the canonical map $\,M'\ot_QN\to M\ot_QN\,$ 
is injective, while its image generates $M\ot_QN$ as a left $Q$-module. Hence 
$$
\Rat(M\ot_QN)\cong\Rat(M'\ot_QN)\qquad\text{and}\qquad
(M\ot_QN)^{\Hd}\!\cong(M'\ot_QN)^{\Hd}
$$
by Lemma 4.1. So we may apply Lemma 4.4 noting that $\,\Rat(M')\cong\Rat(M)\,$ 
and ${M'}^{\Hd}\!\cong M^{\Hd}\!$. If $\,x\in M^{\Hd}\!$ and $y\in N^{\Hd}\!$, 
then $\,x\ot_Qy\in(M\ot_QN)^{\Hd}$ and
$$
\eqalign{
h\triangleright(x\ot_Qy)
=\sum h\1x\ot_QyS(h\2)
&=\sum h\1xS(h\2)\ot_Qh\3yS(h\4)\cr
&=\sum\,(h\1\triangleright x)\ot_Q(h\2\triangleright y)
}\eqno(4.6)
$$
for all $h\in H$. This explains that the second map is an isomorphism 
in the monoidal category of left $H$-modules.
\endproof

\proclaim
Corollary 4.6.
Let $\,M_2,\,M_3\in\HdratQMQ\,,$ and let $M_1$ be an object of either 
$\,\HdratMQ\,$ or $\,\HdratQMQ$. Then
$$
\Rat(M_1\ot_QM_2\ot_QM_3)\cong\Rat(M_1)\ot_H\Rat(M_2)\ot_H\Rat(M_3)
$$
and\quad
$(M_1\ot_QM_2\ot_QM_3)^{\Hd}\cong M_1^{\Hd}\!\ot M_2^{\Hd}\!\ot M_3^{\Hd}\!$.
\endproclaim

\Proof.
By Lemma 4.3 both $M_2$ and $M_3$ are free as left $Q$-modules. Hence so too is 
$M_2\ot_QM_3$, and we may apply Lemmas 4.4 and 4.5.
\endproof

Next we consider $\Hd$-module corings, as defined in section 1. A left 
$\Hd$-module $Q$-coring is called \emph{rationally generated} if it is 
rationally generated as an object of the category $\HdQMQ$.

\proclaim
Proposition 4.7.
The equivalence described in Proposition 4.2 induces an equivalence between 
the category of rationally generated left $\Hd$-module $Q$-corings and the 
category of left $H$-module coalgebras.
\endproclaim

\Proof.
Suppose that $\scrC$ is a rationally generated left $\Hd$-module $Q$-coring. 
By Lemma 4.5 the comultiplication $\De_\scrC:\scrC\to\scrC\ot_Q\scrC$ is taken 
by the functor $\Rat$ to a morphism
$$
\De':\ \Rat(\scrC)\to\Rat(\scrC)\ot_H\Rat(\scrC)
$$ 
in $\HMHH$, and coassociativity of $\De_\scrC$ implies coassociativity of $\De'$ in 
view of Corollary 4.6. Since $\Rat(Q)=H$, the counit $\ep_\scrC:\scrC\to Q$ 
induces a morphism $\ep'\!:\Rat(\scrC)\to H$ in $\HMHH$ which satisfies the 
counit property for $\De'$ by functoriality of $\Rat$. Thus $\Rat(\scrC)$ is an 
$H$-coring whose structure maps are morphisms in $\HMHH$. Therefore $\Rat(\scrC)$ 
is a coalgebra in the monoidal category $\HMHH$, and by (4.2) we have
$$
\scrC\cong Q\ot_H\Rat(\scrC)\ot_HQ.
$$

Conversely, for any $H$-coring $\scrC'$ the $Q$-bimodule $\scrC=Q\ot_H\scrC'\ot_HQ$ 
is a $Q$-coring in a natural way. Furthermore, if $\scrC'$ is a coalgebra in 
$\HMHH$, then $\scrC$ is a rationally generated left $\Hd$-module $Q$-coring, 
and $\Rat(\scrC)\cong\scrC'$.

Thus the functor $\,\Rat\,$ induces an equivalence between the category of 
rationally generated left $\Hd$-module $Q$-corings and the category of 
coalgebras in the monoidal category $\HMHH$.

The equivalence of categories $\HMHH\approx\HM$ given by the functor 
$N\mapsto N\co H$ is monoidal, as is seen by a check essentially repeating 
formula (4.6). Therefore under this equivalence coalgebras in $\HMHH$ 
correspond to coalgebras in $\HM$. But coalgebras in the monoidal category of 
left $H$-modules are precisely left $H$-module coalgebras.
\endproof

\proclaim
Corollary 4.8.
Let $\scrC$ be a rationally generated left $\Hd$-module $Q$-coring, and let 
$C=\scrC^{\Hd}\!$ be the corresponding left $H$-module coalgebra. The 
assignment $\,\scrI\mapsto\scrI^{\Hd}$ gives a bijection between the sets of 
$\Hd$-invariant coideals of $\scrC$ and $H$-invariant coideals of $C$.
\endproclaim

\Proof.
By Proposition 4.2 the $\Hd$-invariant $Q$-subbimodules of $\scrC$ are in a 
bijective correspondence with the $H$-submodules of $C$. Let $\scrI$ be such a 
subbimodule of $\scrC$, and let $I=\scrI^{\Hd}\!$. Note that 
$\scrC\,/\,\scrI\in\HdratQMQ$ and $\,(\scrC\,/\,\scrI\,)^{\Hd}\!\cong C/I$, 
again by the equivalence of Proposition 4.2. By the equivalence of Proposition 
4.7 it follows that $\scrC\,/\,\scrI$ is a factor coring of $\scrC$ if and only 
if $C/I$ is a factor coalgebra of $C$. Hence $\scrI$ is a coideal of $\scrC$ if 
and only if $I$ is a coideal of $C$.
\endproof

\proclaim
Lemma 4.9.
For each right $H$-comodule $U$ the $k$-linear map $\,U\to(U\ot H)\co H\,$ 
defined by the assignment $\,u\mapsto\sum u\0\ot S(u\1)\,$ is bijective.
\endproclaim

\Proof.
Denote by $U\triv$ the underlying vector space of $U$ equipped with the 
trivial right $H$-comodule structure $u\mapsto u\ot1$. There is an isomorphism 
$$
\ph:U\ot H\to U\triv\ot H
$$
in the monoidal category $\MH$ defined by the assignment 
$\,u\ot h\mapsto\sum u\0\ot u\1h$ with the inverse map sending $u\ot h$ to 
$\,\sum u\0\ot S(u\1)h\,$ for $u\in U$ and $h\in H$.

\smallskip
Note that $(U\triv\ot H)\co H\cong U\ot H\co H\cong U$ since $H\co H=k$. Hence 
$\ph^{-1}$ induces a $k$-linear bijection of $U$ onto $(U\ot H)\co H$.
\endproof

\proclaim
Proposition 4.10.
Let $\scrC=Q\ot Q$. The corresponding left $H$-module coalgebra
$\scrC^{\Hd}\!$ is canonically isomorphic to $H$ with the $H$-module 
structure on $H$ given by left multiplications.
\endproclaim

\Proof.
By Lemma 4.1 $\Rat(\scrC)=H\ot H$. Hence $\scrC^{\Hd}\!\cong(H\ot H)\co H\!$, 
and so by Lemma 4.9 the $k$-linear map $\,\psi:H\to\scrC^{\Hd}\!$ defined 
by the rule
$$
\psi(x)=\sum\,x\1\ot S(x\2),\qquad x\in H,\eqno(4.7)
$$
is bijective. Furthermore,
$\ \psi(yx)=\sum y\1x\1\ot S(x\2)S(y\2)=y\triangleright\psi(x)$, 
$$
\eqalign{
\De_\scrC\bigl(\psi(x)\bigr)=\sum\,x\1\ot 1\ot S(x\2)
&=\sum\,x\1\ot S(x\2)\,x\3\ot S(x\4)\cr
&=\sum\,\psi(x\1)\ot_Q\psi(x\2),
}
$$
and $\,\ep_\scrC\bigl(\psi(x)\bigr)=\sum\,x\1\,S(x\2)=\ep(x)\,$ for all 
$x,y\in H$. Thus $\psi$ is an isomorphism of left $H$-module coalgebras.
\endproof

By Corollary 4.8 and Proposition 4.10 the coideal left ideals of the Hopf 
algebra $H$ are in a canonical bijective correspondence with the 
$\Hd$-invariant coideals of the left $\Hd$-module $Q$-coring $Q\ot Q$, as 
discussed in section 1.

\proclaim
Corollary 4.11.
Let $\scrC=(Q\ot Q)/\,\scrI$ where $\scrI$ is a left $\Hd$-invariant coideal of 
the $Q$-coring $Q\ot Q,$ and let $g\in\scrC$ be the image of the element\/ 
$1\ot1\in Q\ot Q$. Then the map $\psi_\scrC:H\to\scrC^{\Hd}\!$ defined by the 
rule
$$
\psi_\scrC(x)=\sum\,x\1\,g\,S(x\2),\qquad x\in H,\eqno(4.8)
$$
is a surjective homomorphism of left $H$-module coalgebras. Its kernel $\,I$ 
is the coideal left ideal of $H$ corresponding to $\,\scrI$ under the 
bijection of Theorem 1.1. Hence $\psi_\scrC$ induces an isomorphism of left 
$H$-module coalgebras $\,H/I\to\scrC^{\Hd}\!$.
\endproclaim

\Proof.
The map $\psi_\scrC$ is the composite of the bijection 
$\,\psi:H\to(Q\ot Q)^{\Hd}\!$ defined by formula (4.7) and the canonical 
surjection $\,\pi:(Q\ot Q)^{\Hd}\!\to\scrC^{\Hd}\!$, both of which are 
homomorphisms of left $H$-module coalgebras. Clearly $\psi$ maps $I$ 
bijectively onto the kernel $\,\scrI^{\Hd}\!$ of $\pi$.
\endproof

Let $\scrC$ be a left $\Hd$-module $Q$-coring. An \emph{$\Hd$-equivariant 
right $\scrC$-comodule} $M$ is an $\Hd$-equivariant right $Q$-module and a right 
$\scrC$-comodule such that the comodule structure map $M\to M\ot_Q\scrC\,$ is 
$\Hd$-linear and thus is a morphism in the category $\HdMQ$. Furthermore, we 
say that such a comodule is \emph{rationally generated} if it is rationally 
generated as an object of $\HdMQ$.

We denote by $\HdMC$ the category of $\Hd$-equivariant right $\scrC$-comodules 
and by $\HdratMC$ its full subcategory of rationally generated objects. The 
left comodule categories $\HdCM$ and $\HdratCM$ are defined similarly.

\proclaim
Proposition 4.12.
Let $\scrC$ be a rationally generated left $\Hd$-module $Q$-coring, and let 
$\,C=\scrC^{\Hd}\!$ be the corresponding left $H$-module coalgebra. The 
functor $\,M\mapsto M^{\Hd}$ gives category equivalences 
$\,\HdratMC\approx\MC$ and $\,\HdratCM\approx\CM$.
\endproclaim

\Proof.
Proposition 3.2 describes a category equivalence $\HdratMQ\approx\Mk$. Given 
an object $M\in\HdratMQ$, we have $M\cong M^{\Hd}\!\ot Q$. It follows that for 
each object $N$ of the category $\HdMQ$, even if $N$ is not rationally 
generated, the morphisms $M\to N$ in $\HdMQ$ are in a bijective 
correspondence with the $k$-linear maps $M^{\Hd}\!\!\to N^{\Hd}\!$.
In particular, in view of Lemma 4.4, the $\HdMQ$-morphisms $M\to M\ot_Q\scrC$ 
are in a bijective correspondence with the $k$-linear maps 
$$
M^{\Hd}\!\!\to(M\ot_Q\scrC)^{\Hd}\cong M^{\Hd}\!\!\ot C.
$$
Taking into account Corollary 4.6, we see that 
commutativity of the diagrams
$$
\diagram{
M&\hidewidth\lmapr8\rho\hidewidth&M\ot_Q\scrC&\qquad\qquad&M&\lmapr2\rho&M\ot_Q\scrC\cr
\noalign{\smallskip}
\lmapd{14}\rho{}&&\lmapd{14}{}{\id\ot\De_\scrC}&&\lmapd{14}\id{}&&\lmapd{14}{}{\id\ot\ep_\scrC}\cr
\noalign{\smallskip}
M\ot_Q\scrC&\lmapr4{\rho\ot\id}&M\ot_Q\scrC\ot_Q\scrC,&&M&\cong&M\ot_QQ\cr
}
$$
required for $\scrC$-comodule structure maps translates into commutativity of 
similar diagrams defining $C$-comodule structures.

The equivalence of left comodule categories follows by symmetry.
\endproof

\section
5. Coaction invariants and the second bijection

Let $\scrC$ be a left $\Hd$-module factor coring of $Q\ot Q$. Thus 
$\scrC=(Q\ot Q)/\,\scrI$ where $\scrI$ is a left $\Hd$-invariant coideal of 
$\,Q\ot Q$. Denote by $\,\pi:Q\ot Q\to\scrC\,$ the canonical surjective 
homomorphism of left $\Hd$-module $Q$-corings.

The $\Hd$-invariant element $g=\pi(1\ot1)$ satisfies $\De_\scrC(g)=g\ot_Qg$ and 
$\ep_\scrC(g)=1$. We call $g$ the \emph{distinguished grouplike} of $\scrC$ 
and use it to define a right $\scrC$-comodule structure 
$\,\rho_Q:Q\to Q\ot_Q\scrC\cong\scrC\,$ on $Q$ by the rule
$$ 
\rho_Q(x)=gx,\qquad x\in Q.\eqno(5.1)
$$
Note that $\rho_Q$ is an $\Hd$-linear map, and so $Q$ is an object of 
the category $\HdMC$ of $\Hd$-equivariant right $\scrC$-comodules introduced 
at the end of section 4. In a similar way $Q$ is an object of the category 
$\HdCM$ of $\Hd$-equivariant left $\scrC$-comodules with respect to the left 
comodule structure determined by $g$.

The set $M\co\scrC$ of coaction invariants of a right $\scrC$-comodule $M$ is 
the equalizer of the two maps
$$
\rho,\,\tau:\,M\longrightarrow M\ot_Q\scrC\eqno(5.2)
$$
where $\rho$ is the comodule structure on $M$ and $\tau$ is defined by the 
rule $\,\tau(x)=x\ot_Qg$ for $x\in M$. In other words,
$$
M\co\scrC=\{x\in M\mid\,\rho(x)=x\ot_Qg\}.\eqno(5.3)
$$
For each left $\scrC$-comodule $M$ the set of coaction invariants 
$\lco\scrC M$ is defined similarly. By (5.1) and (5.3) we have
$$
\lco\scrC Q=Q\co\scrC=\{x\in Q\mid\,gx=xg\}.\eqno(5.4)
$$
This is a subalgebra of $Q$ such that for each right $\scrC$-comodule $M$ both 
maps in (5.2) are right $Q\co\scrC$-linear (the comodule structure map $\rho$ 
is right $Q$-linear by definition), and therefore $M\co\scrC$ is a 
$Q\co\scrC$-submodule of $M$.

\proclaim
Lemma 5.1.
We have $\,Q\co\scrC=\{x\in Q\mid 1\ot x-x\ot1\in\scrI\,\}$. Thus 
$\,Q\co\scrC$ is the left $\Hd$-invariant artinian subalgebra of\/ $Q$ 
corresponding to the coideal\/ $\scrI$ under the assignment of section 1.  
\endproclaim

\Proof.
Indeed, since $\,\pi:Q\ot Q\to\scrC\,$ commutes with the left and right 
actions of $Q$, the equality $gx=xg$ in $\scrC$ means precisely that 
$\,1\ot x-x\ot1\in\Ker\pi=\scrI$.
\endproof

\proclaim
Lemma 5.2.
Suppose that\/ $\scrC=Q\ot_AQ$ where $A$ is a left $\Hd$-invariant subalgebra 
of\/ $Q$. Then $\,Q\co\scrC\,$ coincides with the dominion of $A$ in $Q$. In 
particular, $\,Q\co\scrC=A\,$ when $A$ is artinian.
\endproclaim

\Proof.
In this case the equality $gx=xg$ in $\scrC$ is rewritten as $1\ot_Ax=x\ot_A1$. 
Such elements $x$ constitute the dominion of $A$ in $Q$. If $A$ is artinian, 
then $A$ coincides with its dominion in $Q$ by Proposition 3.6.
\endproof

Lemma 5.2 shows that the composite $\,\Lat_3\to\Lat_2\to\Lat_3$ of the two 
maps defined in section 1 is the identity map. To ensure that these maps are 
bijective it remains to prove that each left $\Hd$-module factor coring of 
$Q\ot Q$ has the form $Q\ot_AQ$ for some left $\Hd$-invariant artinian 
subalgebra $A$ of $Q$.

\smallskip
So we assume further that $\scrC=(Q\ot Q)/\,\scrI$ is an arbitrary left 
$\Hd$-module factor coring of $\,Q\ot Q$. Let $\,A=Q\co\scrC$.

Since $1\ot x-x\ot1\in\scrI$ for all $x\in A$ by Lemma 5.1, the canonical 
homomorphism $\,\pi:Q\ot Q\to\scrC\,$ factors through the coring $Q\ot_AQ$. 
Thus $\pi$ induces a surjective homomorphism of left $\Hd$-module $Q$-corings
$$
\pi':Q\ot_AQ\to\scrC.
$$
Bijectivity of this $\pi'$ is closely related to equivalences between certain 
categories of modules and comodules. By \cite{Br-W, 28.8} there is a pair of 
adjoint functors
$$
\openup1\jot 
\vcenter{\halign{\hfil$#$\hfil&\qquad\qquad\hfil$#$\hfil\cr 
\MA\to\MscrC&\MscrC\to\MA\cr
W\mapsto W\ot_AQ\,,&M\mapsto M\co\scrC.\cr
}}\eqno(5.5)
$$
The first functor makes $\,W\ot_AQ\,$ into a $\scrC$-comodule by means of 
the structure map
$$
W\ot_A\rho_Q:\ W\ot_AQ\to W\ot_A\scrC\cong(W\ot_AQ)\ot_Q\scrC\eqno(5.6)
$$
where $\rho_Q:Q\to\scrC$ is given by (5.1), while $M\co\scrC$ is an 
$A$-submodule of $M$. If $W$ is an object of $\HdMA$, then tensoring with $W$ 
over $A$ gives a functor 
$$
W\ot_A{}?:\HdAMQ\to\HdMQ\,.
$$
So $W\ot_AQ$ has an $\Hd$-module structure, and the comodule structure map
(5.6) is $\Hd$-linear. Conversely, if $M\in\HdMC$, then the $A$-submodule 
$M\co\scrC\sbs M$ is stable also under the action of $\Hd$ since both maps in 
(5.2) are $\Hd$-linear. This shows that the functors given in (5.5) induce 
a pair of adjoint functors
$$
\HdMA\to\HdMC,\qquad\HdMC\to\HdMA.\eqno(5.7)
$$

\proclaim
Lemma 5.3.
If\/ $W$ is a flat right $A$-module, then $\,(W\ot_AQ)\co\scrC\cong W$.
\endproclaim

\Proof.
There is an exact sequence of left $A$-modules $\,0\to A\to Q\lmapr1\ph\scrC\,$
where $\ph$ is defined by the rule $\,\ph(x)=\rho_Q(x)-xg=gx-xg\,$ for $x\in Q$. 
Applying to it the exact functor $\,W\ot_A{}?$, we obtain the desired conclusion.
\endproof

\setitemsize(a)
\proclaim
Lemma 5.4.
The following conditions are equivalent:

\item(a)
\ the homomorphism of left $\Hd$-module $Q$-corings 
$\,\pi':Q\ot_AQ\to\scrC\,$ is bijective,

\item(b)
\ the functors in {\rm(5.7)} are quasi-inverse equivalences,

\item(c)
\ $M=M\co\scrC\!\cdot Q\,$ for each object $\,M\in\HdMC,$

\item(d)
\ $M\co\scrC\!\ne0\,$ for each nonzero $Q$-finite object $\,M\in\HdMC$.

\endproclaim

\Proof.
Suppose (a) holds, i.e., $\scrC\cong Q\ot_AQ$. The category $\MscrC$ 
of right comodules for 
this coring is equivalent to the category of descent data for the ring 
extension $A\sbs Q$ (see \cite{Br-W, 25.4}). By Proposition 3.6 $Q$ is a free 
$A$-module on both side. In particular, $Q$ is faithfully flat over $A$. It is 
now a classical fact of noncommutative descent theory that the functors in 
(5.5) are quasi-inverse equivalences (see \cite{Br-W, 28.19} and \cite{Cae-MZ, 
Prop. 109}). Then (b) also follows.

We have shown that (a)$\,\Rar\,$(b). Condition (b) implies that 
$\,M\cong M\co\scrC\ot_AQ\,$ for each $M\in\HdMC$. Hence (b)$\,\Rar\,$(c). 
The implication (c)$\,\Rar\,$(d) is obvious.

Suppose that condition (d) holds. Since $\Ker\pi'$ is a $Q$-subbimodule of 
$Q\ot_AQ$, it is generated by its intersection with the $(H,Q)$-subbimodule 
$\,HA\ot_AQ\sbs Q\ot_AQ$. Taking various finite-dimensional right coideals $U$ 
of $H$, we obtain a directed family of $\Hd$-invariant $Q$-finite right 
coideals $UA\ot_AQ$ of the $Q$-coring $Q\ot_AQ$ whose union coincides 
with $\,HA\ot_AQ$. Injectivity of $\pi'$ will follow once we show that the 
restriction of $\pi'$ to each of these right coideals is injective.

So let $M=UA\ot_AQ$ be one of right coideals in that family. We may view $M$ as a 
$Q$-finite object of the category $\HdMC$. Then $\pi'|_M:M\to\scrC$ is a morphism 
in $\HdMC$ and its kernel $K=\Ker\pi'|_M$ is an $\HdMC$-subobject of $M$. 
Furthermore, $K$ is $Q$-finite since $Q$ is artinian. We have 
$$
K\co\scrC=K\cap M\co\scrC.
$$
Note that $UA$ is a projective right $A$-module by Theorem 2.1 since $UA$ is 
an $A$-finite object of the category $\HdMA$ and $A$ is an artinian $\Hd$-simple 
$\Hd$-module algebra. Hence $M\co\scrC\cong UA$ by Lemma 5.3, and so 
$\,M\co\scrC\sbs Q\ot_A1\sbs Q\ot_AQ$.

Since $\pi'$ is a homomorphism of corings, the composite $\ep_\scrC\circ\pi'$ 
is the counit of the coring $Q\ot_AQ$, i.e., 
$\ep_\scrC\bigl(\pi'(x\ot_Ay)\bigr)=xy$ for $x,y\in Q$. It follows that 
$\ep_\scrC:\scrC\to Q$ is a retraction of the map $Q\to\scrC$ given by the 
assignment $x\mapsto\pi'(x\ot1)$. Hence $\pi'$ is injective on 
$\,Q\ot_A1$, and therefore on $M\co\scrC$. But this forces 
$K\co\scrC=0$, whence $K=0$ by (d). Thus $\pi'|_M$ is indeed 
injective.

This shows that $\pi'$ is injective. But $\pi'$ is also surjective, and so the 
implication (d)$\,\Rar\,$(a) is also proved.
\endproof

The remaining step in establishing bijectivity of the correspondence described 
in section 1 consists in proving that each left $\Hd$-module factor coring of 
$Q\ot Q$ does satisfy condition (d) of Lemma 5.4. It is not easy and will be 
accomplished in the more general setup of section 7. The final conclusion 
confirming the property $M\co\scrC\ne0$ will be presented in Corollary 7.10. 
As a consequence, we arrive at

\proclaim
Proposition 5.5.
If $\scrC$ is any left $\Hd$-module factor coring of\/ $Q\ot Q$ and 
$A=Q\co\scrC\!,$ then $\,\scrC\cong Q\ot_AQ$.
\endproclaim

In some cases Proposition 5.5 can be proved very quickly without the extra 
work done in section 7. This approach is based on the following observation:

\proclaim
Lemma 5.6.
Suppose that for each $Q$-finite object $M\in\HdratMQ$ and each rational 
$\Hd$-submodule $U\sbs M$ which generates $M$ there exists a basis of $M$ 
over $Q$ contained in $U$. Then the conclusion of Proposition 5.5 is true.
\endproclaim

\Proof.
Let $M=UA\ot_AQ$ where $U$ is a finite-dimensional right coideal of $H$. The 
image $L=\pi'(M)$ of $M$ in $\scrC$ is an $\Hd$-invariant right coideal of 
$\scrC$ and $\,\pi'|_M:M\to L\,$ is an epimorphism in $\HdratMC$. Since 
$x\ot_A1\in M\co\scrC$ for all $x\in UA$, the formula
$$
\ph(x)=\pi'(x\ot_A1),\qquad x\in UA,
$$
defines a map $\ph:UA\to L\co\scrC$ which is $\Hd$-linear and $A$-linear, 
i.e., a morphism in $\HdMA$. The set $\{x\ot_A1\mid x\in U\}$ generates $M$ as 
a $Q$-module. Hence $L=\ph(U)\mskip2mu Q$, i.e., $L$ is generated by its 
rational $\Hd$-submodule $\ph(U)$. By the hypothesis there exists a set of 
elements $e_1,\ldots,e_n\in U$ whose images $e'_i=\ph(e_i)$ form a basis of 
$L$ over $Q$. Each element $y\in L$ can be written as $\,y=\sum e'_i\,y_i\,$ 
with $y_1,\ldots,y_n\in Q$. Since $e'_i\in L\co\scrC$, we have then 
$$
\De_\scrC(y)=\sum e'_i\ot_Qgy_i\in L\ot_Q\scrC,
$$
and $y\in L\co\scrC$ if and only if
$$
\De_\scrC(y)=y\ot_Qg=\sum e'_i\ot_Qy_ig.
$$
Since the chosen basis $e'_1,\ldots,e'_n$ of $L$ provides the decomposition 
of $L\ot_Q\scrC$ as a direct sum of vector subspaces 
$e'_iQ\ot_Q\scrC\cong\scrC$, we deduce that $y\in L\co\scrC$ if and only if 
$y_ig=gy_i$, i.e., $y_i\in Q\co\scrC=A$ for each $i$. This shows that 
$L\co\scrC$ is the $A$-submodule of $L$ generated by $e'_1,\ldots,e'_n$.

It follows that $\ph$ maps the $A$-submodule of $UA$ generated by 
$e_1,\ldots,e_n$ onto the whole $L\co\scrC$. As we have seen in the proof of 
Lemma 5.4 the fact that $\pi'$ is a homomorphism of corings implies that the 
map $\ph$ is injective. Then $e_1,\ldots,e_n$ must generate $UA$ as an 
$A$-module. Hence the elements $e_1\ot_A1,\ldots,e_n\ot_A1$ generate $M$ as a 
$Q$-module. Since their images $e'_1,\ldots,e'_n\in L$ are linearly 
independent over $Q$, we conclude that the restriction of $\pi'$ to $M$ is 
injective. Then $\pi'$ is injective, and therefore bijective, as in the proof 
of the implication (d)$\,\Rar\,$(a) of Lemma 5.4.
\endproof

By Proposition 3.2 we know that all objects of the category $\HdratMQ$ are 
free $Q$-modules. So the question is whether a basis for a free $Q$-module $M$ 
can be found in a vector subspace $U$ such that $M=UQ$. This indeed can be 
done when $Q$ is local or, even more generally, when all simple factor rings 
of the artinian ring $Q$ are skew fields and the base field $k$ is either 
infinite or finite of cardinality not less than the number of maximal ideals 
of $Q$.

Furthermore, by considering the chain of subobjects of $M$ generated by the 
$\Hd$-submodules in a composition series for $U$, the property of $M$ in the 
hypothesis of Lemma 5.6 reduces to the case when $U$ is a simple rational 
$\Hd$-module. If $H$ is pointed, then such a module has dimension 1. Then $M$ 
is a free $Q$-module of rank 1, and any nonzero element of $U$ forms a basis 
of $M$ over $Q$. So Lemma 5.6 applies in this case too.

It is reasonable to ask whether the property in the hypothesis of Lemma 5.6 is 
satisfied for other Hopf algebras, although this property seems to be quite 
strong. For $M=Q$ it implies that each nonzero right coideal of $H$ contains a 
nonzerodivisor of $H$. Specializing even further to the case when 
$\,\dim H<\infty$, we ask

\proclaim
Question 5.7.
Suppose that $H$ is a finite-dimensional Hopf algebra. Is it true that every 
nonzero one-sided coideal of $H$ contains an invertible element of $H$?
\endproclaim

\section
6. Twisting of comodules over $H$-module corings

In this section $H$ stands for an arbitrary Hopf algebra over the base field 
$k$. Let $R$ be a left $H$-module algebra and $\scrC$ a left $H$-module 
$R$-coring. In section 2 we recalled the twisting functors $U\ot\,?$ 
associated with right $H$-comodules. These functors are defined, in 
particular, on the category $\MR$ of right $R$-modules. The present section 
aims to extend the twisting functors to the category $\MscrC$ of right 
$\scrC$-comodules.

Let $U\in\MH$. By Lemma 2.4 there are natural isomorphisms
$$
\ph_{U,V,N}:U\ot(V\ot_RN)\to(U\ot V)\ot_RN\eqno(6.1)
$$
in $\MR$ for objects $V\in\MR$ and $N\in\HRMR$. In particular, $\ph_{U,V,N}$ 
is defined for $N=\scrC$ and for $N=\scrC\ot_R\scrC$.

\proclaim
Lemma 6.1.
Let $\rho:V\to V\ot_R\scrC$ be the structure map of a right $\scrC$-comodule 
$V$. Then the twisted right $R$-module $U\ot V$ is a right $\scrC$-comodule 
with structure map $\rho_U$ defined as the composite
$$
U\ot V\lmapr4{\id\ot\rho}U\ot(V\ot_R\scrC)
\lmapr4{\ph_{U,V,\scrC}}(U\ot V)\ot_R\scrC\eqno(6.2)
$$
\endproclaim

\Proof.
Both maps in (6.2) are morphisms in $\MR$. Hence so too is $\rho_U$. The 
composite of $\rho_U$ and the map 
$$
\id_{U\ot V}\ot\ep_\scrC:\ (U\ot V)\ot_R\scrC\to U\ot V
$$
is fully expressed as the composite of maps in the following diagram:
$$
\cdiagram{
U\ot V&\lmapr4{\id\ot\rho}&
U\ot(V\ot_R\scrC)&\lmapr8{\id\ot(\id\ot\ep_\scrC)}&U\ot(V\ot_RR)&{}\cong{}&U\ot V\cr
\noalign{\smallskip}
&&\lmapd{16}{\ph_{U,V,\scrC}}{}&&\lmapd{16}{}{\ph_{U,V,R}}&&\lmapd{16}{}{\id}\cr
\noalign{\smallskip}
&&(U\ot V)\ot_R\scrC&\lmapr8{(\id\ot\id)\ot\ep_\scrC}&(U\ot V)\ot_RR&{}\cong{}&U\ot V.\cr
}
$$
The middle square in this diagram commutes by naturality of (6.1) since the 
counit $\ep_\scrC:\scrC\to R$ is a morphism in the category $\HRMR$, while the 
right square commutes since $\ph_{U,V,R}$ sends 
$u\ot(v\ot_R1)$ to $(u\ot v)\ot_R1$ for each $u\in U$ and each $v\in V$. Since 
$\,(\id\ot\ep_\scrC)\circ\rho\,$ is the identity endomorphism of $V$, it 
follows that $\,(\id_{U\ot V}\ot\ep_\scrC)\circ\rho_U\,$ is the identity 
endomorphism of $U\ot V$.

By naturality of (6.1) there are also commutative diagrams
$$
\cdiagram{
U\ot(V\ot_R\scrC)&\lmapr9{\id\ot(\id\ot\De_\scrC)}&U\ot(V\ot_R\scrC\ot_R\scrC)\cr
\noalign{\smallskip}
\lmapd{16}{\ph_{U,V,\scrC}}{}&&\lmapd{16}{}{\ph_{U,V,\mskip1mu\scrC\ot_R\scrC}}\cr
\noalign{\smallskip}
(U\ot V)\ot_R\scrC&\lmapr9{(\id\ot\id)\ot\De_\scrC}&(U\ot V)\ot_R\scrC\ot_R\scrC\cr
}\eqno(6.3)
$$
and
$$
\cdiagram{
U\ot(V\ot_R\scrC)&\lmapr6{\id\ot(\rho\ot\id)}&U\ot(V\ot_R\scrC\ot_R\scrC)\cr
\noalign{\smallskip}
\lmapd{16}{\ph_{U,V,\scrC}}{}&&\lmapd{16}{}{\ph_{U,V\ot_R\scrC,\mskip1mu\scrC}}\cr
\noalign{\smallskip}
(U\ot V)\ot_R\scrC&\lmapr6{(\id\ot\rho)\ot\id}&\bigl(U\ot(V\ot_R\scrC)\bigr)\ot_R\scrC
&\lmapr7{\ph_{U,V,\scrC}\ot\id}&(U\ot V)\ot_R\scrC\ot_R\scrC\cr
}\eqno{\hskip-20pt(6.4)}
$$
Note that the composite of the two bottom maps in (6.4) is the map 
$\rho_U\ot\id$, while
$$
(\ph_{U,V,\scrC}\ot\id)\circ\ph_{U,\,V\ot_R\scrC,\,\scrC}=\ph_{U,V,\,\scrC\ot_R\scrC}\,.
$$
To check the last formula let $u\in U$, $v\in V$, and $x,y\in\scrC$. Applying 
the definition of $\ph$ in Lemma 2.4, we see that the element 
$\,u\ot(v\ot_Rx\ot_Ry)\in U\ot(V\ot_R\scrC\ot_R\scrC)\,$ is taken to
$$
\sum\,\bigl(u\0\ot(v\ot_Rx)\bigr)\ot_Ru\1y\in\bigl(U\ot(V\ot_R\scrC)\bigr)\ot_R\scrC 
$$
by $\,\ph_{U,V\ot_R\scrC,\,\scrC}\,$, and this element is taken further to
$$
\sum\bigl(u\0\ot v)\ot_Ru\1x\ot_Ru\2y
=\sum\bigl(u\0\ot v)\ot_Ru\1(x\ot_Ry)\in(U\ot V)\ot_R\scrC\ot_R\scrC
$$
by $\,\ph_{U,V,\scrC}\ot\id$.

Since the image of $\ \id\ot\rho:U\ot V\to U\ot(V\ot_R\scrC)\,$ is contained 
in the equalizer of the top horizontal maps appearing in diagrams (6.3) and 
(6.4), it follows that $\,\rho_U=\ph_{U,V,\scrC}\circ(\id\ot\rho)\,$ has images 
in the equalizer of the bottom maps in these two diagrams. In other words,
$$
(\id_{U\ot V}\ot\De_\scrC)\circ\rho_U=(\rho_U\ot\id_\scrC)\circ\rho_U,
$$
the coassociativity law for $\rho_U$.
\endproof

\proclaim
Lemma 6.2.
The operation of twisting right $\scrC$-comodules by right $H$-comodules 
makes $\MscrC$ into a left module category over the monoidal category $\MH$.
\endproclaim

\Proof.
The construction of the twisted right comodules $U\ot V$ presented in Lemma 6.1 
defines a functor $\MH\times\MscrC\to\MscrC$, additive in each argument. Given 
two right $H$-comodules $U$ and $U'$, the canonical $k$-linear bijection 
$$
U\ot(U'\ot V)\cong(U\ot U')\ot V\eqno(6.5)
$$
is an isomorphism in $\MR$. Moreover, it is compatible with the right 
$\scrC$-comodule structures defined by (6.2), as is seen from the 
commutative diagram
$$
\cdiagram{
U\ot(U'\ot V)&\lmapr7{\id\ot(\id\ot\rho)}&U\ot(U'\ot(V\ot_R\scrC))
&\lmapr7{\id\ot\ph_{U'\!,V,\scrC}}&U\ot((U'\ot V)\ot_R\scrC)\cr
\noalign{\smallskip}
\lmapd{14}{\cong}{}&&\lmapd{14}{\cong}{}&&\lmapd{14}{}{\ph_{U,\,U'\ot V,\scrC}}\cr
\noalign{\smallskip}
(U\ot U')\ot V&\lmapr7{(\id\ot\id)\ot\rho}&(U\ot U')\ot(V\ot_R\scrC)
&\lmapr7{\ph_{U\ot U'\!,V,\scrC}}&(U\ot U'\ot V)\ot_R\scrC.\cr
}
$$
Commutativity of the right square is checked by observing that 
$\,u\ot(u'\ot(v\ot_Rc))\in U\ot(U'\ot(V\ot_R\scrC))\,$ is mapped to
$$
\sum\,(u\0\ot u'\0\ot v)\ot_Ru\1u'\1c\in(U\ot U'\ot V)\ot_R\scrC
$$
by both composite maps in that square. Thus (6.5) is an isomorphism in 
$\MscrC$, and so is the canonical bijection $\,k\triv\ot V\cong V\,$ where 
$k\triv$ is the one-dimensional trivial $H$-comodule.
\endproof

If $U$ is a finite-dimensional right $H$-comodule, then its dual vector space 
$U^*$ has a right $H$-comodule structure which makes $U^*$ the left dual of 
$U$ in the monoidal category $\MH$.

\proclaim
Corollary 6.3.
Let $U$ be a finite-dimensional right $H$-comodule and $U^*$ its left dual in 
$\MH$. For $V,W\in\MscrC$ there are canonical $k$-linear bijections
$$
\MscrC(U^*\ot V,\,W)\cong\MscrC(V,\,U\ot W).
$$
If $\,W$ is injective in $\MscrC,$ then so too is $\,U\ot W$.
\endproclaim

\Proof.
The stated bijection is a standard property of left duals. It shows that the 
twisting endofunctors $U^*\ot\,?$ and $U\ot\,?$ of the category 
$\MscrC$ form an adjoint pair. Since the left adjoint is exact, the 
right adjoint preserves injectives.
\endproof

The category $\HMscrC$ of $H$-equivariant right $\scrC$-comodules is defined 
as its special case in section 4. An object $M\in\HMscrC$ is an 
$H$-equivariant right $R$-module and a right $\scrC$-comodule such that the 
comodule structure map $M\to M\ot_R\scrC$ is $H$-linear.

\proclaim
Lemma 6.4.
Let $M\in\HMscrC$. Consider the twisted object $H\ot M\in\MscrC$ where $H$ is
viewed as a right $H$-comodule with respect to the comultiplication in $H$. 
The $k$-linear map
$$
\mu:H\ot M\to M
$$
afforded by the $H$-module structure on $M$ is a morphism in $\MscrC$.

For a right $H$-comodule $U$ denote by $U\triv$ the $H$-comodule which has the 
same underlying vector space, but the coaction of $H$ on it is trivial. 
Then
$$
U\ot M\cong U\triv\ot M\quad\text{in $\,\MscrC$}.
$$
In particular, $\,U\ot M\cong M^d$ in $\MscrC$ when $\,d=\dim_kU<\infty$.
\endproclaim

\Proof.
By \cite{Sk-Oy06, Lemma 1.2(iii)} $\,\mu$ is a morphism in $\MR$. Let 
$\rho:M\to M\ot_R\scrC$ be the $\scrC$-comodule structure map. Then $\rho_H$ 
defined by (6.2) is the $\scrC$-comodule structure on $H\ot M$. The desired 
equality $\rho\circ\mu=(\mu\ot\id)\circ\rho_H$ is expressed by means of the 
commutative diagram
$$
\cdiagram{
H\ot M&\lmapr4{\id\ot\rho}H\ot(M\ot_R\scrC)\lmapr6{\ph_{H,M,\scrC}}&(H\ot M)\ot_R\scrC\cr
\noalign{\smallskip}
\lmapd{16}{\mu}{}&&\lmapd{16}{}{\mu\ot\id}\cr
\noalign{\smallskip}
M&\hidewidth\ \lmapr{28}{\rho}\hidewidth&M\ot_R\scrC.\cr
}
$$
Note that for $h\in H$, $x\in M$, and $c\in C$ the element 
$\,h\ot(x\ot_Rc)\in H\ot(M\ot_R\scrC)$ is taken by the composite of 
$\ph_{H,M,\scrC}$ and $\mu\ot\id$ to 
$$
\sum h\1x\ot_Rh\2c=h\cdot(x\ot_Rc)\in M\ot_R\scrC,
$$
i.e., $(\mu\ot\id)\circ\ph_{H,M,\scrC}$ is the map defining the $H$-module 
structure on $M\ot_R\scrC$. Since $\rho$ is a homomorphism of left $H$-modules, 
the diagram does indeed commute.

The comodule structure on $U$ yields a map $\,\de:U\to U\triv\ot H$ which is a 
morphism in $\MH$. It follows that the composite
$$
\xi:\ \ U\ot M\lmapr4{\de\ot\id}U\triv\ot H\ot M\lmapr4{\id\ot\mu}U\triv\ot M
$$
is a morphism in $\MscrC$ since so are both factors by functoriality of the 
twisting. For $u\in M$ and $x\in M$ we have $\,\xi(u\ot x)=\sum u\0\ot u_1x$. 
We see that the assignment $\,u\ot x\mapsto\sum u\0\ot S(u_1)x\,$ defines the 
inverse map $\xi^{-1}$. Therefore $\xi$ is bijective, and so $\xi$ is an 
isomorphism in $\MscrC$.
\endproof

\section
7. A coring structure theorem

In this section we assume that $H$ is an arbitrary Hopf algebra with bijective 
antipode over the base field $k$ and $R$ is a right artinian $H$-simple left 
$H$-module algebra. We aim to find a set of conditions on a left $H$-module 
factor coring $\scrC$ of the canonical $R$-coring $R\ot R$ which imply that 
all nonzero $R$-finite objects $M$ of the category $\HMscrC$ satisfy 
$M\co\scrC\ne0$. Under several further assumptions this will be used to show 
that $\,\scrC=R\ot_AR\,$ for some $H$-invariant subalgebra $A\sbs R$.

The main result stated in Theorem 7.11 fills in the gap left in section 5 and 
so it completes the proof of Theorem 1.1. Intermediate results of this 
section are valid for more general $H$-module corings.

For any left $R$-flat $R$-coring $\scrC$ the category $\MscrC$ or right 
$\scrC$-comodules is abelian, and moreover a Grothendieck category 
\cite{Br-W, 18.14}. In particular, this category contains injective hulls of 
its objects. We say that $V\in\MscrC$ is \emph{$R$-finite} if $V$ is finitely 
generated as a right $R$-module.

Since $R$ is right artinian, each $R$-finite right $\scrC$-comodule has finite 
length, even as an object of $\MR$, and an arbitrary object of $\MscrC$ is the 
union of its $R$-finite subobjects by the Finiteness Theorem \cite{Br-W, 18.16}. 
As a consequence, each object of $\MscrC$ is an essential extension of its 
socle. Any injective object of $\MscrC$ is the injective hull of its socle, 
and therefore a direct sum of indecomposable injectives whose socles are 
simple comodules. Note that arbitrary direct sums of injective 
right $\scrC$-comodules are injective. Furthermore, we have

\proclaim
Lemma 7.1.
Any left $R$-flat $R$-coring $\scrC$ viewed as a right $\scrC$-comodule with 
respect to the comultiplication is an injective cogenerator in $\MscrC$.  
\endproclaim

\Proof.
As a special case of \cite{Br-W, 18.10}, there are natural $k$-linear bijections
$$
\MscrC(V,\scrC)\cong\Hom_R(V,R),\qquad V\in\MscrC.
$$
It was proved in \cite{Sk11, Th. 1.1} that every right artinian $H$-simple 
left $H$-module algebra is a quasi-Frobenius ring. This applies to $R$, and so 
$R$, as a right module over itself, is an injective cogenerator in $\MR$. It 
follows that the functor $\MscrC(?,\scrC)$ is faithfully exact.
\endproof

We will have to use 3 abelian groups:

$G_0(R)$, the Grothendieck group of the category of right $R$-modules of 
finite length,

$G_0(\MscrC)$, the Grothendieck group of the category of $R$-finite right 
$\scrC$-comodules,

$K_0(\MscrC)$, the Grothendieck group of the category of injective right 
$\scrC$-comodules with $R$-finite socles.

\smallskip
The group $G_0(\MscrC)$ is generated by symbols $[V]$ associated with the 
isomorphism classes of $R$-finite right $\scrC$-comodules, with the set of 
defining relations
$$
\,[V]=[V']+[V'']\,
$$
corresponding to various short exact sequences $0\to V'\to V\to V''\to0$ of 
$R$-finite right $\scrC$-comodules. The groups $G_0(R)$ and $K_0(\MscrC)$ are 
defined similarly.

Each of these abelian groups is free. Standard bases for $G_0(R)$ and 
$G_0(\MscrC)$ are formed by the isomorphism classes of simple objects of $\MR$ 
and $\MscrC$, respectively. A standard basis $\,\Indinj\MscrC\,$ for the group 
$K_0(\MscrC)$ consists of the isomorphism classes of indecomposable injectives.

We will have to deal with the following two conditions:

\setitemsize(C2)
\item(C1)
\ the set $\Irr\MscrC$ of isomorphism classes of simple right $\scrC$-comodules 
is finite,

\item(C2)
\ all indecomposable injective right $\scrC$-comodules are $R$-finite.

Condition (C1) means that the free abelian groups $G_0(\MscrC)$ and 
$K_0(\MscrC)$ have finite rank. The rank of $G_0(R)$ is finite by the right 
artinian assumption.

The forgetful functor $\MscrC\to\MR$, being exact, induces a group homomorphism
$$
G_0(\MscrC)\to G_0(\MR).\eqno(7.1)
$$
Condition (C2) leads to another homomorphism $K_0(\MscrC)\to G_0(\MscrC)$, and 
then we may take its composite
$$
K_0(\MscrC)\to G_0(\MR)\eqno(7.2)
$$
with the previous one. Condition (C1) is easily verified for the corings we 
have in mind. Condition (C2) is much more tricky to establish, but this will 
be crucial for the final conclusions.

\medskip
Further on we assume that $\scrC$ is a left $R$-flat left $H$-module $R$-coring. 

In section 6 we introduced the operation of twisting right $\scrC$-comodules 
by a right $H$-comodule $U$. It extends the twisting operation on $R$-modules 
recalled in section 2. Twisting by $U$ defines endofunctors $U\ot{}?$ of the 
abelian categories $\MR$ and $\MscrC$ which are clearly exact. 

If $\,\dim_kU<\infty$, then the twisted $R$-module $U\ot V$ is finitely 
generated whenever so is $V$ \cite{Sk-Oy06, Lemma 1.1}. In this case the 
twisting endofunctors $U\ot{}?$ induce endomorphisms of the Grothendieck groups 
$G_0(R)$ and $G_0(\MscrC)$. By Corollary 6.3 the class of injective right 
$\scrC$-comodules is stable under twisting. We thus obtain also an endomorphism 
of $K_0(\MscrC)$ provided that condition (C2) holds. Canonical maps (7.1) 
and (7.2) commute with respective endomorphisms of $G_0(R)$, $G_0(\MscrC)$, 
and $K_0(\MscrC)$.

These endomorphisms depend only on the set of composition factors of $U$. Thus 
we may view $G_0(R)$ and $G_0(\MscrC)$ as left modules over the Grothendieck 
ring $G_0(\MH)$ of the monoidal category of finite-dimensional right 
$H$-comodules. The dual abelian groups $G_0(R)^*$ and $G_0(\MscrC)^*$ are 
right $G_0(\MH)$-modules in a natural way. The same applies to $K_0(\MscrC)$ 
and $K_0(\MscrC)^*$ when (C2) holds.

For a function $f$ in $G_0(R)^*$ or $G_0(\MscrC)^*$ we denote by $fU$ the 
image of $f$ under the action of $\,[U]\in G_0(\MH)$, i.e., $fU$ is defined by 
formula (2.2) with $x$ in $G_0(R)$ or $G_0(\MscrC)$, respectively.

Special properties of $H$-equivariant $\scrC$-comodules we are interested in 
are related in an essential way to the fact that the endomorphisms of the 
above Grothendieck groups induced by the twisting functors for sufficiently 
large right $H$-comodules have positive matrices with respect to standard 
bases. We will apply theorems of Perron and Frobenius by an argument similar 
to the one used in the proof of Proposition 2.2. First we mention a positivity 
result for endomorphisms of $G_0(R)$.

\proclaim
Lemma 7.2.
There exists a finite-dimensional right $H$-comodule $U$ such that the twisted 
$R$-module $U\ot V$ is faithful for each nonzero $V\in\MR,$ and therefore 
the endomorphism of $\,G_0(R)$ induced by the functor $\,U\ot\,?$ has 
a positive matrix.
\endproclaim

\Proof.
Denote by $I$ the annihilator in $R$ of a right $R$-module $V$. Since $I$ is 
an ideal of $R$, so too is the set
$$ 
I_C=\{a\in R\mid Ca\sbs I\}
$$
for each subcoalgebra $C\sbs H$. By formula (2.1) the annihilator of $U\ot V$ 
in $R$ coincides with the ideal $I_{S(C)}$ where $C$ is the smallest 
subcoalgebra of $H$ such that $U\in\MC$.

Denote by $\calF$ the set of all finite-dimensional subcoalgebras of $H$. The
intersection $\,\bigcap_{C\in\calF}I_{S(C)}\,$ is the largest $S(H)$-stable 
ideal of $R$ contained in $I$. It was proved in \cite{Sk11, Lemma 3.2} that 
each right artinian $H$-simple left $H$-module algebra is even $S(H)$-simple. 
Hence $\bigcap_{C\in\calF}I_{S(C)}=0$ provided that $V\ne0$, and therefore 
$I_{S(C)}=0$ for some $C\in\calF$ since $R$ is right artinian. 

Thus for each simple right $R$-module $V$ the twisted $R$-module $U\ot V$ is 
faithful when a finite-dimensional right $H$-comodule $U$ is large enough. 
Since there are only finitely many pairwise nonisomorphic simple right 
$R$-modules, we can find $U$ such that the faithfulness of $U\ot V$ holds for 
all simple right $R$-modules $V$ simultaneously, and this $U$ fulfills the 
desired conclusion. Indeed, each faithful right $R$-module has all simple 
right $R$-modules as its composition factors.
\endproof

\proclaim
Corollary 7.3.
For each finite-dimensional right $H$-comodule $U$ whose composition factors 
include a sufficiently large finite set of simple comodules there exists a 
positive function $f\in G_0(R)^*$ such that $\,fU=(\dim_kU)f$.
\endproclaim

\Proof.
By Lemma 7.2 there exists a finite-dimensional right $H$-comodule $U$ such 
that the endomorphism of $\,G_0(R)$ induced by the twisting functor 
$\,U\ot\,?$ has a positive matrix $M$. This remains true when $U$ is replaced 
by any other finite-dimensional right $H$-comodule which has the same or a 
larger set of composition factors.

The matrix of the dual endomorphism of $G_0(R)^*$ is the transpose of $M$, 
so it is positive too. By the Perron-Frobenius theorem the extension of the 
latter endomorphism to the real vector space $G_0(R)^*\ot_\bbZ\bbR$ admits a 
positive eigenvector $f$. The respective eigenvalue $\la$ is the largest real 
eigenvalue of $M$.

The class $[R]$ of cyclic free right $R$-modules is a positive element of 
$G_0(R)$. Since $U\ot R$ is a free right $R$-module of rank $d=\dim_kU$ by 
\cite{Sk-Oy06, Lemma 1.2(i)}, we have $\,U\ot[R]=d{\mskip2mu}[R]\,$ in 
$G_0(R)$. It follows that $\la=d$. Since $\la\in\bbZ$, the eigenvector 
$f$ can be found inside $G_0(R)^*$.
\endproof

\proclaim
Lemma 7.4.
Given two indecomposable injectives $E,E'\in\MscrC,$ let $V$ and $V'$ be their 
simple socles. Let $U$ be a finite-dimensional right $H$-comodule and\/ $U^*$ 
its left dual in the monoidal category $\MH$. Then $E'$ is isomorphic to a 
direct summand of $\,U\ot E\,$ if and only if\/ $V$ is isomorphic to a 
subfactor of $\,U^*\ot V'$.
\endproclaim

\Proof.
The comodules $E$ and $E'$ are injective hulls of their socles $V$ and $V'$ 
which are simple right $\scrC$-comodules. By Corollary 6.3 the twisted 
$\scrC$-comodule $\,U\ot E$ is injective in $\MscrC$. It follows that $E'$ is 
isomorphic to its direct summand if and only if $\,\MscrC(V',\,U\ot E)\ne0$, 
while $V$ is isomorphic to a subfactor of $\,U^*\ot V'\,$ if and only if 
$\,\MscrC(U^*\ot V',\,E)\ne0$. These two conditions are equivalent to each 
other, again by Corollary 6.3.
\endproof

\proclaim
Corollary 7.5.
Assume that conditions\/ {\rm(C1)} and\/ {\rm(C2)} are satisfied. Then the 
endomorphism of $K_0(\MscrC)$ induced by the functor $U\ot{}?$ has a positive 
matrix if and only if so does the endomorphism of $G_0(\MscrC)$ induced by the 
functor $U^*\ot{}?$.
\endproclaim

This follows immediately from Lemma 7.4.

\proclaim
Proposition 7.6.
Suppose that $N\in\HMscrC$ is an $R$-finite object such that
$$
\textstyle\scrC=\sum_{\,\xi\in\MscrC(N,\,\scrC)}\,\xi(N).\eqno(7.3)
$$
Then condition\/ {\rm(C1)} is satisfied and there is an equivalence relation 
$\sim$ on the standard basis $\,\Indinj\MscrC$ of the group $K_0(\MscrC)$ 
such that for two indecomposable injectives $E,\,E'\in\MscrC$ one has 
$\,[E]\sim[E']\,$ if and only if $E'$ is isomorphic to a direct summand of 
$\,U\ot E$ for some finite-dimensional right $H$-comodule $U$.
\endproclaim

\Proof.
By Lemma 7.1 each right $\scrC$-comodule embeds in a direct sum of copies of 
$\scrC$. Since $\scrC$ is a sum of $\MscrC$-epimorphic images of $N$, each 
right $\scrC$-comodule is isomorphic to a subfactor of a direct sum of copies 
of $N$, i.e., $N$ is a subgenerator of the category $\MscrC$. In particular, 
each simple right $\scrC$-comodule is isomorphic to a subfactor of $N$. Since 
$N$ is $R$-finite, it has an $\MscrC$-composition series. Each simple right 
$\scrC$-comodule is isomorphic to one of the finitely many factors of such a 
series. This verifies condition (C1). Moreover, the class $[N]$ of $N$ is a 
positive element of the Grothendieck group $G_0(\MscrC)$.

The assignment $E\mapsto\soc E$ gives a bijection of the set $\Indinj\MscrC$ 
onto $\Irr\MscrC$. Lemma 7.4 allows us to reformulate the conclusion of 
Proposition 7.6 in terms of the binary relation on $\Irr\MscrC$ defined by 
the set $\Ga$ of all pairs
$$
([V],[V'])\in\Irr\MscrC\times\Irr\MscrC
$$
where $V$ and $V'$ are two simple right $\scrC$-comodules for which there exists 
a finite-dimensional right $H$-comodule $U$ such that $V$ is isomorphic to a 
composition factor of $\,U^*\ot V'$.

Since the antipode of $H$ is bijective, the functor $U\mapsto U^*$ is an 
antiequivalence of the category of finite-dimensional right $H$-comodules. It 
follows that $([V],[V'])\in\Ga$ for two simple right $\scrC$-comodules $V,V'$ 
if and only if $V$ is isomorphic to a composition factor of $U\ot V'$ for some 
finite-dimensional right $H$-comodule $U$. Clearly the binary relation $\Ga$ 
is reflexive and transitive. To prove that $\Ga$ is an equivalence relation it 
remains to show that it is symmetric.

Since the set $\Irr\MscrC$ is finite, we can find a finite-dimensional right 
$H$-comodule $U$ such that for each pair $V,V'$ of simple right 
$\scrC$-comodules we have
$$
([V],[V'])\in\Ga\quad\text{if and only if}\quad
\text{$V$ is isomorphic to a subfactor of $U\ot V'$}.
$$
Taking a larger $U$, if necessary, we may assume that $U$ satisfies also the 
conclusion of Corollary 7.3.

Let $U$ be such a comodule. Since $N$ is an $H$-equivariant $\scrC$-comodule, 
by Lemma 6.4 there is an isomorphism $U\ot N\cong N^d$ in 
$\MscrC$ where $\,d=\dim_kU$. It follows that $U\ot[N]=d\,[N]$ in the 
Grothendieck group $G_0(\MscrC)$. On the other hand, taking the composite of a 
positive additive function $G_0(R)\to\bbZ$ satisfying the conclusion of 
Corollary 7.3 with the map (7.1) we obtain a positive function $f\in 
G_0(\MscrC)^*$ such that $fU=df$.

In other words, both the endomorphism $\al$ of the group $G_0(\MscrC)$ induced 
by the functor $U\!\ot{}?$ and the dual endomorphism of the group 
$G_0(\MscrC)^*$ admit positive eigenvectors. Exactly as in the proof of 
Proposition 2.2, this implies that $\Ga$ is indeed symmetric (by a theorem of 
Frobenius the matrix of $\al$, with respect to a suitable ordering of the 
basis elements, is block diagonal with positive diagonal blocks).
\endproof

\Remark.
Proposition 7.6 is valid even when $H$ is a Hopf algebra whose antipode is not 
bijective. Indeed, since the right artinian $H$-simple left $H$-module algebra 
$Q$ is $S^2(H)$-simple by \cite{Sk11, Lemma 3.2}, the conclusion of Lemma 7.2 
can be established for some comodule which is a left dual in the monoidal 
category $\MH$. It follows then that for each sufficiently large 
finite-dimensional right $H$-comodule $U$ the endomorphism of $G_0(R)$ induced 
by the functor $\,U^*\ot\,?$ has a positive matrix. This is what is needed to 
deduce that the set $\Ga$ in the proof of Proposition 7.6  defines an 
equivalence relation. It is not clear whether bijectivity of the antipode is 
necessary in Proposition 7.9.
\endremark

\proclaim
Corollary 7.7.
If in Proposition 7.6 the comodule $N$ is a simple object of the category 
$\HMscrC,$ then the set $\,\Indinj\MscrC$ is a single equivalence class.
In this case condition\/ {\rm(C2)} holds when there exists at least one nonzero 
$R$-finite injective right $\scrC$-comodule, and then the endomorphism of 
$K_0(\MscrC)$ induced by the functor $U\ot{}?$ has a positive matrix for each 
finite-dimensional right $H$-comodule $U$ whose composition factors include a 
sufficiently large finite set of simple comodules.
\endproclaim

\Proof.
Let $V$ be a simple right $\scrC$-comodule contained in the $\MscrC$-socle of 
$N$. The map $\mu:H\ot V\to N$ afforded by the action of $H$ on $N$ is a 
morphism in $\MscrC$ by Lemma 6.4. Its image is a nonzero $\HMscrC$-subobject 
of $N$. Since $N$ is simple in $\HMscrC$, we deduce that $\mu$ is surjective. 

As we have seen in the proof of Proposition 7.6, each simple right 
$\scrC$-comodule is isomorphic to a subfactor of $N$, and therefore to a 
subfactor of $H\ot V$. Since $H$ is the union of its finite-dimensional right 
coideals, each simple right $\scrC$-comodule is isomorphic to a subfactor of 
$U\ot V$ for some finite-dimensional right $H$-comodule $U$. Thus the 
equivalence class of $[V]$ is the whole set $\Irr\MscrC$. As in Proposition 7.6, 
we use here the binary relation on $\Irr\MscrC$ corresponding to the 
established equivalence relation on $\Indinj\MscrC$ under the natural 
bijection between the two sets.

Thus the set $\Indinj\MscrC$ is also a single equivalence class. Suppose 
that $E$ is a nonzero $R$-finite injective right $\scrC$-comodule. All direct 
summands of $E$ are $R$-finite. So we may assume that $E$ is an indecomposable 
comodule. By \cite{Sk-Oy06, Lemma 1.1} the twisted $\scrC$-comodule $U\ot E$ 
is $R$-finite for any finite-dimensional right $H$-comodule $U$. Each 
indecomposable injective right $\scrC$-comodule is a direct summand of 
$U\ot E$ for a suitable $U$, and therefore it is $R$-finite as well. 
\endproof

\proclaim
Lemma 7.8.
Under the hypothesis of Proposition 7.6 the coring\/ $\scrC$ has a maximal 
right coideal.
\endproclaim

\Proof.
Considering any $\MscrC$-composition series $0=N_0\sbs N_1\sbs\cdots\sbs N_s=N$, 
put
$$
\textstyle\scrC_i=\sum_{\,\xi\in\MscrC(N,\,\scrC)}\,\xi(N_i)
$$
for each $i$. We thus obtain a chain of right coideals 
$0=\scrC_0\sbs\scrC_1\cdots\sbs\scrC_s=\scrC$ of $\scrC$. Let $p$ be the 
smallest integer such that $\scrC_p=\scrC$. The quotient $\scrC\,/\,\scrC_{p-1}$ 
is a nonzero semisimple right $\scrC$-comodule, being a sum of homomorphic 
images of the simple comodule $N_p/N_{p-1}$. Therefore 
$\scrC\,/\,\scrC_{p-1}$ has a maximal subcomodule, and we can take its 
preimage in $\scrC$.
\endproof

For a right $\scrC$-comodule $V$ and a grouplike element $g\in\scrC$ the 
vector subspace of $g$-invariants of $V$ is defined by the formula
$$
V_g\co C=\{v\in V\mid\rho(v)=v\ot_Rg\}
$$
where $\rho:V\to V\ot_R\scrC$ is the comodule structure map.

\setitemsize(iii)
\proclaim
Proposition 7.9.
Suppose $N$ is an $R$-finite simple object of the category $\HMscrC$ such 
that equality\/ \refeq{7.3} is satisfied and there exists a nonzero map $N\to 
R$ which is simultaneously $R$-linear and $H$-linear, i.e., 
$\,\Hom_R(N,R)^H\ne0$.

Then both conditions\/ {\rm(C1)} and\/ {\rm(C2)} are satisfied. Furthermore, 
for each nonzero $R$-finite object $M\in\HMscrC$ the following is true:

\item(i)
\ the class $\,[M]\in G_0(\MscrC)\,$ is a rational multiple of 
$\,[N]\in G_0(\MscrC),$

\item(ii)
\ each simple right $\scrC$-comodule is isomorphic to an $\MscrC$-subobject of 
$M,$

\item(iii)
\ $M_g\co C\ne0\,$ for any grouplike element $g\in\scrC$.

\endproclaim

\Proof.
Each $R$-linear map $\eta:N\to R$ gives rise to a morphism $\eta':N\to\scrC$ 
in $\MscrC$ defined as the composite
$$
N\lmapr2\rho N\ot_R\scrC\lmapr4{\eta\ot\id}R\ot_R\scrC\cong\scrC
$$
where $\rho$ is the comodule structure map. If $\eta$ is $H$-linear, then so 
too is $\eta'$, in which case $\eta'$ is a morphism in $\HMscrC$. Furthermore, 
$\eta'\ne0$ provided that $\eta\ne0$, and then the simplicity of $N$ in 
$\HMscrC$ implies that $\eta'$ is injective, and so $N\cong\eta'(N)$. Thus the 
hypothesis of Proposition 7.9 allows us to assume that $N$ is an 
$\HMscrC$-subobject of $\scrC$, i.e., an $H$-invariant right coideal.

Since $\scrC$ is injective in $\MscrC$ by Lemma 7.1, each morphism $N\to\scrC$ 
in $\MscrC$ extends to an endomorphism of $\scrC$. Equality (7.3) yields
$$
\textstyle\scrC=\sum_{\,\xi\in\End^\scrC\!\scrC}\,\xi(N)\eqno(7.4)
$$
where $\End^\scrC\!\scrC$ is the $\MscrC$-endomorphism ring of $\scrC$. It 
follows from Lemma 7.8 that there exist a simple right $\scrC$-comodule $V$ and 
an epimorphism $\ph:\scrC\to V$ in $\MscrC$. By (7.4) $\ph\bigl(\xi(N)\bigr)\ne0$ 
for some $\xi\in\End^\scrC\scrC$.  Replacing $\ph$ with $\ph\circ\xi$, we may 
assume that the restriction $\ph|_N$ of $\ph$ to $N$ is nonzero.

We fix such an epimorphism $\ph:\scrC\to V$. Its kernel $J=\Ker\ph$ is a right 
coideal of $\scrC$ such that $N\not\sbs J$ and $\scrC\,/J\cong V$ in $\MscrC$. 
Since $V$ is simple, it is $R$-finite. We are going to find a right coideal 
$J'$ of $\scrC$ such that $J'\cap N=0$ and $\scrC\,/J'$ is an $R$-module of 
finite length.

For each right coideal $U$ of the Hopf algebra $H$ put
$$
J_U=\{c\in\scrC\mid hc\in J\text{ for all }h\in U\}.
$$
Suppose that $\,\dim_kU<\infty$. Since the antipode of $H$ is bijective, the 
dual vector space of $U$ has a comodule structure which makes it the right dual 
$\Urd$ of $U$ in the monoidal category $\MH$. Then $U$ is the left dual of 
$\Urd$. Corollary 6.3 shows that under the canonical bijection
$$
\Hom_k(U\ot\scrC,\,V)\cong\Hom_k(\scrC,\,\Urd\ot V)
$$
the $\MscrC$-morphisms $\,U\ot\scrC\to V$ correspond to the 
$\MscrC$-morphisms $\,\scrC\to\Urd\ot V$. Now $J_U$ is the kernel of 
the map $\,\nu:\scrC\to\Urd\ot V\,$ that corresponds to the composite
$$
U\ot\scrC\lmapr2\mu\scrC\lmapr2\ph V
$$
where $\mu$ arises from the action of $H$ on $\scrC$. Since $\mu$ is a 
morphism in $\MscrC$ by Lemma 6.4, so is $\ph\circ\mu$, and therefore also 
$\nu$. Hence $J_U$ is a right coideal of $\scrC$ and there is a monomorphism 
$\scrC\,/J_U\hookrightarrow\Urd\ot V$ in $\MscrC$. Note that $\Urd\ot V$ is an 
$R$-finite comodule since so is $V$ and $\Urd$ is finite-dimensional. It 
follows that $\scrC\,/J_U$ is an $R$-finite comodule as well.

Consider the family $\calF$ of all right coideals $J_U$ of $\scrC$ for various 
finite-dimensional right coideals $U$ of $H$. Since $N$ has finite length in 
$\MR$, and therefore also in $\MscrC$, there exists $J'\in\calF$ such that 
$J'\cap N$ is the smallest among all intersections $J_U\cap N$ for $J_U\in\calF$. 
Then
$$
J'\cap N=\{x\in N\mid Hx\sbs J\}.
$$
Hence $J'\cap N$ is an $H$-submodule of $N$ contained in $J\cap N\ne N$. Since 
$J'\cap N$ is also a $\scrC$-subcomodule, it is an $\HMscrC$-subobject of $N$. 
Since $N$ is simple in $\HMscrC$, we do get $J'\cap N=0$, and we know already 
that $\scrC\,/J'$ is $R$-finite since $J'\in\calF$.

Denote by $E$ the injective hull of $N$ in the category $\MscrC$. Since $\scrC$ 
is injective in $\MscrC$, the inclusion map $N\to\scrC$ extends to a monomorphism 
$E\to\scrC$ in $\MscrC$. It allows us to view $E$ as a right coideal of $\scrC$ 
as well. Since $E$ is an essential extension of $N$ in $\MscrC$, it follows that 
$J'\cap E=0$. The canonical map $E\to\scrC\,/J'$ is therefore injective, 
whence $E$ is $R$-finite. Now condition (C2) follows from Corollary 7.7, while 
(C1) is a consequence of Proposition 7.6.

In the final part of the proof we deal with statements (i), (ii), and (iii).

By Corollaries 7.5 and 7.7 there exists a finite-dimensional right $H$-comodule 
$U$ such that the endomorphism $\al$ of the Grothendieck group $G_0(\MscrC)$ 
induced by the functor $U\ot{}?$ has a positive matrix. Extend $\al$ by 
linearity to the real vector space $G_0(\MscrC)\ot_\bbZ\bbR$. By the 
Perron-Frobenius theorem $\al$ has precisely one, up to scalar multiplication, 
eigenvector with nonnegative coordinates. But $U\ot[M]=d\,[M]\,$ by Lemma 6.4. 
Thus $[M]$ is such an eigenvector of $\al$, and the same holds for $[N]$. 
Hence $[M]=c\,[N]$ for some $c\in\bbR$. Then clearly $c\in\bbQ$ since both 
$[M]$ and $[N]$ lie in the group $G_0(\MscrC)$.

Denote by $E$ the injective hull of $M$ in the category $\MscrC$. Since $M$ is 
$R$-finite, its $\MscrC$-socle is a direct sum of finitely many simple comodules, 
whence $E$ is a direct sum of finitely many indecomposable injective comodules. 
By (C2) $E$ is $R$-finite. By Corollary 7.7 there exists a finite-dimensional 
right $H$-comodule $U$ such that each indecomposable injective right 
$\scrC$-comodule is isomorphic to a direct summand of $U\ot E$, i.e., 
$U\ot E$ is an injective cogenerator in $\MscrC$. Taking a larger $U$, if 
necessary, we may assume that $U$ satisfies the conclusion of Corollary 7.3. 

Put $d=\dim_kU$. The $\scrC$-comodule $U\ot M\cong M^d$ embeds in $U\ot E$. 
Since $E^d$ is an injective hull of $M^d$ in $\MscrC$, there exists also an 
embedding of $E^d$ in $U\ot E$. Hence $U\ot E\cong E^d\oplus E'$ in $\MscrC$ 
for some comodule $E'$, and therefore 
$$
U\ot[E]=d\,[E]+[E']
$$
in $K_0(\MscrC)$. Applying map (7.2), we get such an equality also in $G_0(R)$, 
and then, applying the function $f$ given by Corollary 7.3, we deduce that 
$f([E'])=0$ since $f(U\ot[E])=d\mskip1mu f([E])$. It follows that $E'=0$ by 
positivity of $f$. Hence $U\ot E\cong E^d$. This implies that $E$ is a 
cogenerator in $\MscrC$ since so is $U\ot E$. Therefore each simple right 
$\scrC$-comodule embeds in $E$. Since $E$ is an essential extension of $M$ in 
$\MscrC$, all simple subcomodules of $E$ are contained in $M$. This verifies (ii).

The right $R$-submodule $gR\sbs\scrC$ generated by a grouplike $g\in\scrC$ is 
a right coideal of $\scrC$, and $g\in(gR)_g\co\scrC$. If $V$ is any simple 
factor comodule of the right $\scrC$-comodule $gR$, then the image of $g$ in 
$V$ is contained in the subspace of $g$-invariants $V_g\co\scrC$ and generates 
$V$ as an $R$-module. In particular, $V_g\co\scrC\ne0$. Since $V$ embeds in $M$ 
by (ii), we deduce (iii).
\endproof

We are interested in the special case of Proposition 7.9 when $N$ is a free 
$R$-module of rank 1 generated by an $H$-invariant element, i.e., $N\cong R$ 
in $\HMR$. The right $\scrC$-comodule structure on $R$ is given by the map
$$
\xi:\,R\to R\ot_R\scrC\cong\scrC,\qquad x\mapsto gx\,\text{ for $\,x\in R$},
\eqno(7.5)
$$
where $g$ is an $H$-invariant grouplike element of $\scrC$. Furthermore, 
$$
\MscrC(R,\scrC)\cong\Hom_R(R,R)\cong R.
$$
The identity element $1\in R$ corresponds to the $H$-linear $\MscrC$-morphism 
$\xi:R\to\scrC$ which coincides with the map defining comodule structure (7.5). 
This $\xi$ generates $\MscrC(R,\scrC)$ as a left $R$-module. Condition (7.3) 
implies that the assignment
$$
a\ot b\mapsto a\,\xi(b)=agb,\qquad a,b\in R,
$$
defines a surjective homomorphism of left $H$-module $R$-corings 
$\,\pi:R\ot R\to\scrC$. It allows us to identify $\scrC$ with a factor coring 
of the canonical coring $R\ot R$ associated with the ring extension $k\sbs R$. 
Note that $\,g=\pi(1\ot1)$. This is the \emph{distinguished grouplike} of the 
factor coring. Thus we are led to the following conclusion:

\proclaim
Corollary 7.10.
Let $\,\scrC=(R\ot R)/\,\scrI$ where $\scrI$ is an $H$-invariant coideal of 
the canonical $R$-coring $R\ot R$. Suppose that $\scrC$ is left $R$-flat and 
$R$ has no $H$-invariant right ideals except for the zero ideal and the whole 
$R$. Then $\,M\co\scrC\ne0\,$ for each nonzero $R$-finite object 
$\,M\in\HMscrC$.
\endproclaim

\Proof.
We use the distinguished grouplike $g\in\scrC$ to define a $\scrC$-comodule 
structure on $R$, as shown in (7.5). By the hypothesis $R$ is simple as an 
object of the category $\HMR$, and therefore simple as an object of $\HMscrC$. 
Hence Proposition 7.9 applies with $N=R$.
\endproof

\setitemsize(3)
\proclaim
Theorem 7.11.
Let $B$ be an $H$-module subalgebra of a right artinian left $H$-module 
algebra $Q$. Suppose that

\item(1)
\ $B$ is the sum of its finite-dimensional $H$-submodules,

\item(2)
\ each nonzero $H$-stable right ideal of $B$ contains a nonzerodivisor of $B,$

\item(3)
\ $Q$ is a classical two-sided quotient ring of $B$.

If\/ $\scrC$ is any left $H$-module factor coring of the canonical $Q$-coring 
$Q\ot Q,$ then $\,\scrC\cong Q\ot_AQ\,$ where $\,A=Q\co\scrC\!$.
\endproclaim

\Proof.
We repeat several earlier arguments. Conditions (2) and (3) imply that $Q$ has 
no $H$-invariant right ideals except for the zero ideal and the whole $Q$. In 
particular, $Q$ is an $H$-simple algebra.

Let $\pi:Q\ot Q\to\scrC$ be the canonical surjective homomorphism of left 
$H$-module $Q$-corings. Since the left coideal $\pi(Q\ot B)$ of $\scrC$ 
is a locally $Q$-finite object of the category $\HQM$, it is projective in 
$\QM$ by the $H\cop$-variant of Theorem 2.1. Since $Q$ is a classical right 
quotient ring of $B$, and so, in particular, $Q$ is left $B$-flat, we deduce 
that $\scrC\cong\pi(Q\ot B)\ot_BQ$ is left $Q$-flat. (Then $\scrC$ is 
projective in $\QM$ and even free by the $H\cop$-variants of Lemmas 2.10, 
2.11.) Thus $Q$ and $\scrC$ satisfy the hypothesis of Corollary 7.10, and 
therefore $M\co\scrC\ne0$ for each nonzero $Q$-finite object $M\in\HMscrC$.

The algebra $A$ consists of all elements $x\in Q$ such that 
$1\ot x-x\ot1\in\Ker\pi$ by the argument in the proof of Lemma 5.1. It 
follows that $A$ is a dominion subalgebra of $Q$. By the $H\cop$-variant of 
Proposition 2.13 $A$ is right artinian and $H$-simple. By the $H\cop$-variant 
of Proposition 2.12 $Q$ is a free left $A$-module.

The homomorphism $\pi$ factors through $\pi':Q\ot_AQ\to\scrC$. Arguing as in 
the proof of the implication (d)$\,\Rar\,$(a) of Lemma 5.4, we deduce that the 
restriction of $\pi'$ to $UA\ot_AQ$ is injective for each finite-dimensional 
$H$-submodule $U\sbs Q$. By condition (1) in the hypothesis this implies that 
the restriction of $\pi'$ to $BA\ot_AQ$ is injective. Since $Q$ is a classical 
left quotient ring of $B$, each nonzero left $Q$-submodule of $Q\ot_AQ$ has a 
nonzero intersection with the $(B,Q)$-subbimodule $BA\ot_AQ$. It follows that 
$\Ker\pi'=0$, i.e., $\pi'$ is bijective. Thus $\pi'$ is an isomorphism of left 
$H$-module $Q$-corings.
\endproof

When $H=k$ is the trivial one-dimensional Hopf algebra, an algebra 
$Q$ without nontrivial $H$-invariant right ideals is just a skew field, and 
conditions (1)--(3) in the hypothesis of Theorem 7.11 are satisfied with $B=Q$. 
In this case Theorem 7.11 reduces to Sweedler's Fundamental Theorem 
\cite{Sw75} which describes all factor corings of the canonical coring 
$\,Q\ot Q$.

Suppose now that $Q$ is a finite-dimensional central simple algebra. Since the 
algebra $Q\ot Q\op$ is isomorphic to $\End_kQ$, the category of $Q$-bimodules 
$\QMQ$ is equivalent to the category $\Mk$ of vector spaces. Up to isomorphism, 
$Q$ is the only simple bimodule. An equivalence is given by the functor 
$Z:\QMQ\to\Mk$ such that
$$
Z(M)=\{\,x\in M\mid\,ax=xa\,\text{ for all $a\in Q$}\,\}
$$
for each $Q$-bimodule $M$. Moreover, since $Z(M\ot_QN)\cong Z(M)\ot Z(N)$, 
naturally in $M,N\in\QMQ$, this equivalence is monoidal. As an immediate 
consequence we get

\proclaim
Proposition 7.12.
For any finite-dimensional central simple algebra $Q$ the functor $Z$ induces 
an equivalence between the category of $Q$-corings and the category of 
coalgebras over the base field.
\endproclaim

Denote by $X(M)$ the vector space of all $Q$-bimodule homomorphisms $M\to Q$. 
Then $X(M)\cong Z(M)^*$, and it follows that the canonical map 
$$
X(M)\ot X(N)\to X(M\ot_QN)
$$
is bijective whenever $M$ and $N$ are finitely generated $Q$-bimodules.

\proclaim
Corollary 7.13.
For any finite-dimensional central simple algebra $Q$ the functor $X$ induces 
an antiequivalence between the category of $Q$-corings which are finitely 
generated $Q$-bimodules and the category of finite-dimensional algebras over 
the base field. Under this antiequivalence the canonical coring $\,Q\ot Q$ 
corresponds to $Q,$ and therefore factor corings of $\,Q\ot Q$ correspond 
to subalgebras of $Q$.
\endproclaim

\Proof.
For each finitely generated $Q$-coring $\scrC$ the algebra $X(\scrC)$ is the 
dual of the coalgebra $Z(\scrC)$. Let $\scrC=Q\ot Q$. In this case for each 
element $a\in Q$ there is a unique bimodule homomorphism $f_a:\scrC\to Q$ 
sending $1\ot1$ to $a$. It is given by the formula
$$
f_a(x\ot y)=xay,\qquad x,y\in Q.\eqno(7.6)
$$
Under the convolution product
$$
\scrC\lmapr2{\De_\scrC}\scrC\ot_Q\scrC\lmapr4{f_a\ot f_b}Q\ot_QQ\cong Q
$$
the element $1\ot1\in\scrC$ goes to $ab\in Q$, i.e., $f_a*f_b=f_{ab}$ for 
all $a,b\in Q$. Also, $f_1=\ep_\scrC$ is the identity element of $X(\scrC)$. 
Thus the assignment $a\mapsto f_a$ gives an isomorphism of algebras 
$\,Q\cong X(Q\ot Q)$.
\endproof

By Corollary 7.13 each factor coring $\scrC$ of the $Q$-coring $Q\ot Q\,$ is 
reconstructed from the corresponding subalgebra $\,X(\scrC)\sbs Q\,$ 
as $\,\scrC=(Q\ot Q)/\,\scrI\,$ where
$$
\scrI=\{\,t\in Q\ot Q\mid\,f_a(t)=0\,\text{ for all $\,a\in X(\scrC)$}\}
$$
and $\,f_a:\,Q\ot Q\to Q\,$ is defined by formula (7.6).

If $\,\scrC=Q\ot_AQ\,$ where $A$ is some subalgebra of $Q$, then $X(\scrC)$ 
coincides with the centralizer of $A$ in $Q$. Incidentally, a subalgebra of a 
central simple algebra is a dominion subalgebra if and only if it coincides 
with its double centralizer (see \cite{Scho, Lemma 7.16}). Thus, under the 
antiequivalence of Corollary 7.13 factor corings of the form $Q\ot_AQ$ 
correspond to dominion subalgebras of $Q$. If $Q$ is not a skew field, there 
may exist subalgebras of $Q$ which are not centralizers, and then not all 
factor corings of $Q\ot Q$ can be written as $Q\ot_AQ$.

Already in the case $H=k$ we see that the condition on $H$-stable right ideals 
used in the hypothesis of Theorem 7.11 cannot be weakened to a similar 
condition only on $H$-stable two-sided ideals.

\section
8. Equivalences between module and comodule categories

Starting from this section and up to the end of the paper we assume again that 
$H$ is a Hopf algebra satisfying assumptions (A1) and (A2). Let $C$ be a left 
$H$-module factor coalgebra of $H$. We are going to describe the categories 
$\MC$ and $\CM$ of right and left $C$-comodules in terms of the corresponding 
left $\Hd$-invariant artinian subalgebra $A$ of the quotient ring $Q=Q(H)$. 
Let $\scrC=Q\ot_AQ$ be the corresponding left $\Hd$-module factor coring of 
the canonical $Q$-coring $Q\ot Q$. By Corollary 4.11 $C$ is canonically 
isomorphic to the left $H$-module coalgebra $\,\scrC^{\Hd}\!$.

An object $W\in\HdMA$ will be called \emph{rationally extendible} if 
the $\Hd$-equivariant $Q$-module $W\ot_AQ$ is rationally generated, as defined 
in section 3. We will denote by $\HdxratMA$ the full subcategory of rationally 
extendible objects of $\HdMA$.

If $W\in\HdMA$ is \emph{rationally generated} in the sense that $W=\Rat(W)A$, 
then so is $W\ot_AQ$. Therefore the category $\HdxratMA$ contains all 
rationally generated objects of $\HdMA$. However, subobjects of a rationally 
generated object are not necessarily rationally generated.

The full subcategory $\HdxratAM$ of rationally extendible objects of $\HdAM$ 
is defined similarly.

\proclaim
Proposition 8.1.
For each exact sequence $0\to W'\to W\to W''\to0$ in $\HdMA$ the objects $W'$ 
and $W''$ are rationally extendible whenever so is $W$. In particular, the 
category $\HdxratMA$ is abelian and the inclusion functor $\,\HdxratMA\to\HdMA\,$ 
is exact.
\endproclaim

\Proof.
By Proposition 3.6 $Q$ is a free $A$-module. Applying the functor $?{}\ot_AQ$ 
we obtain then an exact sequence
$$
0\to W'\ot_AQ\to W\ot_AQ\to W''\ot_AQ\to0
$$
in $\HdMQ$. If $W\ot_AQ$ is rationally generated, then so are the other terms 
of that sequence by Corollary 3.5, whence the desired conclusion.
\endproof

\proclaim
Theorem 8.2.
There are quasi-inverse functors
$$
\displaylines{
\Phi:\HdxratMA\to\MC\quad\text{and}\quad\Psi:\MC\to\HdxratMA
\quad\text{such that}\cr
\noalign{\smallskip}
\Phi(W)=(W\ot_AQ)^{\Hd},\qquad\Psi(V)=(V\ot Q)\co\scrC
\hphantom{such that}
}
$$
for $\,W\in\HdxratMA$ and $\,V\in\MC\!$. Similarly, there are quasi-inverse 
functors 
$$
\displaylines{
\Phi':\HdxratAM\to\CM\quad\text{and}\quad\Psi':\CM\to\HdxratAM
\quad\text{such that}\cr
\noalign{\smallskip}
\Phi'(W)=(Q\ot_AW)^{\Hd},\qquad\Psi'(V)=\lco\scrC(Q\ot V)
\hphantom{such that}
}
$$
for $\,W\in\HdxratAM$ and $\,V\in\CM$. Thus 
$$
\MC\approx\HdxratMA\qquad\text{and}\qquad\CM\approx\HdxratAM.
$$
\endproclaim

\Proof.
By Lemma 5.4 the functors in (5.7) given by the assignments 
$\,W\mapsto W\ot_AQ\,$ and $M\mapsto M\co\scrC$ are quasi-inverse equivalences. 
Since $W\in\HdxratMA$ if and only if $W\ot_AQ\in\HdratMQ$, these functors 
induce a pair of quasi-inverse functors between the categories $\,\HdxratMA\,$ 
and $\,\HdratMC$.

By Proposition 4.12 the functor $\HdratMC\to\MC$ such that $M\mapsto M^{\Hd}\!$ 
is also an equivalence of categories. Here the right $\scrC^{\Hd}\!$-comodule 
$M^{\Hd}$ is viewed as a right $C$-comodule by means of the canonical 
isomorphism $C\cong\scrC^{\Hd}$. The quasi-inverse functor $\MC\to\HdratMC$ is 
given by the assignment $V\mapsto V\ot Q$ with the right $\scrC$-comodule 
structure
$$
V\ot Q\longrightarrow(V\ot Q)\ot_Q\scrC\cong V\ot\scrC
$$
obtained as the $Q$-linear extension of the composite 
$V\to V\ot C\hrar V\ot\scrC$ where the first map is given by the $C$-comodule 
structure on $V$, and the second map arises from the isomorphism 
$\,C\cong\scrC^{\Hd}\!$. Now $\Phi$ is the composite 
$$
\HdxratMA\longrightarrow\HdratMC\longrightarrow\MC
$$
of the functors just described, while $\Psi$ is the composite of quasi-inverse 
functors.

Since $Q$ is a free $A$-module on both sides, the faithfully 
flat descent shows also that there are quasi-inverse functors $\AM\to\scrCM$ 
and $\scrCM\to\AM$ given by the assignments $W\mapsto Q\ot_AW$ and 
$M\mapsto\lco\scrC M$. They induce quasi-inverse functors between the 
categories $\HdAM$ and $\HdCM$, and between the categories $\HdxratAM$ and 
$\HdratCM$.

By Proposition 4.12 the functor $\,\HdratCM\to\CM$, $\,M\mapsto M^{\Hd}\!$, 
has a quasi-inverse functor $\,\CM\to\HdratCM$, $\,V\mapsto Q\ot V$. The 
functors $\Phi'$ and $\Psi'$ are obtained again by composing the functors just 
described.
\endproof

\proclaim
Corollary 8.3.
Each object $W\in\HdxratMA$ is the union of its $A$-finite subobjects. 
Moreover, $W$ is a free $A$-module of rank equal to the vector space dimension 
of the corresponding $C$-comodule $\Phi(W)$.
\endproclaim

\Proof.
By the construction of $\Phi$ and $\Psi$ in Theorem 8.2 there is an isomorphism
$$
W\ot_AQ\cong\Phi(W)\ot Q\quad\text{in $\,\HdratMC$}.
$$
This object is a free $Q$-module of rank equal to the dimension of $\Phi(W)$. 
Since $Q$ is a faithfully flat ring extension of $A$, the $A$-module $W$ is 
finitely generated if and only if so is the $Q$-module $W\ot_AQ$, and this 
holds precisely when $\,\dim\Phi(W)<\infty$. By Theorem 8.2 the 
$\HdxratMA\,$-subobjects of $W$ are in a bijective correspondence with the 
subcomodules of $\Phi(W)$. It follows that the $A$-finite subobjects correspond 
to the finite-dimensional subcomodules, and the union of $A$-finite subobjects 
gives the whole $W$ since the union of finite-dimensional subcomodules gives 
the whole $\Phi(W)$. Since $W$ is locally $A$-finite, its freeness as an 
$A$-module follows from Theorem 2.1 and Lemma 2.14.
\endproof

The image $1_C\in C$ of the identity element $1\in H$ is a grouplike element 
which we call the \emph{distinguished grouplike} of $C$. It generates $C$ as a 
left $H$-module and spans a one-dimensional subcoalgebra $k1_C\sbs C$. The 
canonical surjection $H\to C$ is given by the assignment $h\mapsto h1_C$ for 
$h\in H$. Note that the distinguished grouplike of $\scrC$ is the image of 
$1_C$ under the canonical isomorphism $\,C\to\scrC^{\Hd}\!$.

\proclaim
Lemma 8.4.
We have $\Phi(HA)\cong C$ and $\Phi(A)\cong k1_C$ in $\MC$. Also,
$\Phi'(AH)\cong C$ and $\Phi'(A)\cong k1_C$ in $\CM$.
\endproclaim

\Proof.
At the first step we form the object $M=HA\ot_AQ$ of the category $\HdratMC\!$. 
Since $HA\sbs Q$, this object is identified with an $\Hd$-invariant right 
coideal of the coring $\scrC=Q\ot_AQ$. Furthermore, $QM=\scrC$ since $1\in 
HA$. By Lemma 4.1 $\scrC^{\Hd}\!\sbs M$, and therefore 
$\,\scrC^{\Hd}\!\!=M^{\Hd}\!$. The right $\scrC^{\Hd}\!$-comodule structure 
on $M^{\Hd}$ is provided by the comultiplication of $\scrC^{\Hd}\!$, and the 
right $C$-comodule structure on $\,\Phi(HA)=M^{\Hd}$ arises from the coalgebra 
isomorphism $\,\psi:C\to\scrC^{\Hd}$ described in Corollary 4.11. Considering 
$C$ as a right $C$-comodule, the map $\,\psi:C\to\Phi(HA)\,$ is an isomorphism 
in $\MC$.

The object $A\ot_AQ\cong Q\in\HdratMC$ is identified with the $\Hd$-invariant 
right coideal $gQ$ of the coring $\scrC$ where $\,g=1\ot_A1\,$ is the 
distinguished grouplike of $\scrC$. Since $Q^{\Hd}\!\!=k$, we have 
$\Phi(A)=(gQ)^{\Hd}\!\!=kg$, which is the image of $k1_C$ under $\psi$. This 
yields the evaluation of $\Phi$ at $A$.

In the case of the functor $\Phi'$ we proceed similarly, working now with the 
left coideals $\,Q\ot_AAH$ and $\,Qg\,$ of the coring $\scrC$.
\endproof

\proclaim
Corollary 8.5.
An object $W\in\HdMA$ is rationally extendible if and only if there exists a 
monomorphism $W\to W'$ in $\HdMA$ where $W'$ is rationally generated.
\endproclaim

\Proof.
It is well-known that injective $C$-comodules are direct summands of direct 
sums of copies of $C$. It follows then from Theorem 8.2 and Lemma 8.4 that 
injective objects of the category $\HdxratMA$ are direct summands of direct 
sums of copies of $HA$. Since $HA$ is rationally generated, so are all 
injective objects of $\HdxratMA$. Furthermore, an arbitrary object of 
$\HdxratMA$ embeds in an injective one.

This shows that each object of $\HdxratMA$ embeds in a rationally generated 
object. Conversely, if $W'\in\HdMA$ is rationally generated, then 
$W'$ is rationally extendible, whence so are all subobjects of $W'$ by 
Proposition 8.1.
\endproof

The next lemma provides another explicit evaluation of the functor $\Phi$. For 
each right $H$-comodule $U$ we regard $U\ot A$ as a rationally generated object 
of $\HdMA$ with respect to the action of $A$ on the second tensorand by right 
multiplications and the tensor product of $\Hd$-module structures on $U$ and 
$A$. The assignment $u\mapsto u\ot1$ defines an isomorphism of $U$ onto a 
rational $\Hd$-submodule of this object which generates $\,U\ot A\,$ freely as 
a right $A$-module. In a similar way $A\ot U\,$ is an object of the category 
$\HdAM$.

\proclaim
Lemma 8.6.
For objects $U\in\MH$ there are natural isomorphisms 
$\,\Phi(U\ot A)\cong\nobreak U$ 
in $\MC$ where each right $H$-comodule $U$ is regarded as a $C$-comodule 
with respect to the canonical coalgebra homomorphism $H\to C$.
\endproclaim

\Proof.
Put $W=U\ot A$ and $M=W\ot_AQ\cong U\ot Q$. Then $\Phi(W)=M^{\Hd}\!$. As in 
Lemma 4.9, we have a $k$-linear map $\,\ph:U\to M^{\Hd}\!$ defined by the rule
$$
\ph(u)=\sum u\0\ot S(u\1),\qquad u\in U.
$$
Since the $\MQ$-endomorphism of $M$ defined by the assignment
$$
u\ot q\mapsto\ph(u)q=\sum u\0\ot S(u\1)q,\qquad u\in U,\ q\in Q,
$$
is invertible, we deduce that $\ph$ is injective and its image $\ph(U)$ 
generates $M$ as a right $Q$-module freely, i.e., $M\cong\ph(U)\ot Q$. Since 
$Q^{\Hd}\!\!=k$, it follows that $\ph$ maps $U$ bijectively onto $M^{\Hd}\!$. 

For convenience let us identify $U$ and $W$ with their canonical images in $M$. 
This allows us to write elements omitting the tensor product signs. With this 
convention $\ph(u)=\sum u\0\,S(u\1)$, and the right $\scrC$-comodule structure 
$\,\rho:M\to M\ot_Q\scrC\,$ is defined by the formula
$$
\rho(xq)=x\ot_Qg{\mskip1mu}q,\qquad x\in W,\ q\in Q,
$$
where $g$ is the distinguished grouplike of $\scrC$. Since $U\sbs W$, we get
$$
\eqalign{
\rho\bigl(\ph(u)\bigr)
=\sum u\0\ot_QgS(u\1)
&=\sum u\0S(u\1)\ot_Qu\2gS(u\3)\cr
&=\sum\,\ph(u\0)\ot_Q\psi(u\1)
}
$$
for all $u\in U$ where $\psi:H\to\scrC^{\Hd}\!$ is the homomorphism of left 
$H$-module coalgebras defined by formula (4.8). The map $\rho$ induces the 
right $\scrC^{\Hd}\!$-comodule structure
$$
M^{\Hd}\!\to M^{\Hd}\!\ot\,\scrC^{\Hd},\qquad
\ph(u)\mapsto\sum\,\ph(u\0)\ot\psi(u\1)\quad\text{for $u\in U$}.
$$
Since the canonical isomorphism $C\cong\scrC^{\Hd}\!$ is induced by $\psi$, 
it is clear that under the bijection $\ph$ the above comodule structure on 
$M^{\Hd}\!$ corresponds precisely to the right $C$-comodule structure on $U$ 
defined by the composite map $\,U\to U\ot H\to U\ot C$, and we are done.  
\endproof

Next we will determine the largest rational $\Hd$-submodules of the objects in 
the categories $\HdxratMA$ and $\HdxratAM$ in terms of the cotensor product 
$\sqC$. Recall that there are functors $\,?{}\sqCH:\MC\to\MH\,$ and 
$\,H\sqC{}?:\CM\to\vphantom{\calM}^H\!\calM\,$ (see \cite{Br-W}). 
The coaction invariants of $C$-comodules are defined by means of the 
distinguished grouplike $1_C$.

\proclaim
Proposition 8.7.
We have $\,\Rat\bigl(\Psi(V)\bigr)=V\sqCH\sbs V\ot H\,$ for each $V\in\MC$ and 
$$
\Rat\bigl(\Psi'(V)\bigr)=(S\ot\id)(H\sqC V)\sbs H\ot V\quad
\text{for each $V\in\CM$}.
$$
In particular, $\ \Rat(A)=\lco CH=S(H\co C)$.
\endproclaim

\Proof.
Let $M=V\ot Q\in\HdratMC$ where $V\in\MC$. Then $\Psi(V)=M\co\scrC$. Since 
$\Hd$ acts trivially on $V$ and $\Rat(Q)=H$ by Corollary 3.4, we have 
$\Rat(M)=V\ot H$. Let $\psi:C\to\scrC^{\Hd}\!$ be the canonical isomorphism. 
The action of $H$ on $C$ corresponds under $\psi$ to the action $\,\triangleright\,$ 
on $\,\scrC^{\Hd}\!$ given by formula (4.1). The distinguished grouplike of the 
coring $\scrC$ is the element $g=\psi(1_C)$. Recall that $M\co\scrC$ is the 
equalizer of the pair of right $A$-linear maps
$$
\rho,\,\tau:\,M\longrightarrow M\ot_Q\scrC\cong V\ot\scrC
$$
where $\rho$ is the comodule structure on $M$ and $\tau$ is defined by the 
rule $\tau(x)=x\ot_Qg$ for $x\in M$. Denote by $\ka$ the composite
$$
V\ot C\ot H\lmapr7{\id\ot\psi\ot\id}V\ot\scrC^{\Hd}\!\!\ot H
\lmapr4{\id\ot\mu}V\ot\scrC
$$
where $\mu:\scrC^{\Hd}\!\!\ot H\to\scrC$ is afforded by the right $Q$-module 
structure on $\scrC$. Recall that $\Rat(\scrC)\in\MHH$ and $\mu$ is a 
$k$-linear bijection of $\,\scrC^{\Hd}\!\!\ot H$ onto $\Rat(\scrC)$ by the 
fundamental theorem on Hopf modules. It follows that $\ka$ is injective. 
If $v\in V$ and $h\in H$, then
$$
\eqalign{
\rho(v\ot h)&=\sum v\0\ot\psi(v\1)h=\sum\ka(v\0\ot v\1\ot h),\cr
\tau(v\ot h)&=v\ot hg=\sum\ka(v\ot h\11_C\ot h\2)
}
$$
since $\,hg=\sum\,(h\1\triangleright g)\,h\2=\sum\psi(h\11_C)\,h\2\,$.

\smallskip
Hence $\,\Rat\bigl(\Psi(V)\bigr)=\Rat(M)\cap M\co\scrC\,$ coincides with the 
equalizer of the two $k$-linear maps $\,V\ot H\to V\ot C\ot H\,$ defined by the 
assignments
$$
v\ot h\mapsto\sum v\0\ot v\1\ot h\qquad\text{and}\qquad
v\ot h\mapsto\sum v\0\ot h\11_C\ot h\2\,.
$$
This equalizer is nothing else but precisely $V\sqCH$ since the map $H\to C\ot H$ 
defined by the rule $h\to\sum h\11_C\ot h\2$ is the left $C$-comodule structure 
on $H$.

In the second case where $V\in\CM$ we put $M=Q\ot V\in\HdratCM$ and proceed 
similarly. The left $A$-module $\,\Psi'(V)=\lco\scrC M\,$ is the equalizer 
of the pair of maps
$$
\la,\,\tau:\,M\longrightarrow\scrC\ot_QM\cong\scrC\ot V
$$
where $\la$ is the comodule structure on $M$ and $\tau$ is defined by the 
rule $\tau(x)=g\ot_Qx$ for $x\in M$. If $v\in V$ and $h\in H$, then
$$
\eqalign{
\la(h\ot v)&=\sum\,h\,\psi(v\ng)\ot v\0,\cr
\tau(h\ot v)&=gh\ot v
=\sum\,h\2\bigl(S^{-1}(h\1)\triangleright g\bigr)\ot v
=\sum\,h\2\,\psi\bigl(S^{-1}(h\1)1_C\bigr)\ot v.
}
$$
Now $\Rat(M)=H\ot V$, and $\,\Rat\bigl(\Psi'(V)\bigr)=\Rat(M)\cap\lco\scrC M\,$ 
coincides with the equalizer of the two $k$-linear maps 
$\,H\ot V\to H\ot C\ot V$ defined by the assignments
$$
h\ot v\mapsto\sum h\ot v\ng\ot v\0\qquad\text{and}\qquad
h\ot v\mapsto\sum h\2\ot S^{-1}(h\1)1_C\ot v.
$$
The invertible linear endomorphism $\,S^{-1}\ot\id\,$ of the vector space 
$\,H\ot V\,$ transforms this equalizer into the equalizer of the first of 
the two previous maps and another linear map defined by the assignment
$$
h\ot v\mapsto\sum h\1\ot h\21_C\ot v.
$$
This last equalizer is precisely the cotensor product $\,H\sqC V\sbs H\ot V$, 
and we get the claimed description of the rational $\Hd$-module 
$\,\Rat\bigl(\Psi'(V)\bigr)$.

For the final conclusion take $V=k1_C$ regarded as either right or left 
$C$-comodule. In this case, respectively, $V\ot Q\cong Q$ or $Q\ot V\cong Q$ 
with the $\scrC$-comodule structure obtained by means of the distinguished 
grouplike $g\in\scrC$. Hence $\,\Psi(V)\cong Q\co\scrC=A\,$ or 
$\,\Psi'(V)\cong\lco\scrC Q=A$. As we have proved, 
$$
\eqalign{
\Rat\bigl(\Psi(k1_C)\bigr)&=k1_C\sqCH\cong\lco CH,\cr
\noalign{\smallskip}
\Rat\bigl(\Psi'(k1_C)\bigr)&=(S\ot\id)(H\sqC k1_C)\cong S(H\co C),
}
$$
which yields the two required expressions for the rational part $\,\Rat(A)\,$ 
of $A$.
\endproof

Given a left $\Hd$-module algebra $R$, the objects $M\in\HdMC$ equipped with a 
left $R$-module structure with respect to which $M$ is an object of the category 
$\HdRMQ$ and the comodule structure map $M\to M\ot_Q\scrC$ is $R$-linear form 
a category which we denote by $\HdRMC$.

The objects of another category $\HdCMR$ are $\Hd$-equivariant left 
$\scrC$-comodules equipped with a right $R$-module structure satisfying 
similar conditions.

In this section of the paper the categories $\HdRMC$ and $\HdCMR$ will be used 
for $R=H$. Another case where $R=A$ will appear in section 10.

\proclaim
Lemma 8.8.
For any left $\Hd$-module algebra $R$ the functors {\rm(5.7)} induce a 
category equivalence $\,\HdRMA\approx\HdRMC$. Similarly, $\,\HdAMR\approx\HdCMR$.  
\endproclaim

\Proof.
If $M$ is an object of $\HdRMC$, then the two maps in (5.2) are left $R$-linear, 
as well as $\Hd$-linear and $A$-linear. Their equalizer $M\co\scrC$ is an 
$R$-submodule of $M$, and therefore $M\co\scrC\in\HdRMA$. 

Conversely, given an object $W\in\HdRMA$, its extension $W\ot_AQ$ is 
an object of $\HdRMQ$, and the comodule structure map (5.6) is $R$-linear.
\endproof

Denote by $\HdratHMQ$ and $\HdratHMC$ the full subcategories, respectively, of 
the categories $\HdHMQ$ and $\HdHMC$ consisting of those objects that are 
rationally generated as objects of the category $\HdMQ$. 

The full subcategories $\,\HdratQMH\sbs\HdQMH\,$ and $\,\HdratCMH\sbs\HdCMH\,$ 
are defined similarly by the rational generation in $\HdQM$.

\proclaim
Lemma 8.9.
There are four category equivalences
$$
\vcenter{\halign{\hfil$\displaystyle#$\hfil&\qquad\qquad\hfil$\displaystyle#$\hfil\cr
\HdratHMQ\to\HM\,,&\HdratHMC\to\HMC,\cr
\noalign{\smallskip}
\noalign{\smallskip}
\HdratQMH\to\HM\,,&\HdratCMH\to\HCM\,,\cr
}}
$$
each of which is defined by the assignment $\,M\mapsto M^{\Hd}\!$.
\endproclaim

\Proof.
If $M$ is an object of either $\HdratHMQ$ or $\HdratQMH$, then its largest 
rational $\Hd$-submodule $\Rat(M)$ is an object of $\HMHH$. Exactly as in the 
proof of Proposition 4.2 this implies that the subspace of $\Hd$-invariants 
$V=M^{\Hd}\!$ is a left $H$-module with respect to the action defined by 
formula (4.1).

If $M\in\HdratHMQ$, then $M\cong V\ot Q$ by Proposition 3.2, and the 
corresponding left action of $H$ on $V\ot Q$ is as follows:
$$
h\cdot(v\ot q)=\sum\,(h\1\triangleright v)\ot h\2q,\qquad 
h\in H,\ v\in V,\ q\in Q.
$$
Conversely, given any left $H$-module $V$, the action of $H$ defined by this 
formula makes $V\ot Q$ an object of $\HdratHMQ$ such that 
$V\cong(V\ot Q)^{\Hd}\!$ and the action (4.1) recovers the original $H$-module 
structure on $V$. In this way we obtain a quasi-inverse of the functor 
$\HdratHMQ\to\HM$.

If $M\in\HdratQMH$ and $V=M^{\Hd}\!$, then $M\cong Q\ot V$ with the right 
action of $H$
$$
(q\ot v)\cdot h=\sum qh\2\ot\bigl(S^{-1}(h\1)\triangleright v\bigr),\qquad 
h\in H,\ v\in V,\ q\in Q.
$$
Using this formula for an arbitrary left $H$-module $V$, we obtain a 
quasi-inverse of the functor $\HdratQMH\to\HM$.

By Proposition 4.12 the right $\scrC$-comodule structure $\rho:M\to M\ot_Q\scrC$ 
of an object $M\in\HdratMC$ corresponds to a right $C$-comodule structure 
$\rho_V:V\to V\ot C$ on the subspace $V=M^{\Hd}\!$. These two structures are 
related by means of the identity
$$
\rho(vq)=\sum v\0\ot_Q\psi(v_1)\,q,\qquad v\in V,\ q\in Q.
$$
where $\psi:C\to\scrC^H$ is the canonical isomorphism. The property of $\psi$ 
being $H$-linear is expressed as
$$
\psi(hc)=h\triangleright\psi(c)=\sum\,h\1\,\psi(c)\,S(h\2),\qquad 
h\in H,\ c\in C.
$$
Suppose that $M$ also has a left $H$-module structure which makes $M$ into an 
object of $\HdratHMQ$. Then $V$ is a left $H$-module with respect to the action 
(4.1). If $\rho$ is left $H$-linear, then
$$
\eqalign{
\rho(h\triangleright v)
=\sum\,h\1\,\rho(v)\,S(h\2)
&=\sum\,h\1v\0\ot_Q\psi(v_1)S(h\2)\cr
&=\sum\,h\1v\0S(h\2)\ot_Qh\3\psi(v_1)S(h\4)\cr
&=\sum\,(h\1\triangleright v\0)\ot_Q\psi(h\2v_1),
}
$$
whence
$$
\rho_V(h\triangleright v)=\sum\,(h\1\triangleright v\0)\ot h\2v_1
$$
for all $h\in H$ and $v\in V$. This is the compatibility condition required 
for objects of the category $\HMC$. Conversely, if $\rho_V$ satisfies this 
identity, then
$$
\eqalign{
\rho(hvq)=\sum\,\rho\bigl((h\1\triangleright v)h\2q\bigr)
&=\sum\,(h\1\triangleright v\0)\ot_Q\bigl(h\2\triangleright\psi(v_1)\bigr)h\3q\cr
&=\sum\,h\1v\0S(h\2)\ot_Qh\3\psi(v_1)q\cr
&=\sum\,hv\0\ot_Q\psi(v_1)q=h\,\rho(vq)
}
$$
for all $h\in H$, $v\in V$, and $q\in Q$, i.e., $\rho$ is left $H$-linear. 
Thus $M\in\HdratHMC$ if and only if $V\in\HMC$.

Similarly, the $\scrC$-comodule structure $\la:M\to\scrC\ot_QM$ of an object 
$M\in\HdratCM$ corresponds to a $C$-comodule structure $\la_V:V\to C\ot V$ such 
that
$$
\la(qv)=\sum\,q\,\psi(v\ng)\ot_Qv\0,\qquad v\in V,\ q\in Q.
$$
Suppose that $M$ is also an object of $\HdratQMH$. If $\la$ is right $H$-linear, 
then
$$
\la(h\triangleright v)=\sum\,h\1\,\la(v)\,S(h\2)
=\sum\psi(h\1v\ng)\ot_Q(h\2\triangleright v\0),
$$
whence
$$
\la_V(h\triangleright v)=\sum\,h\1v\ng\ot(h\2\triangleright v\0)
$$
for all $h\in H$ and $v\in V$. Conversely, this identity implies that
$$
\eqalign{
\la(qvh)&=\sum\,\la\bigl(qh\2(S^{-1}(h\1)\triangleright v)\bigr)\cr
&=\sum\,qh\3\bigl(S^{-1}(h\2)\triangleright\psi(v\ng)\bigr)
\ot_Q\bigl(S^{-1}(h\1)\triangleright v\0\bigr)\cr
&=\sum\,q\,\psi(v\ng)\ot_Qv_0h=\la(qv)\,h,
}
$$
i.e., $\la$ is right $H$-linear. Hence $M\in\HdratCMH$ if and only if 
$V\in\HCM$.
\endproof

\proclaim
Theorem 8.10.
The functors in the statement of Theorem 8.2 induce equivalences
$$
\HMC\approx\HdxratHMA\,,\qquad\qquad\HCM\approx\HdxratAMH
$$
where $\HdxratHMA$ and $\HdxratAMH$ are the full subcategories, respectively, 
of the categories $\HdHMA$ and $\HdAMH$ consisting of those objects that are 
rationally extendible as objects of the categories $\HdMA$ and $\HdAM$.
\endproclaim

\Proof.
An object $W\in\HdHMA$ is rationally extendible if and only if $\,W\ot_AQ\,$ 
lies 
in $\HdratHMC$. Therefore the equivalence $\HdHMC\approx\HdHMA$ of Lemma 8.8 
induces an equivalence
$$
\HdratHMC\approx\HdxratHMA.
$$
Composed with the equivalence $\HMC\approx\HdratHMC$ of Lemma 8.9 it gives the 
first equivalence in the statement of Theorem 8.10. The second one is obtained 
similarly as the composite 
$\displaystyle\ \HCM\approx\HdratCMH\approx\HdxratAMH$.
\endproof

\section
9. Exactness of induction and flatness over coideal subalgebras

It was proved by Cline, Parshall, Scott \cite{Cl-PS77} and Oberst \cite{Ob77} 
that for a closed subgroup $K$ of an affine algebraic group $G$ the quotient 
$G/K$ is affine if and only if the induction functor from the category of 
rational $K$-modules to the category of rational $G$-modules is exact. Let 
$k[G]$ and $k[K]$ be the commutative Hopf algebras representing $G$ and $K$. 
The induction functor is given by the cotensor product 
$?{}\mathbin{\square}_{k[K]}k[G]$, and its exactness means that $k[G]$ is 
coflat over $k[K]$. The algebra of regular functions on the quotient $G/K$ may 
be identified with the left coideal subalgebra $B$ of $k[G]$ consisting of all 
regular functions on $G$ invariant with respect to the action of $K$ on $k[G]$ 
induced by right translations of $G$. Affineness of $G/K$ means that $k[G]$ is 
a faithfully flat $B$-module and 
$k[K]\cong k[G]\,/{\mskip1mu}k[G]{\mskip1mu}B^+$ where 
$B^+$ is the maximal ideal of $B$ consisting of all functions lying in $B$ 
which vanish at the identity element of $G$. By the Takeuchi correspondence 
described in \cite{Tak79} the affine case occurs precisely when $k[G]$ is 
faithfully coflat over $k[K]$, i.e., the induction functor is faithfully 
exact.

Takeuchi's results suggest that for an arbitrary Hopf algebra $H$, say with 
bijective antipode, left $H$-module factor coalgebras $C$ over which $H$ is 
faithfully coflat should be considered as analogs of closed subgroups of 
algebraic groups producing affine homogeneous spaces. This raises the question 
whether coflatness of $H$ over $C$ implies faithful coflatness. This is not 
true in general \cite{Sk25b, Cor. 1.4}, but the affirmative answer has long 
been known in the case when $C$ is a factor Hopf algebra of $H$ \cite{Doi83, 
Remark on p. 247}.

Suppose that $H$ satisfies our basic assumptions (A1) and (A2). Then it was 
shown in \cite{Sk10, Th. 0.3} that coflatness does imply faithful coflatness 
for a certain class of left $H$-module factor coalgebras $C$ of $H$. Employing 
the correspondence of Theorem 0.1 we are now able to remove any restriction on 
$C$ in that result.

\proclaim
Theorem 9.1.
Any left $H$-module factor coalgebra $C$ of $H$ is a simple object of the 
category $\HMC$ and a simple object of the category $\HCM$. Left or right 
coflatness of $H$ over $C$ implies faithful coflatness.
\endproclaim

\Proof.
Let $A$ be the left $\Hd$-invariant artinian subalgebra of the quotient ring $Q$ 
of $H$ corresponding to $C$ under the bijection of Theorem 0.1. By Lemma 8.4 we 
have $\Phi(HA)=C$ where $\Phi:\HdxratMA\to\MC$ is the equivalence of Theorem 8.2. 
Moreover, this evaluation remains true even if $\Phi$ is understood as the 
equivalence
$$
\HdxratHMA\to\HMC
$$
of Theorem 8.10. Therefore $C$ is a simple object of the category $\HMC$ if 
and only if $HA$ is a simple object of the category $\HdHMA$. But $HA$ is an 
$\HdHMA$-subobject of $Q$. Since $Q$ is a classical left quotient ring of $H$, 
its nonzero left $H$-submodules have nonzero intersections with $H$. 
Suppose that $W$ is an $\HdHMA$-subobject of $HA$. Then $W\cap H$ is a left 
$\Hd$-invariant left ideal of $H$, i.e., a left ideal which is also a right 
coideal. Therefore $W\cap H$ is either 0 or the whole $H$ since $H$ is a 
simple object of the category $\HMH$. Hence either $W=0$ or $W=HA$.

Thus $C$ is indeed a simple object of the category $\HMC$. Suppose that $H$ is 
left coflat over $C$, i.e., the cotensor product functor $?\sqCH:\MC\to\MH$ 
is exact. Since this functor is a right adjoint of the functor $\MH\to\MC$ 
induced by the canonical coalgebra homomorphism $H\to C$, the equality 
$V\sqCH=0$ for a right $C$-comodule $V$ means precisely that $\MC(W,V)=0$ for 
each right $H$-comodule $W$. Denote by $\calR$ the class of all right 
$C$-comodules $V$ with this property. By coflatness of $H$ over $C$ this class 
is closed under arbitrary colimits in $\MC$.

The coalgebra homomorphism $H\ot C\to C$ given by the action of $H$ on $C$ 
allows us to view $\MC$ as a left module category over the monoidal category 
$\MH$. Each finite-dimensional comodule $U\in\MH$ has a left dual $U^*$ in 
$\MH$, whence
$$
\MC(W,\,U\ot V)\cong\MC(U^*\ot W,\,V)
$$
for all $V,W\in\MC$. If $V\in\calR$ and $W\in\MH$, then the right hand side of 
the above formula vanishes since $U^*\ot W\in\MH$, and therefore 
$\MC(W,\,U\ot V)=0$ as well. This shows that $\calR$ is closed under the 
tensoring functors $U\ot{}?$ for finite-dimensional comodules in $\MH$, but 
then for all $U\in\MH$ by the local finiteness of comodules. In particular, 
$H\ot V\in\calR$ for each $V\in\calR$ where $H$ is a right $H$-comodule with 
respect to the comultiplication $\De$. 

If $V$ is a right coideal of $C$, then there is a morphism $H\ot V\to C$ given 
by the action of $H$ on $C$. Since its image $HV$ is an $\HMC$-subobject of 
the simple object $C$, we must have either $HV=0$ or $HV=H$. Furthermore, if 
$V\in\calR$, then $HV\in\calR$ too. In this case $HV\ne H$, and therefore 
$HV=0$, i.e., $V=0$.

Thus the class $\calR$ contains no nonzero right coideals of $C$. Since 
$\calR$ is closed under factor objects and $C$ is a cogenerator in $\MC$, it 
follows that $\calR=0$. This means that the exact functor $?\sqCH:\MC\to\MH$ 
is faithful.

It is deduced similarly that $C$ is a simple object of the category 
$\HCM$ because $AH$ is a simple object of the category $\HdAMH$, and right 
coflatness of $H$ over $C$ implies faithful coflatness.
\endproof

\proclaim
Corollary 9.2.
If $C$ is a left $H$-module factor coalgebra of $H,$ then $Hg=C$ for each 
grouplike $g\in C$.
\endproclaim

\Proof.
Indeed, $Hg$ is a nonzero $\HMC$-subobject of $C$.
\endproof

The properties of left $H$-module factor coalgebras established in Theorem 9.1 
translate into flatness of $H$ over right coideal subalgebras.

\proclaim
Theorem 9.3.
Let $A$ be a right coideal subalgebra of $H$. If the inclusion map 
$A\to\nobreak H$ admits a right (or left) $A$-linear retraction 
$H\to A,$ then $H$ is right (respectively, left) faithfully flat over $A$.
\endproclaim

\Proof.
Put $C=H/HA^+$ where $A^+=\{a\in A\mid\ep(a)=0\}$. By the hypothesis $A$ is an 
$\MA$-direct summand of $H$. Hence $H$ is a generator in $\MA$, and therefore 
$H$ is injective in $\MC$ \cite{Sk25a, Prop. 2.1}. This means that $H$ is a 
coflat right $C$-comodule, and so by Theorem 9.1 $H$ is right faithfully 
coflat over $C$. Furthermore, $A=\lco CH$ since $A$ is a dominion subalgebra 
of $H$ \cite{Ch24, Prop. 2.2}. Then $H$ is right faithfully flat over $A$ by 
Takeuchi's theorem \cite{Tak79, Th. 2}.
\endproof

\Remark.
For an arbitrary Hopf algebra $H$ over a field the conclusion of Theorem 9.3 
does not hold in general. However, if the map $A\to H$ admits an $A$-bimodule 
retraction, then $H$ is right faithfully flat over its right coideal subalgebra 
$A$. This was established by Chirvasitu \cite{Ch14, Prop. 1.4, Prop. 1.6}, 
applying in an essential way Mesablishvili's criterion for effective descent 
ring extensions \cite{Mes06, Th. 8.1}. A more elementary approach was proposed 
by Bichon \cite{Bich23}.

There is a short proof based on the fact that separable functors preserving 
epimorphisms reflect projective objects \cite{Nas-BO89, Prop. 1.2}. If $R$ is 
a subring of a ring $S$ such that $R$ is an $R$-bimodule direct summand of 
$S$, then the extension of scalars functor between module categories is 
separable \cite{Nas-BO89, Prop. 1.3}. In this case a right $R$-module $M$ is 
projective whenever so is the $S$-module $M\ot_RS$. If $M\in\MAH$, then 
$$
M\ot_AH\cong M/MA^+\ot H
$$
is a free $H$-module \cite{Tak79, p.~456}. Under the assumption that $A$ is an 
$A$-bimodule direct summand of $H$ it follows that all objects of the category 
$\MAH$ are projective in $\MA$, and then nonzero objects are projective 
generators in $\MA$.
\endremark

\section
10. Biideals are Hopf ideals

In a 1978 paper \cite{Nich78} Nichols pointed out several cases in which a 
biideal $I$ of a Hopf algebra $H$ is proved to be a Hopf ideal. Assuming that 
the base ring is a field, any weakly finite factor bialgebra $H/I$ is a Hopf 
algebra. This statement was given in \cite{Er-Sk09, Cor. 2.4}, although it 
can be easily proved along the line of reasoning found in \cite{Nich78, Th. 1} 
for the cases when $H$ is either finite-dimensional or commutative. It was 
also proved in \cite{Nich78} that each biideal of $H$ is a Hopf ideal when $H$ 
is either pointed or cocommutative.

In terms of the Takeuchi correspondence between coideal left ideals and right 
coideal subalgebras of $H$ the right coideal subalgebra corresponding to a 
biideal is stable under the right adjoint action of $H$ on itself \cite{Tak94, 
Prop. 1.3}. Conversely, if a right coideal subalgebra $R$ is stable under the 
right adjoint action of $H$, then the left ideal $HR^+$ generated by the 
augmentation ideal $R^+$ of $R$ is a Hopf ideal \cite{Tak94, Prop. 1.4}. As a 
consequence, a biideal $I$ is a Hopf ideal whenever $H$ is either right or 
left faithfully coflat over $H/I$ \cite{Tak94, Cor. 1.5}.

It will be proved in this section that all biideals are Hopf ideals for a Hopf 
algebra $H$ satisfying our basic assumptions (A1) and (A2). Such a Hopf 
algebra is weakly finite. However, weak finiteness does not necessarily pass 
over to factor rings, and for this reason we cannot apply the result from 
\cite{Er-Sk09}. The proof we provide is much more intricate. Essential work 
will be done on the level of equivariant modules over $\Hd$-invariant artinian 
subalgebras of the quotient ring $Q$.

\setitemsize(a)
\proclaim
Theorem 10.1.
For a coideal left ideal $I$ of $H$ and the corresponding subalgebra $A$ of 
$Q$ the following conditions are equivalent:

\item(a)
\ $I$ is a biideal of $H,$

\item(b)
\ $I$ is a Hopf ideal of $H,$

\item(c)
\ $A$ is stable under the right adjoint action of $H$ on $Q$ defined by the rule
$$
q\triangleleft h=\sum\,S(h\1)\,q\,h\2\,,\qquad q\in Q,\ h\in H.\eqno(10.1)
$$

\endproclaim

\Proof.
The most difficult part of this theorem is to show that (a) implies (b). 
Assuming that $I$ is a biideal, the quotient $C=H/I$ is a factor bialgebra of 
$H$. The category $\MC$ is then monoidal with respect to the tensor product 
of comodules. By Theorem 8.2 $\MC\approx\HdxratMA$. In Proposition 10.11 we 
will describe the corresponding monoidal structure on the category 
$\HdxratMA$. By Corollary 8.3 finite-dimensional $C$-comodules correspond to 
$A$-finite objects of $\HdxratMA$. In Proposition 10.12 it will be shown that 
each $A$-finite object has a left dual in $\HdxratMA$. Hence each 
finite-dimensional right $C$-comodule has a left dual in $\MC$. But this is a 
necessary and sufficient condition for a bialgebra to be a Hopf algebra (see 
Ulbrich \cite{Ulb90}), and so (a)$\,\Rar\,$(b).

Thus all arguments which prove the implication (a)$\,\Rar\,$(b) will be given 
in the rest of this section. For now we will prove only that (a)$\,\Lrar\,$(c).

Since $I$ is a coideal and a left ideal, it is a biideal of $H$ if and only if 
it is also a right ideal. By Proposition 4.10 there is a $k$-linear bijection 
$\,\psi:H\to(Q\ot Q)^{\Hd}\!$ such that $\,\psi(x)=\sum x\1\ot S(x\2)\,$ for 
$x\in H$. Define a right $H$-module structure on $\,Q\ot Q\,$ by the rule
$$
(a\ot b)\prec h=\sum\,ah\1\ot S(h\2)b,\qquad a,b\in Q,\ h\in H.\eqno(10.2)
$$
Since $\psi(xy)=\psi(x)\prec y$ for all $x,y\in H$, the subspace of 
$\Hd$-invariants $(Q\ot Q)^{\Hd}$ is \emph{$\prec$-stable}, i.e. stable under 
the action of $H$ defined in (10.2). We see also that $I$ is a right ideal of 
$H$ if and only if its image $\psi(I)$ is $\prec$-stable.

The coideal $\scrI$ of the $Q$-coring $Q\ot Q$ corresponding to $I$ is the 
$Q$-subbimodule generated by $\psi(I)$, and also 
$\psi(I)=\scrI^{\Hd}\!\!=\scrI\cap(Q\ot Q)^{\Hd}\!$. Since the action (10.2) 
commutes with both the left and right actions of $Q$ defining the bimodule 
structure on $Q\ot Q$, it follows that $\psi(I)$ is $\prec$-stable if and 
only if so is $\scrI$.

On the other hand, $\scrI$ corresponds to the algebra $A$, and so $\scrI$ 
is the subbimodule of $\,Q\ot Q\,$ generated by $\de(A)$ where the map 
$\,\de:Q\to Q\ot Q\,$ is defined by the rule $\,\de(x)=1\ot x-x\ot 1\,$ for 
$x\in Q$. Therefore $\scrI$ is $\prec$-stable if and only if $\scrI$ contains
$$
\eqalign{
\de(a)\prec h
&=\sum\,\bigl(h\1\ot S(h\2)a-ah\1\ot S(h\2)\bigr)\cr
&=\sum\,h\1\,\bigl(1\ot S(h\2)ah\3-S(h\2)ah\3\ot1\bigr)\,S(h\4)\cr
&=\sum\,h\1\,\de(a\triangleleft h\2)\,S(h\3)
}
$$
for all $a\in A$ and $h\in H$. This containment does hold when 
$A\triangleleft H\sbs A$. Hence (c) implies (a). Conversely, suppose that 
$\,\de(a)\prec h\in\scrI\,$ for all $a\in A$ and $h\in H$. Then
$$
\de(a\triangleleft h)
=\sum\,S(h\1)\,\bigl(\de(a)\prec h\2)\,h\3\in\scrI,
$$
whence $\,a\triangleleft h\in A\,$ since $\,A=\{x\in Q\mid\de(x)\in\scrI\}\,$ 
by the correspondence of Theorem 1.1. This shows that (a)$\,\Rar\,$(c).
\endproof

\proclaim
Corollary 10.2.
If $I$ is a biideal of $H,$ then $H/I$ is a Hopf algebra with bijective 
antipode.
\endproclaim

\Proof.
By Theorem 10.1 $S(I)\sbs I$. Applying this conclusion to the Hopf algebra 
$H\op$ in which $I$ is still a biideal, we also get $S^{-1}(I)\sbs I$. 
Since the antipode of $H/I$ is induced by the bijective antipode $S$ of $H$, it 
is bijective too.
\endproof

Exactly as in the already mentioned paper of Nichols \cite{Nich78}, we can 
weaken slightly the assumption about the ideal $I$ by omitting the condition 
$\ep(I)=0$ required for biideals.

\proclaim
Corollary 10.3.
Suppose that $I$ is an ideal of $H$ such that $\,\De(I)\sbs I\ot H+H\ot I$ 
and $I\ne H$. Then $I$ is a Hopf ideal.
\endproclaim

\Proof.
The ideal $I^+=I\cap\Ker\ep$ is also a coideal, and so a biideal of $H$. By 
Theorem 10.1 $H/I^+$ is a Hopf algebra. Suppose that $I\ne I^+$. Then $I/I^+$ 
is a one-dimensional ideal of $H/I^+$ spanned by an element $x$ such that 
$\ep(x)=1$ and $\De(x)=a\ot x+x\ot b$ for some $a,b\in H/I$. Since 
$\,(\ep\ot\id)\circ\De=(\id\ot\ep)\circ\De=\id$, both $a$ and $b$ are scalar 
multiples of $x$, and it follows that $x$ is a grouplike. But then $x$ is 
invertible in $H/I$. Hence $I/I^+=H/I^+$, which contradicts the assumption $I\ne H$.
\endproof

Let $A$ be a left $\Hd$-invariant artinian subalgebra of $Q$. Consider the 
corresponding $Q$-coring $\scrC=Q\ot_AQ$. The next lemma provides a 
reformulation of condition (c).

\proclaim
Lemma 10.4.
The algebra $A$ is stable under the right adjoint action of $H$ if and only if 
$\,ca=ac\,$ for all $a\in A$ and $c\in\scrC^{\Hd}\!$.
\endproclaim

\Proof.
As an intermediate step in the proof of the equivalence (a)$\,\Lrar\,$(c) of 
Theorem 10.1 we have shown that (c) holds if and only if
$$
\sum\,\bigl(h\1\ot S(h\2)a-ah\1\ot S(h\2)\bigr)\in\scrI\quad
\text{for all $a\in A$ and $h\in H$}
$$
where the coideal $\scrI$ is the kernel of the canonical homomorphism 
$Q\ot Q\to\scrC$. This condition is equivalent to the identity
$$
h\1\,g\,S(h\2)\,a=a\,h\1\,g\,S(h\2),\qquad a\in A,\ h\in H,\eqno(10.3)
$$
in the coring $\scrC$ where $g=1\ot_A1$ is the distinguished grouplike of 
$\scrC$. By Corollary 4.11 $\,\scrC^{\Hd}\!$ consists precisely of those 
elements $c\in\scrC$ which can be written as
$$
c=h\1\,g\,S(h\2)
$$
for some $h\in H$. Hence (10.3) amounts to the condition that $\,ca=ac\,$ 
for all $a\in A$ and $c\in\scrC^{\Hd}\!$.
\endproof

\proclaim
Lemma 10.5.
Each object $M\in\HdratMQ$ has a uniquely determined left $Q$-module structure 
such that
$$
qv=vq\quad\text{for all $\,q\in Q\,$ and $\,v\in M^{\Hd}\!$}.\eqno(10.4)
$$
This additional structure makes $M$ into an object of the category $\,\HdratQMQ$. 

In this way we obtain a functor $\,\HdratMQ\to\HdratQMQ\,$ and an isomorphism 
of categories $\,\HdratMQ\to\HdratQM\,$ which leave all morphisms unchanged.

Moreover, this functor gives an isomorphism of $\HdratMQ$ onto a monoidal 
subcategory of the category $\HdQMQ,$ and the functor $\,\HdratMQ\to\Mk\,,$ 
$\,M\mapsto M^{\Hd}\!,$ is monoidal.
\endproclaim

\Proof.
By Proposition 3.2 $M\cong V\ot Q$ where $V=M^{\Hd}\!$. The required left 
$Q$-module structure on $M$ is obtained by means of left multiplications on 
the ring $Q$. The resulting left $Q$-linear map $Q\ot V\to M$ is bijective and 
is also $\Hd$-linear since the action of $\Hd$ on $V$ is trivial. Hence $M$ is 
isomorphic to $Q\ot V$ as an object of the category $\HdratQM$. Since the left 
action of $Q$ commutes with the right one, $M$ becomes an object of the 
category $\HdratQMQ$.

Each morphism $\,\ph:M\to N\,$ in $\,\HdratMQ\,$ is the right $Q$-linear extension of 
some $k$-linear map $\,M^{\Hd}\!\to N^{\Hd}\!$. It is clear from (10.4) that 
$\ph$ is also left $Q$-linear, and therefore a morphism in $\HdratQMQ$. Thus, 
considering each object $M$ of the category $\HdratMQ$ with the additional 
left $Q$-module structure defined above we obtain a functor 
$\,\HdratMQ\to\HdratQMQ\,$ which is an isomorphism onto the full subcategory 
of $\,\HdratQMQ\,$ consisting of objects satisfying (10.4).

In view of Proposition 3.3 the functor $\,\HdratMQ\to\HdratQM\,$ obtained by 
forgetting the right $Q$-module structure has the inverse functor 
$\,\HdratQM\to\HdratMQ\,$ constructed similarly.

The category $\HdQMQ$ is monoidal with respect to the tensor product $\ot_Q\,$, 
and we may identify $\HdratMQ$ with its full subcategory. The unit object $Q$ 
lies in $\HdratMQ\,$, and $\,Q^{\Hd}\!\!=k\,$ is the unit object 
of $\Mk$. If $\,M,\,N$ are two objects of the category $\,\HdratMQ\,$, then
$$
M\ot_QN\cong(M^{\Hd}\!\ot Q)\ot_Q(N^{\Hd}\!\ot Q)
\cong(M^{\Hd}\!\ot N^{\Hd})\ot Q.
$$
It is clear that the canonical $k$-linear map 
$\,M^{\Hd}\!\ot N^{\Hd}\!\to(M\ot_QN)^{\Hd}\!$ is bijective, and 
a straightforward check shows that the rationally generated object $M\ot_QN$ 
of the category $\HdQMQ\,$ satisfies (10.4), i.e., $\,M\ot_QN\in\HdratMQ\,$.
\endproof

Further on in this section we assume that the left $\Hd$-invariant artinian 
subalgebra $A$ of $Q$ is stable under the right adjoint action of $H$. 
By the already proved part of Theorem 10.1 this implies that the corresponding 
coideal $I$ of $H$ is a biideal, and therefore $H/I$ is a factor bialgebra 
of $H$.

\proclaim
Lemma 10.6.
For $M\in\HdratMC$ the comodule structure map $\,\rho:M\to M\ot_Q\scrC\,$
is left $A$-linear with respect to the left action of $A$ on $M$ obtained from 
the left action of $Q$ that was defined in Lemma 10.5.
\endproclaim

\Proof.
By Proposition 4.12 $\,V=M^{\Hd}\!$ has a right $C$-comodule structure such 
that $\,M\cong V\ot Q\,$ and
$$
\rho(v)=\sum v\0\ot_Q\psi(v\1)\quad\text{for all $v\in V$}.
$$
where $\,\psi:C\to\scrC^{\Hd}\!$ is the canonical isomorphism. 
Making use of Lemma 10.4 and formula (10.4), we get
$$
\eqalign{
\rho(avq)=\rho(vaq)&=\sum v\0\ot_Q\psi(v\1)\,aq\cr
&=\sum v\0\ot_Qa\,\psi(v\1)\,q\cr
&=\sum v\0a\ot_Q\psi(v\1)\,q\cr
&=\sum av\0\ot_Q\psi(v\1)\,q=a\,\rho(vq)
}
$$
for all $a\in A$, $v\in V$, and $q\in Q$. So $\rho$ is indeed left $A$-linear.
\endproof

\proclaim
Corollary 10.7.
For each object $M\in\HdratMC$ its $\HdMA$-subobject $M\co\scrC$ is also a 
left $A$-submodule. By means of the left $A$-module structures obtained in 
this way the assignment $\,M\mapsto M\co\scrC\,$ defines a functor 
$\,\HdratMC\to\HdAMA$.
\endproclaim

\Proof.
By Lemma 10.6 each object $M\in\HdratMC$ with its left $A$-module structure is 
an object of the category $\HdAMC$ introduced in section 8. Each morphism 
in $\HdratMC$, being a morphism in $\HdratMQ$, is left $Q$-linear, and so, in 
particular, left $A$-linear by Lemma 10.5. This means that all morphisms in 
$\HdratMC$ are morphisms in $\HdAMC$. 

We thus obtain a functor $\HdratMC\to\HdAMC$ which adds the left $A$-module 
structure defined by formula (10.4) to each object of $\HdratMC$. As a special 
case of Lemma 8.8, the assignment $\,M\mapsto M\co\scrC\,$ gives a functor 
$\,\HdAMC\to\HdAMA$ which is even an equivalence. The latter is just 
restricted to $\HdratMC$.
\endproof

\proclaim
Corollary 10.8.
Let $\,F:\HdAMA\to\HdMA\,$ be the forgetful functor. There is a functor 
$$
E:\HdxratMA\to\HdAMA
$$
such that $F\circ E$ is the inclusion functor $\,\HdxratMA\to\HdMA\,$ and $E$ is 
naturally isomorphic to the composite
$$
\HdxratMA\to\HdratMC\to\HdAMA\eqno(10.5)
$$
of the functors described in Theorem 8.2 and Corollary 10.7.

This functor gives an isomorphism of the category $\HdxratMA$ onto a full 
subcategory $\calE$ of $\HdAMA$ which has the property that for each object of 
$\HdAMA$ lying in $\calE$ all its subobjects and factor objects also lie in 
$\calE$.
\endproclaim

\Proof.
The first functor in (10.5) is induced by the functor $\,\HdMA\to\HdMC\,$ 
shown in (5.7). It takes an object $\,W\in\HdxratMA\,$ to 
$\,W\ot_AQ\in\HdratMC$. Hence the composite (10.5) takes $W$ further to 
$\,(W\ot_AQ)\co\scrC$. By Lemma 5.4 there are natural isomorphisms
$$
W\cong(W\ot_AQ)\co\scrC\quad\text{in $\HdMA$}.
$$
The functor $E$ is obtained by equipping each object of $\HdxratMA$ with the 
additional left $A$-module structure arising from those isomorphisms. Thus 
$E(W)$, regarded as an object of $\HdMA$, coincides with $W$, i.e., 
$F\bigl(E(W)\bigr)=W$, and the morphisms do not change either.

Denote by $\calE$ the full subcategory of $\HdAMA$ consisting of those objects 
$X$ which are rationally extendible as objects of the category $\HdMA$ and 
fulfill the equality $X=E\bigl(F(X)\bigr)$. Then $E$ is a functor 
$\,\HdxratMA\to\calE$, and the forgetful functor $\,\calE\to\HdxratMA\,$ is 
its inverse. Hence $E$ is an isomorphism of $\HdxratMA$ onto $\calE$.

Suppose that $X\in\calE$ and $Y$ is an $\HdAMA$-subobject of $Y$. We may 
identify $Y$ with its set-theoretic image in $X$. Since $X$ is rationally 
extendible as an object of $\HdMA$, so too is $Y$ by Proposition 8.1. The 
equality $F\bigl(E(F(Y))\bigr)=F(Y)$ means that $E\bigl(F(Y)\bigr)$ and $Y$ 
coincide as objects of $\HdMA$. But these two objects are both left 
$A$-submodules of $X$ sharing the same set of elements, and therefore they 
have the same left action of $A$ too. Hence $E\bigl(F(Y)\bigr)=Y$, which shows 
that $Y\in\calE$. Then $X/Y\in\calE$ as well since the functors are exact.
\endproof

Thus each object $W\in\HdxratMA$ has a canonical left $A$-module structure. 
This structure is not defined straightforwardly in terms of $W$ itself, as one 
has to pass to the extension $W\ot_AQ$ of $W$, and then use formula (10.4). 
Nevertheless the left action of $A$ on rationally generated objects is 
described easily.

\proclaim
Lemma 10.9.
If $U$ is a rational $\Hd$-submodule of an object $W\in\HdxratMA,$ then
$$
au=\sum\,u\0(a\triangleleft u\1)\quad
\text{for all $a\in A$ and $u\in U$}.\eqno(10.6)
$$
As a consequence, $\,A\,U=UA$. If\/ $W$ is freely generated by\/ $U\!$ as a 
right $A$-module, then\/ $W$ is also freely generated by\/ $U\!$ as a left 
$A$-module.
\endproclaim

\Proof.
Put $M=W\ot_AQ$. There is a map $\,\ph:U\to M^{\Hd}\!$ defined by the rule 
$$
\ph(u)=\sum u\0\ot_AS(u\1),\qquad u\in U.
$$
Then $\,u\ot_A1=\sum\ph(u\0)u\1$. Taking $a\in A$ and applying (10.4), 
we get
$$
\eqalign{
au\ot_A1=\sum\ph(u\0)au\1&=\sum\ph(u\0)u\1S(u\2)au\3\cr
&=\sum u\0\ot_A(a\triangleleft u\1)=\sum u\0(a\triangleleft u\1)\ot_A1
}
$$
in $M$, which implies equality (10.6) in $W$. It shows that $\,A\,U\sbs UA$. 
The opposite inclusion is also true since (10.6) is equivalent to the identity
$$
ua=\sum\,\,(a\triangleleft S^{-1}u\1)\,u\0\quad
\text{for all $a\in A$ and $u\in U$}.\eqno(10.7)
$$
Moreover, the $k$-linear map $\,\tau:A\ot U\to U\ot A\,$ defined by the 
assignment
$$
a\ot u\mapsto\sum\,u\0\ot(a\triangleleft u\1)
$$
is bijective. The map $\,\la:A\ot U\to W\,$ afforded by the left module 
structure is the composite of $\tau$ with the map $\,\mu:U\ot A\to W\,$ 
afforded by the right module structure. Hence $\la$ is bijective if and only if 
so is $\mu$.  
\endproof

\proclaim
Lemma 10.10.
Each object $\,W\!\in\HdxratMA\,$ is free both as a right $A$-module and as a 
left $A$-module. The right rank\/ $\rrk_AW\!$ and the left rank\/ $\lrk_AW\!$ 
of\/ $W\!$ over $A$ are equal.
\endproclaim

\Proof.
We know from Corollary 8.3 that $W$ is a free right $A$-module. Suppose first 
that its rank $\rrk_AW\!$ is finite. By Corollary 8.5 $W$ embeds in a 
rationally generated object $W'$ of the category $\HdMA$. Since the rational 
$\Hd$-module $\Rat(W')$ is the union of the directed family $\calF$ of its 
finite-dimensional submodules, there exists $U\in\calF$ such that $W\sbs UA$. 
By Lemma 10.9 $\,UA=\!A\,U$. Hence $UA$ is a finitely generated left 
$A$-module, and therefore so is $W$. Since $W$ is a left $A$-finite object of 
the category $\HdAMA$, it is left $A$-free by Lemma 2.11 applied to the left 
$(\Hd)\cop$-module algebra $A\op$.

Put $M=W\ot_AQ\in\HdratMC$. This object is a free right $Q$-module of rank 
equal to $\rrk_AW$. By Lemma 10.5 it is also a free left $Q$-module of the same 
rank. The left $A$-module structure on $W\cong M\co\scrC$ is the restriction 
of the left $Q$-module structure on $M$. The map $\,Q\ot_AW\to M\,$ arising from 
the left action of $Q$ is a morphism in $\HdQM$. Hence its image $N$ is an 
$\HdQM$-subobject of $M$. However, by Lemma 10.5 the $\HdQM$-subobjects of 
$M$ coincide with the $\HdMQ$-subobjects. Since $W$ generates $M$ as a right 
$Q$-module, we get $N=M$. This implies that the rank of the free left 
$Q$-module $\,Q\ot_AW$ equal to $\,\lrk_AW$ cannot be less than that of $M$, 
i.e., 
$$
\lrk_AW\ge\rrk_AW.
$$
Suppose that $0\to W'\to W\to W''\to0$ is an exact sequence in $\HdxratMA$. 
Then
$$
\lrk_AW=\lrk_AW'+\lrk_AW''\qquad\text{and}\qquad\rrk_AW=\rrk_AW'+\rrk_AW''.
$$
As we have proved, $\,\lrk_AW'\ge\rrk_AW'$ and $\,\lrk_AW''\ge\rrk_AW''$. It 
follows that these inequalities become equalities whenever $\,\lrk_AW=\rrk_AW$. 
In other words, if the equality of the left and right ranks holds for $W$, 
then it holds for all subobjects and factor objects of $W$.

By Lemma 10.9 both $\,\lrk_AW$ and $\,\rrk_AW$ are equal to $\,\dim_kU$ when 
$W\cong U\ot A$ is freely generated as a right $A$-module by a 
finite-dimensional rational $\Hd$-module. Each $A$-finite rationally generated 
object of $\HdxratMA\,$ is isomorphic to a factor object of such an object 
$U\ot A$, and an arbitrary $A$-finite object is isomorphic to a subobject of a 
rationally generated one by Corollary 8.5. It follows that
$$
\lrk_AW=\rrk_AW
$$
for all $A$-finite objects of the category $\,\HdxratMA$.

Even if $W$ is not $A$-finite, it is still locally $A$-finite by Corollary 8.3. 
Consider the set of all quadruples $(X,B_l,B_r,\si)$ where $X$ is a subobject 
of $W$, $\,B_l$ and $B_r$ are bases for $X$ as a left and a right $A$-module, 
respectively, and $\si$ a bijection of $B_l$ onto $B_r$. Define an order 
relation on this set by setting
$$
(X,B_l,B_r,\si)\le(X',B'_l,B'_r,\si')\quad\text{if}
\ \,X\sbs X',\ \,B_l\sbs B_l',\ \,B_r\sbs B_r',\ \,\si=\si'|_{B_l}
$$
If $X\ne W$, then there exists a subobject $X'$ of $W$ such that $X$ is 
properly contained in $X'$ and $X'/X$ is an $A$-finite object of the category 
$\,\HdxratMA$. As we have proved already, then $X'/X$ is left and right 
$A$-free of equal ranks. This implies that for any quadruple with the first 
component $X$ there exists a larger quadruple with the first component $X'$. 
By Zorn's lemma this set of quadruples has a maximal element, and it follows 
that the first component of any maximal quadruple has to be the whole $W$. 
Hence $\,\lrk_AW=\rrk_AW$ as an equality of cardinals.
\endproof

Recall that the category $\HdAMA$ is monoidal with respect to the tensor product 
functor $\ot_A$. The functor described in Corollary 10.8 allows us to identify 
$\HdxratMA$ with a full subcategory of $\HdAMA$, which we do further on.

\proclaim
Proposition 10.11.
The category $\HdxratMA$ is a monoidal subcategory of the category $\HdAMA$.
The functor\/ $\Phi:\HdxratMA\to\MC$ of Theorem 8.2 is a monoidal equivalence.
\endproclaim

\Proof.
The object $A\in\HdxratMA$ is generated as a right $A$-module by its vector 
subspace $U=k$ on which $\Hd$ acts trivially. Formula (10.6) shows that the 
left action of $A$ on itself provided by the functor $\HdxratMA\to\HdAMA$ 
coincides with the action by left multiplications. With respect to this 
structure $A$ is the unit object of the monoidal category $\HdAMA$.
Furthermore, $\Phi$ takes it to the trivial $C$-comodule 
$\,\Phi(A)\cong k1_C\,$ by Lemma 8.4.

We also have to show that $W\ot_AW'\in\HdxratMA$ for all $W,W'\in\HdxratMA$. 
If $X$ is a subobject of $W$ and $X'$ is a subobject of $W'$, then the 
canonical map $X\ot_AX'\to W\ot_AW'$ is a monomorphism in the category 
$\HdAMA$. Injectivity of this map is a consequence of the fact that all 
objects of $\HdxratMA$ are left and right $A$-free. By Corollary 10.8 
$X\ot_AX'$ lies in $\HdxratMA$ whenever so does $W\ot_AW'$. Therefore 
Corollary 8.5 reduces the verification to the case when $W$ and $W'$ are 
rationally generated.

So let $W=UA$ and $W'=U'A$ for some rational $\Hd$-submodules $U\sbs W$ and 
$U'\sbs W'$. Since $AU'=U'A$ by Lemma 10.9, we have $W\ot_AW'=TA$ where $T$ is 
the image of the canonical $\Hd$-linear map $\,U\ot\,U'\to W\ot_AW'$. Since 
$T$ is a rational $\Hd$-module, we see that $W\ot_AW'$ is rationally generated, 
and therefore rationally extendible, as an object of $\HdMA$. Furthermore,
$$
\eqalign{
au\ot_A\!v
=\sum\,u\0(a\triangleleft u\1)\ot_A\!v
&=\sum\,u\0\ot_A(a\triangleleft u\1)\,v\cr
&=\sum\,u\0\ot_Av\0\bigl(a\triangleleft(u\1v\1)\bigr)
}
$$
for all $u\in U$, $v\in U'$, and $a\in A$. This shows that formula (10.6) is 
satisfied for the left action of $A$ on elements of $T$, and so the left 
$A$-module structure of the object $\,W\ot_AW'\in\HdAMA\,$ does indeed 
coincide with the one provided by the functor $\,\HdxratMA\to\HdAMA$.

Recall that $\,\Phi(W)=(W\ot_AQ)^{\Hd}$ with the $C$-comodule structure induced 
by the $\scrC$-comodule structure on $W\ot_AQ$. Using the left $Q$-module 
structures defined by formula (10.4) on objects of the category $\HdratMQ$, 
we get a canonical bijection
$$
(W\ot_AQ)\ot_Q(W'\ot_AQ)\to W\ot_A(W'\ot_AQ)\cong(W\ot_AW')\ot_AQ\eqno(10.8)
$$
which is an isomorphism in $\HdMQ$. By Lemma 10.5 the category $\HdratMQ$ is 
monoidal with respect to the tensor product $\ot_Q$. Hence (10.8) is in fact 
an isomorphism in $\HdratMQ$. Applying the monoidal functor 
$$
\HdratMQ\to\Mk\,,\qquad M\mapsto M^{\Hd}\!,
$$
we obtain a $k$-linear bijection
$$
\al_{W,W'}:\ \Phi(W)\ot\Phi(W')\to\Phi(W\ot_AW')
$$
for each pair of objects $W,W'\in\HdxratMA$. These bijections are natural in 
both arguments and compatible with the associativity and unit constraints. We 
still have to check that they are in fact morphisms in the category $\MC$.

If $\al_{W,W'}$ is $C$-colinear for some pair $W$ and $W'$, then it follows by 
naturality of those bijections that $\al_{X,X'}$ is $C$-colinear when either 
$X$ is a subobject of $W$ and $X'$ is a subobject of $W'$ or
$X$ is a factor object of $W$ and $X'$ is a factor object of $W'$.

Therefore it suffices to consider the case when $W=U\ot A$ and $W'=U'\ot A$ 
for some rational $\Hd$-modules $\,U,\,U'$. In this case 
$\,W\ot_AW'\cong(U\ot U')\ot A$, and
$$
W\ot_AQ\cong U\ot Q,\quad
W'\ot_AQ\cong U'\ot Q,\quad
(W\ot_AW')\ot_AQ\cong(U\ot U')\ot Q.
$$
There are $k$-linear bijections $\ph:U\to\Phi(W)$ and $\ph':U'\to\Phi(W')$ 
such that
$$
\ph(u)=\sum\,u\0\ot S(u\1),\qquad
\ph'(v)=\sum\,v\0\ot S(v\1)
$$
for $u\in U$ and $v\in U'$. Another use of formula (10.4) yields
$$
\eqalign{
\ph(u)\ot_Q\ph'(v)&=\sum\,\bigl(u\0\ot S(u\1)\bigr)\ot_Q\ph'(v)\cr
&=\sum\,(u\0\ot1)\ot_Q\ph'(v)S(u\1)\cr
&=\sum\,(u\0\ot1)\ot_Q\bigl(v\0\ot S(v\1)S(u\1)\bigr)
}
$$
in $\,(W\ot_AQ)\ot_Q(W'\ot_AQ)$, which implies that
$$
\al_{W,W'}\bigl(\ph(u)\ot\ph'(v)\bigr)
=\sum\,(u\0\ot v\0)\ot S(v\1)S(u\1).
$$
Each rational $\Hd$-module is a right $H$-comodule, and we may view it as a 
right $C$-comodule via the canonical homomorphism of coalgebras $H\to C$. With 
this convention $\ph$ and $\ph'$ are isomorphisms in $\MC$ by Lemma 8.6, 
and we see from the last formula that $\,\al_{W,W'}\circ(\ph\ot\ph')\,$ is a 
similar isomorphism
$$
U\ot U'\to\Phi(W\ot_AW')
$$
in $\MC$. Hence $\al_{W,W'}$ is an isomorphism in $\MC$ as well.
\endproof

For an object $W\in\HdxratMA$ we consider its dual as a right $A$-module
$$
W_A^*=\MA(W,A).
$$
Lemma 2.8 provides $W_A^*$ with a left $\Hd$-module structure, and $W_A^*$ is 
an $A$-bimodule with respect to the two actions of $A$ defined by the rule
$$ 
(afb)(x)=af(bx),\qquad a,b\in A,\ f\in W_A^*,\ x\in W.
$$
These module structures are compatible so that $W_A^*$ is an object of $\,\HdAMA$.

\proclaim
Proposition 10.12.
If\/ $W$ is an $A$-finite object of the category $\HdxratMA,$ then so is 
$W_A^*$. Moreover, $W_A^*$ is the left dual of\/ $W$ in the monoidal category 
$\HdAMA$.
\endproclaim

\Proof.
The assignment $W\mapsto W_A^*$ defines a contravariant functor 
$$
\HdxratMA\to\HdAMA
$$
which is exact by the $A$-freeness of all objects of the category $\HdxratMA$. 
Under this functor subobjects of $W$ go to factor objects of $W_A^*$, and factor 
objects of $W$ go to subobjects of $W_A^*$. If the containment 
$W_A^*\in\HdxratMA$ holds for some object $W$, then it holds also for all 
subobjects and factor objects of $W$ by the property of the subcategory 
$\HdxratMA\sbs\HdAMA$ stated in Corollary 10.8.

Therefore it suffices to prove that $W_A^*$ is an $A$-finite object of 
$\HdxratMA$ in the case when $W=U\ot A$ for some finite-dimensional 
right $H$-comodule. The dual vector space $U^*$ has a right $H$-comodule 
structure such that
$$
\sum\,\xi\0(u)\ot\xi\1=\sum\,\xi(u\0)\ot S(u\1)\quad\text{for all $\xi\in U^*$ 
and $u\in U$}.\eqno(10.9)
$$
Equipped with this structure $U^*$ is the left dual of $U$ in the monoidal 
category $\MH$. There is an obvious $\Hd$-linear map $U^*\to W_A^*$ which 
extends to an isomorphism
$$
A\ot U^*\to W_A^*
$$
in the category $\HdAM$. For convenience we identify $U$ with the subspace 
$U\ot k$ of $W$ and $U^*$ with its image in $W_A^*$. Applying (10.6) and 
(10.9), we get
$$
\xi(au)
=\sum\,\xi(u\0)(a\triangleleft u\1)
=\sum\,(a\triangleleft u\1)\,\xi(u\0)
=\sum\,\bigl(a\triangleleft S^{-1}(\xi\1)\bigr)\,\xi\0(u).
$$
for all $\xi\in W_A^*$, $a\in A$, and $u\in U$. This means that
$$
\xi a=\sum\,\bigl(a\triangleleft S^{-1}(\xi\1)\bigr)\,\xi\0,
$$
i.e., the rational $\Hd$-submodule $U^*$ of $W_A^*$ satisfies identity (10.7). 
But (10.7) is equivalent to (10.6). Hence $U^*\!A=A\,U^*=W_A^*$. It follows 
that $W_A^*$ is rationally generated and $A$-finite as an object of the 
category $\HdMA$. Furthermore, the left $A$-module structure on $W_A^*$ 
coincides with the one obtained by means of the functor $\HdxratMA\to\HdAMA$ 
from the right module structure.

Thus for each $A$-finite object $W\in\HdxratMA$ we have proved that $W_A^*$ 
lies in the subcategory $\HdxratMA$ of $\HdAMA$, and $W_A^*$ is also 
$A$-finite. Since $W$ is a free right $A$-module, the canonical map
$$
W\ot_AW_A^*\to\MA(W,W)\eqno(10.10)
$$
is bijective. Here $\MA(W,W)$ is the endomorphism ring of $W$ as a right 
$A$-module. The left action of $A$ on $W$ provides a ring homomorphism 
$A\to\MA(W,W)$ which allows us to view $\MA(W,W)$ as an $A$-bimodule. The 
action of $\Hd$ defined in Lemma 2.8 makes $\MA(W,W)$ into a left $\Hd$-module, 
and with respect to these module structures the map in (10.10) is an isomorphism 
in the category $\HdAMA$.

The identity transformation $\,\Id_W:W\to W$ is an $\Hd$-invariant $A$-central 
element of $\MA(W,W)$. The corresponding element $z\in W\ot_AW_A^*$ is then 
also $\Hd$-invariant and \emph{$A$-central} in the sense that $az=za$ for 
all $a\in A$. It follows that there is a morphism
$$
\coev:A\to W\ot_AW_A^*
$$
in $\HdAMA$ such that $1\mapsto z$, and so $1\mapsto\Id_W$ under the composite 
of $\coev$ with (10.10). The evaluation map
$$
\ev:W_A^*\ot_AW\to A,\qquad f\ot_Ax\mapsto f(x),\eqno(10.11)
$$
is also a morphism in $\HdAMA$. Both composites
$$
\displaylines{
W\cong A\ot_AW\lmapr6{\coev\ot\id}W\ot_AW_A^*\ot_AW
\lmapr6{\id\ot\ev}W\ot_AA\cong W,\cr
\noalign{\smallskip}
W_A^*\cong W_A^*\ot_AA\lmapr6{\id\ot\coev}W_A^*\ot_AW\ot_AW_A^*
\lmapr6{\ev\ot\id}A\ot_AW_A^*\cong W_A^*
}
$$
are the identity maps. This is straightforward to see, and one may also 
observe that map (10.11) factors through the canonical bijection
$\,W_A^*\ot_{\MA(W,W)}W\to A\,$
which forms, together with (10.10), the Morita context connecting $A$ and 
$\MA(W,W)$. Thus the maps $\ev$ and $\coev$ satisfy all the conditions 
required for left duals in a monoidal category.
\endproof

Propositions 10.11 and 10.12 complete the proof of the implication 
(a)$\,\Rar\,$(b) in Theorem 10.1. At the end of this section we prove yet 
another result which compares the functors $\Phi$ and $\Phi'$ of Theorem 8.2.

By Lemma 10.9 each rationally generated object of $\HdMA$ is rationally 
generated also as an object of $\HdAM$. It follows then from Corollary 8.5 
that rational extendibility of objects is also preserved. We get thus a functor
$$
\HdxratMA\to\HdxratAM.\eqno(10.12)
$$
It will be clear soon that this is an isomorphism of categories.
Recall that the factor bialgebra $C=H/I$ is a Hopf algebra with bijective 
antipode by Corollary 10.2.

\proclaim
Proposition 10.13.
There is a commutative, up to a natural equivalence, diagram
$$
\diagram{
\HdxratMA&\lmapr4\Phi&\MC\cr
\noalign{\smallskip}
\lmapd{18}{}{}&&\lmapd{18}{}{S^{-1}}\cr
\noalign{\smallskip}
\HdxratAM&\lmapr4{\Phi'}&\CM\cr
}
$$
where the right vertical arrow represents the functor which transforms each 
right $C$-comodule $V$ into a left comodule by using the structure map 
$$
V\to C\ot V,\qquad v\mapsto\sum S^{-1}(v\1)\ot v\0.
$$
\endproclaim

\Proof.
Given $W\in\HdxratMA$, we have $W\ot_AQ\in\HdratMC$, and also $Q\ot_A\!W$ 
is an object of $\HdratCM$ since $W$ is rationally extendible as an object of 
$\HdAM$. By Lemma 10.5 both $W\ot_AQ$ and $Q\ot_AW$ are in fact objects of 
$\HdQMQ$ with the bimodule structures satisfying (10.4). Since the left action 
of $A$ on $W$ is the restriction of the left action of $Q$ on $W\ot_AQ$, there 
is a left $Q$-linear map
$$
Q\ot_AW\to W\ot_AQ\eqno(10.13)
$$
such that $1\ot_Ax\mapsto x\ot_A1$ for all $x\in W$. This map is also 
$\Hd$-linear, and therefore a morphism in $\HdratQM$. By Lemma 10.5 it is then 
right $Q$-linear too.

By Theorem 8.2 and the $\HdratCM$-variant of Corollary 10.7 
$\,W\cong\lco\scrC(Q\ot_AW)$ is stable under the right action of $A$ on 
$Q\ot_AW$ obtained by restriction of the right action of $Q$, and then the 
induced right action of $A$ on $W$ must coincide with the initial one since 
this holds in $W\ot_AQ$. It follows that there exists a right $Q$-linear map 
$W\ot_AQ\to Q\ot_AW$ such that $x\ot_A1\mapsto 1\ot_Ax$ for all $x\in W$. 
Clearly this map is the inverse of (10.13). This argument shows also that the 
functor
$$
\HdxratAM\to\HdxratMA
$$
constructed similarly to (10.12) is the inverse of (10.12).

Thus we have proved that (10.13) is an isomorphism in $\HdratMQ$ for each $W$. 
It induces a $k$-linear bijection
$$
\io_W:\quad\Phi'(W)=(Q\ot_AW)^{\Hd}\longrightarrow(W\ot_AQ)^{\Hd}=\Phi(W).
$$
We have to check that the left $C$-comodule structure on $\Phi'(W)$ 
corresponds under $\io_W$ to the right $C$-comodule structure on $\Phi(W)$ 
transformed by $S^{-1}$. It suffices to do this for the objects $U\ot A$ 
associated with right $H$-comodules $U$. By naturality of $\io_W$ the desired 
correspondence of comodule structures will then hold for subobjects of factor 
objects of such objects, and by Corollary 8.5 for all objects.

So let $W=U\ot A$. For convenience we identify $U$ with the $\Hd$-submodule 
$U\ot k$ of $W$, and we identify $W$ with its canonical images in $W\ot_AQ$ 
and $Q\ot_AW$. This allows us to write elements omitting the tensor product 
signs. By Lemma 8.6 there is an isomorphism $\ph:U\to\Phi(W)$ in $\MC$ such 
that
$$
\ph(u)=\sum u\0\,S(u\1)\in W\ot_AQ,\qquad u\in U.
$$
For the left $Q$-module structure on $W\ot_AQ$ defined by formula (10.4) we have
$$
qu=\sum\,q\,\ph(u\0)u\1=\sum\,\ph(u\0)\,qu\1,\qquad q\in Q,\ u\in U.
$$
Hence
$$
\sum\,S^{-1}(u\1)\,u\0=\sum\,\ph(u\0)\,S^{-1}(u\2)u\1=\ph(u)
$$
for all $u\in U$. It follows that $\io_W^{-1}\circ\ph=\ph'$ where 
$\ph':U\to\Phi'(W)$ is the $k$-linear map defined by the formula
$$
\ph'(u)=\sum\,S^{-1}(u\1)\,u\0\in Q\ot_AW,\qquad u\in U.
$$
The left $\scrC$-comodule structure 
$\la:Q\ot_AW\to\scrC\ot_Q\!(Q\ot_AW)$ is defined by the formula
$$
\la(qx)=qg\ot_Qx,\qquad q\in Q,\ x\in W
$$
where $g$ is the distinguished grouplike of $\scrC$. In particular, 
this formula holds for $x\in U$. We get
$$
\eqalign{
\la\bigl(\ph'(u)\bigr)=\sum\,S^{-1}(u\1)g\ot_Qu\0
&=\sum\,S^{-1}(u\3)gu\2\ot_QS^{-1}(u\1)u\0\cr
&=\sum\,\psi\bigl(S^{-1}(u\1)\bigr)\ot_Q\ph'(u\0)
}
$$
where $\,\psi:H\to\scrC^{\Hd}$ is the homomorphism of coalgebras defined 
in Corollary 4.11 which induces the canonical isomorphism $C\cong\scrC^{\Hd}$. 
The induced left $C$-comodule structure on $\Phi'(W)$ is given by the map
$$
\Phi'(W)\to C\ot\Phi'(W),\qquad
\ph'(u)\mapsto\sum\,S^{-1}(\overline u\1)\ot\ph'(u\0),
$$
where $\overline h=h+I$ stands for the image of an element $h\in H$ in the 
factor Hopf algebra $C=H/I$. Under the bijection $\io_W$ this corresponds 
precisely to the transformation of the right $C$-comodule structure on 
$\Phi(W)$ by means of $S^{-1}$.
\endproof

\section
11. Conormal factor coalgebras

The Abe-Kanno theorem in \cite{Ab-K59} includes the statement that normal 
subgroups of a connected algebraic group $G$ correspond to subfields of the 
field $k(G)$ of rational functions on $G$ which are both left and right 
invariant. An analog of this conclusion in the Hopf algebraic setup of our 
paper is also true.

There is the notion of conormal quotients of a Hopf algebra $H$ (see 
\cite{Schn93} and \cite{Tak94}). Consider the \emph{left adjoint coaction} of 
$H$ on itself defined by means of the map
$$
\adc:H\to H\ot H,\qquad x\mapsto\sum\,x_1S(x\3)\ot x\2\quad\text{for $x\in H$}.
\eqno(11.1)
$$
A factor coalgebra $C=H/I$ is \emph{left conormal} if the corresponding 
coideal $I$ of $H$ is stable under the left adjoint coaction. If $C$ 
is a left conormal left $H$-module factor coalgebra of $H$, then the 
corresponding right coideal subalgebra $\lco CH$ of $H$ is a Hopf subalgebra 
by \cite{Tak94, Prop. 1.4}.

Assuming that $H$ satisfies conditions (A1) and (A2), we will describe in 
Theorem 11.2 the subalgebras of the quotient ring $Q$ of $H$ corresponding to 
the left conormal left $H$-module factor coalgebras of $H$ under the bijection 
of Theorem 0.1. This result is dual to Theorem 10.1, but it has a much shorter 
proof.

Here we will need the extension to $Q$ of both the left and right actions of 
$\Hd$ on $H$ defined by formulas (1.2). The obvious interrelation between the 
two actions leads to the following fact which will be used in the proof of 
Theorem 11.2.

\proclaim
Lemma 11.1.
The antiautomorphism $S$ of $Q$ extending the antipode of $H$ gives a 
bijection between the left $\Hd$-invariant and the right $\Hd$-invariant 
subalgebras of $Q$.
\endproclaim

\Proof.
The dual Hopf algebra $\Hd$ has a bijective antipode which sends each linear 
function $f\in\Hd$ to the composite $fS$ of $f$ and the antipode of $H$. 
The conclusion of Lemma 11.1 will follow from the identity
$$
S(x)\lhu f=S\bigl((fS)\rhu x\bigr),\qquad x\in Q,\ \,f\in\Hd.\eqno(11.2)
$$
To check it note first that for $x\in H$ both sides are equal to 
$\,\sum f\bigl(S(x\2)\bigr)S(x\1)$. Viewing $\Hom_k(\Hd\!,Q)$ as an algebra 
with respect to the convolution multiplication, the two maps 
$\,\al,\be:\,Q\to\Hom_k(\Hd\!,Q)\,$ defined by the formulas
$$
\al(x)(f)=S(x)\lhu f,\qquad\be(x)(f)=S\bigl((fS)\rhu x\bigr)
$$
for $x\in Q$ and $f\in\Hd$ are algebra antihomomorphsisms which agree on $H$. 
Hence $\al=\be$ everywhere, i.e., (11.2) holds for all $x\in Q$.
\endproof

\setitemsize(a)
\proclaim
Theorem 11.2.
For a coideal left ideal $I$ of $H$ and the corresponding subalgebra $A$ of 
$Q$ the following conditions are equivalent:

\item(a)
\ $I$ is stable under the left adjoint coaction of $H,$ i.e., 
$\adc(I)\sbs H\ot I,$

\item(b)
\ $A$ is both left and right $\Hd$-invariant,

\item(c)
\ $S(A)=A$ where $S:Q\to Q$ is the extension of the antipode of $H$.

\endproclaim

\Proof.
Consider $Q\ot Q$ as an $\Hd$-bimodule with respect to the left $\Hd$-module 
structure defined in section 1 and the tensor product of the right $\Hd$-module 
structures on two copies of $Q$ given similarly by the action $\lhu$. 
Let $\scrI$ be the left $\Hd$-invariant coideal of the $Q$-coring $Q\ot Q$ 
corresponding to $I$.

Let $\psi:H\to Q\ot Q$ be the map defined in Lemma 1.2. 
Since $\,\scrI=Q{\mskip2mu}\psi(I){\mskip1mu}Q\,$ and $\,\psi(I)=\scrI^{\Hd}\!$ 
(the subspace of invariants with respect to the left action of $\Hd$), we 
deduce that $\scrI$ is a right $\Hd$-submodule of $Q\ot Q$ if and only if so 
is $\psi(I)$. The subspace $H\ot H$ is a rational right $\Hd$-submodule of 
$Q\ot Q$. The right action of $\Hd$ on this subspace corresponds to the left 
$H$-comodule structure 
$$
\la:H\ot H\to H\ot(H\ot H)
$$
obtained as the tensor product of the comodule structures on two copies of $H$ 
arising from the comultiplication. Note that
$$
\eqalign{
\la\bigl(\psi(x)\bigr)&=\sum\,\la\bigl(x\1\ot S(x\2)\bigr)\cr
&=\sum\,x_1S(x\4)\ot\bigl(x\2\ot S(x\3)\bigr)\cr
&=\sum\,x_1S(x\3)\ot\psi(x\2)=(\id\ot\psi)\bigl(\adc(x)\bigr)
}
$$
for each $x\in H$. The condition that $\psi(I)$ is stable under the right 
action of $\Hd$ means precisely that $\la\bigl(\psi(I)\bigr)\sbs H\ot\psi(I)$. 
Since $\psi$ is injective, this is equivalent to condition (a) of Theorem 11.2.

The map $\,\de:Q\to Q\ot Q\,$ in the proof of Lemma 1.3 is left and right 
$\Hd$-linear. Since $\,\scrI=Q{\mskip2mu}\de(A){\mskip1mu}Q\,$ and 
$\,A=\de^{-1}(\,\scrI\,)$, we deduce that $\scrI$ is a right $\Hd$-submodule 
of $Q\ot Q$ if and only if so is $A$. Thus both (a) and (b) are equivalent to 
the condition that $\scrI$ is right $\Hd$-invariant, and therefore 
(a)$\,\Lrar\,$(b).

In view of Lemma 11.1 condition (c) implies (b). Conversely, suppose that (b) 
holds. Then (a) is also true. By Lemma 11.1 $S(A)$ is a left $\Hd$-invariant 
subalgebra of $Q$. Under the bijection of Theorem 1.1 $S(A)$ corresponds to
a coideal left ideal, say $I'$, of $H$. For each $x\in I$ we have 
$\,\sum\,x\1\ot_A\!S(x\2)=0\,$ in $\,Q\ot_A\!Q\,$ by the description of 
correspondence in section 1, and therefore, by (11.1),
$$
\sum\,x_1S(x\4)\ot\bigl(x\2\ot_A\!S(x\3)=0\quad\text{in $\,H\ot(Q\ot_A\!Q)$}
$$
since $\adc(x)\in H\ot I$. There is a well-defined $k$-linear map
$\,Q\ot_A\!Q\to Q\ot_{S(A)}\!Q\,$ such that 
$\,x\ot y\mapsto S(y)\ot_{S(A)}S(x)\,$ 
for all $x,y\in Q$. Applying it, we get
$$
\sum\,x_1S(x\4)\ot\bigl(S^2(x\3)\ot_{S(A)}S(x\2)\bigr)=0\quad
\text{in $\,H\ot(Q\ot_{S(A)}Q)\,$}
$$
and, making use of the left $Q$-module structure,
$$
\sum\,x\1\ot_{S(A)}S(x\2)
=\sum\,x_1S(x\4)S^2(x\3)\ot_{S(A)}S(x\2)=0\quad
\text{in $\,Q\ot_{S(A)}Q$}
$$
for each $x\in I$. Hence $I\sbs I'$, again by the description of 
correspondence between coideal left ideals of $H$ and left $\Hd$-invariant 
artinian subalgebras of $Q$. Since this correspondence preserves inclusions, 
we infer that $A\sbs S(A)$. This inclusion applied to the left and right 
$(H\op)^\circ$-invariant artinian subalgebra $A\op$ of $Q\op$ shows that 
$A\sbs S^{-1}(A)$ as well. Hence (b)$\,\Rar\,$(c).
\endproof

\section
12. Comparison with the classical case

For a connected algebraic group $G$ the subfield of the field of rational 
functions $k(G)$ corresponding to a closed subgroup $K\sbs G$ in the Abe-Kanno 
theorem is obtained very easily by taking invariants for the action of $K$ on 
$k(G)$ induced by left translations of $G$. Our correspondence in Theorem 0.1 
is much less obvious as it involves intermediate corings. Proposition 12.1 
provides an alternative description of this correspondence for a certain class 
of left $H$-module factor coalgebras of the Hopf algebra $H$. It shows that in 
the case when $H$ represents an affine algebraic group $G$ our construction 
fully agrees with that of Abe and Kanno \cite{Ab-K59}.

Given a left $H$-module factor coalgebra $C=H/I$ of $H$, the canonical 
surjective homomorphism of coalgebras $H\to C$ induces an injective 
homomorphism of dual algebras $C^*\to H^*$. Denote by $\Cd$ the vector 
subspace of $C^*$ consisting of all linear functions $C\to k$ vanishing on 
an $H$-submodule of finite vector space codimension in $C$.

Since each left ideal of finite codimension in $H$ contains a two-sided ideal 
of finite codimension, it is clear that $\Cd$ is mapped to $\Hd$, and thus 
the image of $\Cd$ consists of all linear functions $f\in\Hd\!$ such that 
$\,I\sbs\Ker f$. Since $I$ is a coideal left ideal of $H$, this image of $\Cd$ 
is in fact a left coideal subalgebra of the dual Hopf algebra $\Hd$. In 
particular, $\Cd$ is a subalgebra of $C^*$. If $\Cd$ is dense in $C^*$, then 
$\Cd\!$ completely determines the coalgebra $C$, and therefore also the 
corresponding subalgebra of the quotient ring $Q$ of $H$. Moreover, any left 
coideal subalgebra of $\Hd$ which is dense in the image of $\Cd$ can be used 
for this purpose. Here we employ the right action of $\Hd$ on $Q$ which 
extends the right action on $H$ defined in (1.2).

\proclaim
Proposition 12.1.
For a left coideal subalgebra $L$ of the dual Hopf algebra $\Hd$ put
$$
\displaylines{
L^\perp=\{x\in H\mid\ \ell(x)=0\,\text{ for all }\ell\in L\},\cr
\noalign{\smallskip}
\LQ=\{x\in Q\mid\ x\lhu\ell=\ell(1)x\,\text{ for all }\ell\in L\}
}
$$
Then $L^\perp$ is a coideal left ideal of $H$ and $\LQ$ is a left 
$\Hd$-invariant artinian subalgebra of $Q$ which correspond to each other 
under the bijection of Theorem 1.1.
\endproclaim

\Proof.
The subspace $L^\perp$ is a coideal of $H$ since $L$ is a subalgebra of $\Hd$, 
and $L^\perp$ is a left ideal of $H$ since $L$ is a left coideal of $\Hd$. 
For each linear function $\ell\in\Hd$ define a right $Q$-linear map 
$\,\al_\ell:Q\ot Q\to Q\,$ by the rule
$$
\al_\ell(x\ot y)=(x\lhu\ell){\mskip1mu}y,\qquad x,y\in Q,\eqno(12.1)
$$
and put
$$
\scrI=\bigcap_{\ell\in L}\,\Ker\al_\ell\sbs Q\ot Q\,.\eqno(12.2)
$$
We consider $Q\ot Q$ as a left $\Hd$-module with respect to the diagonal action 
of $\Hd$. Since the right action $\lhu$ of $\Hd$ on $Q$ commutes with the left 
action $\rhu$, it follows that the map $\al_\ell$ is left $\Hd$-linear, and 
therefore $\Ker\al_\ell$ is an $\HdMQ$-subobject of $Q\ot Q$, for each 
$\ell\in\Hd$. Since $Q$ is a right $\Hd$-module algebra with respect to $\lhu$, 
we have the identity
$$
\al_\ell(qt)=\sum\,(q\lhu\ell\1)\,\al_{\ell\2}(t),\qquad 
q\in Q,\ \,t\in Q\ot Q.\eqno(12.3)
$$
If $\ell\in L$, then $\sum\ell\1\ot\ell\2\in\Hd\ot L$ since $L$ is a left 
coideal subalgebra, whence $\al_\ell(qt)=0$ for all $q\in Q$ and $t\in\scrI$. 
This shows that $\scrI$ is also a left $Q$-submodule, and therefore an 
$\HdQMQ$-subobject of the left $\Hd$-module $Q$-coring $Q\ot Q$.

Proposition 4.2 ensures that $\scrI$ is generated as a $Q$-bimodule by its 
subspace of $\Hd$-invariants $\scrI^{\Hd}\!$. Recall the bijection 
$\psi:H\to(Q\ot Q)^{\Hd}\!$ defined in Proposition 4.10 by formula (4.7). 
If $x\in H$, then
$$
\al_\ell\bigl(\psi(x)\bigr)=\sum\,(x\1\lhu\ell)\,S(x\2)
=\sum\,\ell(x\1)\,x\2S(x\3)=\ell(x).\eqno(12.4)
$$
It follows that $\psi(x)\in\scrI$ if and only if $x\in L^\perp$. Hence $\scrI$ 
is generated by $\psi(L^\perp)$, which means that $\scrI$ is the left 
$\Hd$-invariant coideal of $\,Q\ot Q\,$ corresponding to $\,L^\perp$ under the 
bijection of Theorem 1.1. Next,
$$
\al_\ell(1\ot x-x\ot1)=\ell(1){\mskip1mu}x-(x\lhu\ell)
$$
for each $x\in Q$. It follows that $1\ot x-x\ot1\in\scrI$ if and only if 
$x\in\LQ$. Hence $\LQ$ is the subalgebra of $Q$ corresponding to $\scrI$ under 
the bijection of Theorem 1.1.
\endproof

\proclaim
Corollary 12.2.
Suppose that $L$ and $L'$ are two left coideal subalgebras of $\Hd$. Then
$L^\perp={L'}^\perp$ if and only if $\,\LQ=\LpQ$.
\endproclaim

\Proof.
This is immediately clear from bijectivity of the correspondence established 
by Theorem 1.1.
\endproof

\proclaim
Corollary 12.3.
Suppose that $C$ is a left $H$-module factor coalgebra of $H$ such that $\Cd$ 
is dense in $C^*$. Then $\CdQ$ is the left $\Hd$-invariant artinian subalgebra 
of $Q$ corresponding to $C$ under the bijection of Theorem 0.1.
\endproclaim

\Proof.
Identifying $\Cd$ with its canonical image $L$ in $\Hd$, we will have 
$\,C=H/L^\perp$ by the density assumption.
\endproof

\proclaim
Corollary 12.4.
Let $H$ be a noetherian PI Hopf algebra, finitely generated as an ordinary 
algebra. Then each left $\Hd$-invariant artinian subalgebra of $Q$ coincides 
with $\,\LQ\,$ for some left coideal subalgebra $L$ of $\Hd$.
\endproclaim

\Proof.
It was proved by Anan'in \cite{An92, Lemma 1} that each finitely generated 
right module over a finitely generated right noetherian PI algebra is 
residually finite-dimensional in the sense that its submodules of finite 
vector space codimension have zero intersection. We use this result with the 
right and left sides interchanged. In particular, the Hopf algebra $H$ in the 
hypothesis satisfies our basic assumption (A1), and by \cite{Sk21, Th. 5.5} it 
satisfies also (A2). Moreover, each left $H$-module factor coalgebra $C$ of 
$H$ is residually finite-dimensional as a left $H$-module. This means that 
$\Cd$ is dense in $C^*$. By Corollary 12.3 the subalgebra of $Q$ corresponding 
to $C$ is $\LQ$ where $L$ is the canonical image of $\Cd$ in $\Hd$. By Theorem 
0.1 each left $\Hd$-invariant artinian subalgebra of $Q$ corresponds to some 
left $H$-module factor coalgebra of $H$.
\endproof

Let $C=H/L^\perp$ where $L$ is a left coideal subalgebra of $\Hd$. We may 
identify $\Cd$ with its image in $\Hd$ and $L$ with a dense subalgebra of 
$\Cd$. This allows us to view elements of $L$ as linear functions on $C$. 
In Proposition 12.7 we will describe the functors $\Psi$ and $\Psi'$ of 
Theorem 8.2 in terms of certain actions of $L$. For this result we need a few 
preparations.

The antipode of the Hopf algebra $\Hd$ is given by the assignment 
$\ell\mapsto\ell S$ where $\ell S$ is the composite of $\ell\in\Hd$ and the 
antipode $S$ of $H$. Since $S$ is bijective, the inverse map is given by 
the assignment $\ell\mapsto\ell S^{-1}$. In addition to the already used left 
action $\rhu$ of $\Hd$ on $Q$ we will need another left action $\rhd$ 
defined by the rule
$$
\ell\rhd q=q\lhu(\ell S^{-1}),\qquad\ell\in\Hd\!,\ \,q\in Q.\eqno(12.5)
$$
It commutes with the action $\rhu$ and makes $Q$ into a left 
$(\Hd)\cop$-module algebra.

Let $\scrC=(Q\ot Q)/\,\scrI$ be the left $\Hd$-module factor coring of the 
$Q$-coring $Q\ot Q$ corresponding to the left $H$-module factor coalgebra $C$ 
of $H$ under the bijection of Theorem 1.1, and let $\psi:C\to\scrC^{\Hd}$ be 
the canonical isomorphism of left $H$-module coalgebras defined in Corollary 
4.11.

\proclaim
Lemma 12.5.
For each $\ell\in L$ there are a right $Q$-linear map $\al_\ell:\scrC\to Q$ 
and a left $Q$-linear map $\be_\ell:\scrC\to Q$ such that
$$
\eqalignno{
\al_\ell\bigl(\psi(x)\bigr)&=\ell(x),\qquad
\al_\ell(qt)=\sum\,(q\lhu\ell\1)\,\al_{\ell\2}(t),&(12.6)\cr
\be_\ell\bigl(\psi(x)\bigr)&=\ell(x),\qquad
\be_\ell(tq)=\sum\,\be_{\ell\2}(t)\,(\ell\1\rhd q)&(12.7)
}
$$
for all $x\in C,$ $q\in Q,$ and $t\in\scrC$. Moreover,
$$
\bigcap_{\ell\in L}\,\Ker\al_\ell=\bigcap_{\ell\in L}\,\Ker\be_\ell=0.\eqno(12.8)
$$
\endproclaim

\Proof.
The isomorphism $C\to\scrC^{\Hd}\!$ is induced by the isomorphism 
$H\to(Q\ot Q)^{\Hd}$ of Proposition 4.10. So there is a commutative diagram
$$
\diagram{
H&\lmapr2\psi&(Q\ot Q)^{\Hd}\cr
\noalign{\smallskip}
\mapd{}{}&&\mapd{}{}\hphantom{{}^{\Hd}}\cr
\noalign{\smallskip}
C&\lmapr4\psi\hidewidth&\scrC^{\Hd}\cr
}
$$
where both vertical arrows are the canonical surjections. By abuse of notation 
the same letter $\psi$ is used here to denote two different but related maps.

Since the right $Q$-linear map $\,Q\ot Q\to Q\,$ defined by formula (12.1) 
vanishes on the coideal $\scrI$ in view of (12.2), it induces a right 
$Q$-linear map $\,\scrC\to Q$. Identities (12.3) and (12.4) for the former map 
amount to (12.6) for the latter, while (12.2) shows that the first 
intersection in (12.8) is zero.

Similarly, we first consider $\be_\ell$ as a left $Q$-linear map 
$\,Q\ot Q\to Q\,$ defined by the rule
$$
\be_\ell(x\ot y)=x{\mskip1mu}(\ell\rhd y),\qquad x,y\in Q.\eqno(12.9)
$$
It is a morphism in $\HdQM$ such that
$$
\be_\ell(tq)=\sum\,\be_{\ell\2}(t)\,(\ell\1\rhd q)\quad
\text{for all $q\in Q$ and $t\in Q\ot Q$}
$$
and
$$
\be_\ell\bigl(\psi(x)\bigr)
=\sum x\1\,\bigl(\ell\rhd S(x\2)\bigr)
=\sum \ell(x\3)\,x\1\,S(x\2)=\ell(x)
$$
for all $x\in H$. As in the proof of Proposition 12.1, it follows that the 
intersection
$$
\scrI'=\bigcap_{\ell\in L}\,\Ker\be_\ell
$$
is an $\HdQMQ$-subobject of $\,Q\ot Q\,$ such that 
$\scrI'^{{\mskip1mu}\Hd}\!=\psi(L^\perp)$. Hence $\scrI'$ is the left 
$\Hd$-invariant coideal of the left $\Hd$-module $Q$-coring $\,Q\ot Q\,$ 
corresponding to $L^\perp$ under the bijection of Theorem 1.1, i.e., 
$\scrI'\!=\scrI$. As a consequence, the map defined by (12.9) induces a left 
$Q$-linear map $\scrC\to Q$ satisfying (12.7) and (12.8).
\endproof

Given a right $\scrC$-comodule $M$ with structure map $\,\rho:M\to M\ot_Q\scrC$, 
we define the action of a left $Q$-linear map $\,f:\scrC\to Q\,$ on $M$ by the 
rule
$$
f\rhu m=(\id_M\ot f)\bigl(\rho(m)\bigr),\qquad m\in M.
$$
These actions make $M$ into a left module over the left dual algebra of $\scrC$ 
consisting of all left $Q$-linear maps $\scrC\to Q$ (see \cite{Br-W, 19.4}). 
Similarly, each right $Q$-linear map $f:\scrC\to Q$ acts on a left 
$\scrC$-comodule $M$ by the rule
$$
m\lhu f=(f\ot\id_M)\bigl(\la(m)\bigr),\qquad m\in M
$$
where $\,\la:M\to\scrC\ot_QM\,$ is the comodule structure map.

\proclaim
Lemma 12.6.
Let $C=H/L^\perp$ where $L$ is a left coideal subalgebra of $\Hd,$ and let $\scrC$ 
be the left $\Hd$-module factor coring of $Q\ot Q$ corresponding to $C$. Then
$$
\eqalign{
M\co\scrC&=\{\,m\in M\mid\ \be_\ell\rhu m=\ell(1)m\,
\text{ for all }\ell\in L\,\}\cr
\noalign{\smallskip}
\lco\scrC M&=\{\,m\in M\mid\ m\lhu\al_\ell=\ell(1)m\,
\text{ for all }\ell\in L\,\}
}
$$
for all objects $M$ of the categories, respectively, $\MscrC$ and $\scrCM$ 
which are projective, respectively, in $\MQ$ and $\QM$.
\endproclaim

\Proof.
Let $M\in\MscrC$. Recall that $\,M\co\scrC=\{m\in M\mid\rho(m)=m\ot_Qg\}\,$ 
where $g$ is the distinguished grouplike of $\scrC$. Since $g=\psi(1_C)$ where 
the element $1_C\in C$ is the image of $1\in H$, we have $\be_l(g)=\ell(1)$ 
for each $\ell\in L$ by (12.7).

Let $m\in M$. Assuming that $M$ is projective in $\MQ$, 
the functor $M\ot_Q{}?$ commutes with arbitrary intersections. Since
$$
(\be_\ell\rhu m)-\ell(1)m=(\id_M\ot\be_\ell)\bigl(\rho(m)-m\ot_Qg\bigr),
$$
the equality $\,\be_\ell\rhu m=\ell(1)m\,$ holds for all $\ell\in L$ 
if and only if
$$
\rho(m)-m\ot_Qg\in\bigcap_{\ell\in L}\,\Ker(\id_M\ot\be_\ell)
=M\ot_Q\bigcap_{\ell\in L}\Ker\be_\ell=0
$$
by (12.8), which means that $m\in M\co C$. The second equality for objects of 
the category $\scrCM$ is proved similarly.
\endproof

Given a right $C$-comodule $V$, the canonical embedding of $L$ in the dual 
algebra $C^*$ allows us to view $V$ as a left $L$-module. We denote the action 
of $L$ on $V$ by the symbol $\rhu$. Now $V\ot Q$ becomes an $(L,\Hd)$-bimodule 
with respect to the left action $\rhu$ of $L$ on $V$ and the right action $\lhu$ 
of $\Hd$ on $Q$.

Similarly, each left $C$-comodule $V$ is a right $L$-module in a natural way, 
and we view $Q\ot V$ as an $(\Hd\!,L)$-bimodule with respect to the left 
action $\rhd$ of $\Hd$ on $Q$ and the natural right action $\lhu$ of $L$ on 
$V$.

\proclaim
Proposition 12.7.
Let $C=H/L^\perp$ where $L$ is a left coideal subalgebra of $\Hd\!$. Then
$$
\Psi(V)=\{\,t\in V\ot Q\mid\ \ell\rhu t=t\lhu\ell\,\text{ for all }\ell\in L\,\}
$$
for each right $C$-comodule $V$ and
$$
\Psi'(V)=\{\,t\in Q\ot V\mid\ \ell\rhd t=t\lhu\ell\,\text{ for all }\ell\in L\,\}
$$
for each left $C$-comodule $V$.
\endproclaim

\Proof.
Let $V\in\MC$. Recall that $\Psi(V)=(V\ot Q)\co\scrC$. We claim that
$$
\be_\ell\rhu t=\sum\,\ell\2\rhu t\lhu\ell\1S^{-1}\quad\text{for all 
$\,\ell\in L\,$ and $\,t\in V\ot Q$}.
$$
When checking this equality we may assume that $t=v\ot q$ where $v\in V$ and 
$q\in Q$. Then $\,\rho(t)=\sum v\0\ot\psi(v\1)q\,$ where 
$\,\rho:V\ot Q\to(V\ot Q)\ot_Q\scrC\cong V\ot\scrC\,$ is the respective right 
$\scrC$-comodule structure map. We get, making use of (12.7),
$$
\eqalign{
\be_\ell\rhu t&=\sum v\0\ot\be_\ell\bigl(\psi(v\1)q\bigr)\cr
&=\sum v\0\ot\ell\2(v\1)(\ell\1\rhd q)\cr
&=\sum\,(\ell\2\rhu v)\ot(q\lhu\ell\1S^{-1})
=\sum\,\ell\2\rhu t\lhu\ell\1S^{-1}
}
$$
as claimed. By Lemma 12.6 it follows that $\Psi(V)$ consists of all elements 
$t\in V\ot Q$ such that
$$
\sum\,\ell\2\rhu t\lhu\ell\1S^{-1}=\ell(1){\mskip2mu}t\quad\text{ for all }\ell\in L.
$$
This condition on $t$ is equivalent to the condition in the statement of 
Proposition 12.7 since 
$$
\sum\,\ell\2\,(\ell\1S^{-1})=\sum\,(\ell\2S^{-1})\,\ell\1=\ell(1){\mskip2mu}\ep
$$
in the Hopf algebra $\Hd$ by the properties of the antipode. Here $\ep$ is the 
identity element of $\Hd$, and so $\,t\lhu\ep=t$.

For a left comodule $V\in\CM$ we have $\Psi'(V)=\lco\scrC(Q\ot V)$ and
$$
t\lhu\al_\ell=\sum\,(\ell\1S)\rhd t\lhu\ell\2\quad\text{for all 
$\,\ell\in L\,$ and $\,t\in Q\ot V$}.
$$
Indeed, if $\,t=q\ot v$, then $\,\la(t)=\sum q\,\psi(v\ng)\ot v\0\,$ for the 
left $\scrC$-comodule structure 
$\,\la:Q\ot V\to\scrC\ot_Q(Q\ot V)\cong\scrC\ot V\,$ on $Q\ot V$, whence
$$
\eqalign{
t\lhu\al_\ell&=\sum\,\al_\ell\bigl(q\,\psi(v\ng)\bigr)\ot v\0\cr
&=\sum\,(q\lhu\ell\1)\,\ell\2(v\ng)\ot v\0\cr
&=\sum\,\bigl((\ell\1S)\rhd q\bigr)\ot(v\lhu\ell\2)
=\sum\,(\ell\1S)\rhd t\lhu\ell\2.
}
$$
Thus, an element $\,t\in Q\ot V$ lies in $\Psi'(V)$ if and only if
$$
\sum\,(\ell\1S)\rhd t\lhu\ell\2=\ell(1){\mskip2mu}t\quad\text{for all 
$\,\ell\in L\,$ and $\,t\in Q\ot V$},
$$
and this is equivalent to the condition in the statement of Proposition 12.7.
\endproof

Proposition 12.7 shows that for the evaluation of $\Psi$ and $\Psi'$ one may 
take $L$ to be any left coideal subalgebra of $\Hd$ which is dense in $C^*$. 
For example, suppose that $H=k[G]$ is the commutative Hopf algebra representing 
a reduced affine algebraic group $G$ of finite type over an algebraically closed 
field $k$, and $C=k[K]$ is its factor algebra representing a reduced closed 
subgroup $K$ of $G$. The elements of $K$ are associated with algebra 
homomorphisms $C\to k$, and in this way the group algebra of $K$ embeds in 
$\Cd$ as a dense subalgebra. Let $V$ be a rational $K$-module. The elements of 
$V\ot k(G)$ are rational maps $f:G\to V$ defined on dense Zariski open subsets 
of $G$. The maps $f$ that constitute the $\Hd$-module $\Psi(V)$ are 
characterized by the additional condition
$$
f(ax)=a\,f(x)\qquad a\in K,\ x\in G,
$$
and the action of $G$ on $\Psi(V)$ is induced by right translations of $G$. 
The rational $G$-module $\Rat\bigl(\Psi(V)\bigr)$ consists of those maps 
$f\in\Psi(V)$ that are regular on the whole group $G$. This is the rational 
$G$-module \emph{induced} from $V$. The elements of $\Cd$ are known as 
distributions with finite support on $K$. If $K$ is connected, then 
all distributions with support at the identity element of the group form a 
dense Hopf subalgebra of $\Cd$. This subalgebra can be used for an alternative 
description of the functor $\Psi$.

If $U$ is a right $H$-comodule regarded as a right $C$-comodule with respect 
to the canonical homomorphism of coalgebras $H\to C$, then 
$\Psi(U)\cong U\ot A$ where $A=Q\co\scrC$ since $\Phi(U\ot A)\cong U$ by Lemma 
8.6. More precisely, considering the automorphism $\nu$ of the right 
$Q$-module $U\ot Q$ defined by the rule
$$
\nu(u\ot q)=\sum u\0\ot u\1q,\qquad u\in U,\ q\in Q,
$$
we have $\,\Psi(U)=(U\ot Q)\co\scrC=\nu(U\ot A)$. It follows that
$$
\Psi(V)=(V\ot Q)\co\scrC=(V\ot Q)\cap\nu(U\ot A),
$$
and therefore
$$
\Psi(V)\cong\nu^{-1}(V\ot Q)\cap(U\ot A)\sbs U\ot A
$$
for each $C$-subcomodule $V\sbs U$. The map $\Psi(V)\to U\ot A$ provided by 
this $k$-linear bijection is a monomorphism in $\HdMA$.

A striking similarity to the operation $\natural$ introduced by Moeglin and 
Rentschler \cite{Moe-R86, I.4} is not accidental. The only difference lies in 
the fact that these authors work with left $G$-invariant, i.e., right 
$\Hd$-invariant subalgebras of the ring $k(G)$ of rational functions. Put 
$B=S(A)$ where $S$ is the antiautomorphism of $Q$ induced by the antipode of 
$H$. Then
$$
\mu=(\id\ot S)\circ\nu\circ(\id\ot S^{-1})
$$
is an automorphism of $U\ot Q$, now regarded as a left $Q$-module, such that
$$
V^\natural=\mu^{-1}(V\ot Q)\cap(U\ot B)\sbs U\ot B
$$
is the $B$-submodule defined in \cite{Moe-R86, I.4} in the case when $H=k[G]$, 
$C=k[K]$, and $B=k(G/K)$ is the ring of rational functions on the quotient 
$G/K$.

So defined $V^\natural$ makes sense for an arbitrary Hopf algebra satisfying 
assumptions (A1) and (A2), and $V^\natural$ is also an $\Hd$-submodule of the 
tensor product $U\ot B$ where $\Hd$ acts naturally on $U$, while $B$ is a left 
$\Hd$-module with respect to the action $\rhd$ defined in (12.5). The two 
module structures make $V^\natural$ an object of the category 
$(\Hd)\cop\hbox{-}\BM$. There is an isomorphism of $\Hd$-modules
$$
\Psi(V)\cong V^\natural,
$$
and the right action of $A$ on $\Psi(V)$ is transferred by means of the 
antiisomorphism of algebras $S|_A:A\to B$ from the left action of $B$ on 
$V^\natural$. The category equivalence $\Psi$ allows us to generalize the 
fundamental fact established in \cite{Moe-R86, I.4}: the assignment 
$V\mapsto V^\natural$ gives a bijection between the set of $C$-subcomodules 
$V\sbs U$ and the set of $\Hd$-invariant left $B$-submodules of $\,U\ot B$.

\section
13. Quasiprojective homogeneous spaces

Hopf algebraic interpretation of quasiprojective homogeneous spaces proposed 
in \cite{Sk10} makes use of graded subalgebras of the Laurent polynomial ring 
$H[t,t^{-1}]$ stable under the natural right coaction of $H$. Such subalgebras 
correspond to a certain class of left $H$-module factor coalgebras of $H$. In 
this section we will describe the corresponding subalgebras of the quotient 
ring $Q(H)$ of $H$ obtained by the bijection of Theorem 0.1.

The essential argument is given in the next lemma. For each right coideal 
$U$ of $H$ we denote by $U^+$ its subspace consisting of all elements $x\in U$ 
such that $\ep(x)=0$.

\proclaim
Lemma 13.1.
Let $U$ be a right coideal of $H$ containing a nonzerodivisor $w$ of $H$. The 
left $\Hd$-invariant artinian subalgebras of $Q(H)$ corresponding to the 
coideal left ideals $HU^+\!$ and $HS^{-1}(U^+)$ are the dominions of the 
subalgebras $R$ and $R'$ of $Q(H)$ generated, respectively, by the sets
$$
\{w^{-1}x\mid\,x\in U\}\qquad\text{and}\qquad\{xw^{-1}\mid\,x\in U\},
$$
\endproclaim

\Proof.
Note that $U^+$ is a coideal of $H$, and therefore $HU^+$ is indeed a coideal 
left ideal. The corresponding left $\Hd$-invariant coideal $\scrI$ of the 
canonical $Q$-coring $Q\ot Q$ is the $Q$-subbimodule generated by $\psi(U^+)$ 
where $\psi:H\to Q\ot Q$ is the map defined in Lemma 1.2.

Replacing $w$ with its scalar multiple we may assume that $\ep(w)=1$. 
There is an invertible endomorphism $\xi$ of the right $Q$-module $U\ot Q$ 
such that
$$
\xi(x\ot q)=\psi(x)q=\sum\,x\0\ot S(x\1)q,\qquad
\xi^{-1}(x\ot q)=\sum\,x\0\ot x\1q
$$
for all $x\in U$ and $q\in Q$. Define endomorphisms $\eta$ and $\ze$ of the 
same $Q$-module by the rules
$$
\eta(x\ot q)=\ep(x)\,\psi(w)\,q,\qquad
\ze(x\ot q)=w\ot w^{-1}xq.
$$
Since $U^+=\{x-\ep(x)w\mid x\in U\}$, we have
$$
\psi(U^+)Q=\xi(U^+\ot Q)=\Img\,(\xi-\eta)
=\Img\,\bigl((\xi-\eta)\circ\xi^{-1}\bigr)
=\Img\,(\Id-\eta\circ\xi^{-1}).
$$
Note that $\,(\eta\circ\xi^{-1})(x\ot1)=\sum\,\eta(x\0\ot x\1)=\psi(w)x$. 
Hence the $Q$-module $\psi(U^+)Q$ is generated by the set
$$
\{x\ot1-\psi(w)x\mid x\in U\}.
$$
It follows that $\,x\ot x^{-1}-\psi(w)\in\psi(U^+)Q\,$ whenever $x\in U$ is 
invertible in $Q(H)$. In particular, $\,w\ot w^{-1}\equiv\psi(w)\,$ modulo 
$\psi(U^+)Q$, and we get
$$
x\ot1-w\ot w^{-1}x\in\psi(U^+)Q\quad
\text{for all $x\in U$}.
$$
Thus $\,\Img\,(\Id-\ze)\sbs\psi(U^+)Q\,$. On the other hand, the $Q$-module
$$
\Img\,(\Id-\ze)=\Img\,\bigl((\Id-\ze)\circ\xi\bigr)=\Img\,(\xi-\ze\circ\xi)
$$
is generated by the set $\,\{\psi(x)-\ep(x)w\ot w^{-1}\mid x\in U\}\,$ since
$$
(\ze\circ\xi)(x\ot1)=\sum\,\ze\bigl(x\0\ot S(x\1)\bigr)=\ep(x)w\ot w^{-1}.
$$
Hence $\psi(x)\in\Img\,(\Id-\ze)$ for all $x\in U^+$, and therefore 
$\Img\,(\Id-\ze)=\psi(U^+)Q$. This implies that $\scrI$ is generated as a 
$Q$-bimodule by the set
$$
\{w^{-1}x\ot1-1\ot w^{-1}x\mid\,x\in U\}.
$$
It follows that the $Q$-coring $(Q\ot Q)/\,\scrI$ is isomorphic to $Q\ot_RQ$. 
For $x\in Q$ the inclusion $\,1\ot x-x\ot1\in\scrI\,$ holds if and only if 
$\,1\ot_Rx=x\ot_R1\,$ in $Q\ot_RQ$, which means that $x$ lies in the dominion of 
$R$ in $Q$. This determines the subalgebra of $Q$ corresponding to $\scrI$ 
under the bijection of Theorem 1.1.

The left $\Hd$-invariant coideal $\scrI'$ of the canonical $Q$-coring $Q\ot Q$ 
corresponding to the coideal left ideal $HS^{-1}(U^+)$ is generated 
as a $Q$-bimodule by the set
$$
\psi\bigl(S^{-1}(U^+)\bigr)=\{\,\sum S^{-1}(x\2)\ot x\1\mid\,x\in U^+\,\}.
$$
Applying the already proved part of Lemma 13.1 with $H$ and $Q$ changed to 
$H\op$ and $Q\op$, we deduce that
$$
\{\,1\ot xw^{-1}-xw^{-1}\ot1\mid\,x\in U\,\}.
$$
is another generating set for this bimodule. Hence 
$(Q\ot Q)/\,\scrI'\cong Q\ot_{R'}\!Q$, and the subalgebra of $Q$ corresponding 
to $\scrI'$ is described similarly.
\endproof

\proclaim
Corollary 13.2.
Suppose that $U$ is a right coideal of $H$ such that $1\in U$. Then the left 
$\Hd$-invariant artinian subalgebra of $Q(H)$ corresponding to the coideal 
left ideal $HU^+\!$ of $H$ is the dominion of the subalgebra of $Q(H)$ 
generated by $U$. Moreover, $HS^{-1}(U^+)=HU^+$.
\endproclaim

\Proof.
By Lemma 13.1 the two coideal left ideals $HU^+$ and $HS^{-1}(U^+)$ correspond 
to the same left $\Hd$-invariant artinian subalgebra of $Q(H)$. Certainly, the 
equality $\,U^+H=S(U^+)H\,$ is actually known for every Hopf algebra by a 
direct proof \cite{Kop93, Lemma 3.1}.
\endproof

Let now $\,A=\bigoplus_{i\in\bbZ}A_it^i\,$ be a graded subalgebra of 
the Laurent polynomial ring $H[t,t^{-1}]$ where each $A_i$ is a right coideal 
of $H$ and $A_1\ne0$. If $A$ has a right artinian classical right quotient 
ring $Q(A)$, then the set of homogeneous nonzerodivisors of $A$ satisfies the 
right Ore condition, and so one can build the ring of fractions $Q\gr(A)$ with 
respect to this Ore set \cite{Sk10, section 3}. Moreover, by \cite{Sk10, Prop. 
3.10} the degree 0 component $Q_0(A)$ of $Q\gr(A)$ is a right artinian 
$\Hd$-simple left $\Hd$-module algebra which embeds in $Q(H)$ as a left 
$\Hd$-invariant subalgebra. By Proposition 3.6 $\,Q_0(A)$ is even two-sided 
artinian.

Similarly, if $A$ has a left artinian classical left quotient ring, then 
there is the graded classical left quotient ring $Q\gr(A)$ of $A$ whose 
component $Q_0(A)$ is again an artinian $\Hd$-simple left $\Hd$-module 
subalgebra of $Q(H)$.

\proclaim
Proposition 13.3.
Let $\,A=\bigoplus_{i\in\bbZ}A_it^i\sbs H[t,t^{-1}]\,$ be a graded subalgebra 
where each $A_i$ is a right coideal of $H$ and $A_1\ne0$. If $A$ has an 
artinian classical either left or right quotient ring, then $Q_0(A)$ is the 
left $\Hd$-invariant artinian subalgebra of $Q(H)$ corresponding to the coideal 
left ideal $\,\sum_{i\in\bbZ}HA_i^+$ or $\,\sum_{i\in\bbZ}HS^{-1}(A_i^+),$ 
respectively.
\endproclaim

\Proof.
Denote by $\Ga$ the set of all integers $i$ such that $A_i$ contains a 
nonzerodivisor of $H$. By \cite{Sk10, Prop. 3.10} there exists an integer $r$ 
such that $i\in\Ga$ for all $i>r$. Note that $A_jA_i^+\sbs A_{i+j}^+$ and 
$A_i^+A_j\sbs A_{i+j}^+$. If $j>0$, then $A_j\ne0$, whence $HA_j=H$ since 
$HA_j$ is a nonzero left ideal and a right coideal of $H$, and it follows that
$$
HA_i^+=HA_jA_i^+\sbs HA_{i+j}^+.
$$
Similarly, $HS^{-1}(A_i^+)\sbs HS^{-1}(A_{i+j}^+)$ for $j>0$ since $HS^{-1}(A_j)=H$. 
Taking $j$ large enough we will also have $i+j>r$. This shows that
$$
\sum_{i\in\bbZ}HA_i^+=\sum_{i\in\Ga}HA_i^+\qquad\text{and}\qquad
\sum_{i\in\bbZ}HS^{-1}(A_i^+)=\sum_{i\in\Ga}HS^{-1}(A_i^+).\eqno(13.1)
$$
The ring $Q_0(A)$ consists of all elements of $Q(H)$ which can be written as 
$w^{-1}x$ in the case of left quotient rings or as $xw^{-1}$ in the case of 
right quotient rings where $x,w\in A_i$ for some $i\in\Ga$ and $w$ is a 
nonzerodivisor of $H$. Since $Q_0(A)$ is an artinian left $\Hd$-invariant
subalgebra of $Q(H)$, it coincides with its own dominion in $Q(H)$ by 
Proposition 3.6. The left $\Hd$-invariant artinian subalgebras of $Q(H)$ 
corresponding to the coideal left ideals $HA_i^+$ and $HS^{-1}(A_i^+)$ for 
each $i\in\Ga$ are determined by Lemma 13.1. In view of this lemma the 
conclusion of Proposition 13.3 follows from equalities (13.1).
\endproof

\proclaim
Corollary 13.4.
If the algebra $A$ in Proposition 13.3 has an artinian classical two-sided 
quotient ring, then $\,\sum_{i\in\bbZ}HA_i^+=\sum_{i\in\bbZ}HS^{-1}(A_i^+)$.
\endproclaim

\section
Appendix. Semiprimary rings arising naturally

Schofield's proof of the fact that dominions of subrings in a semiprimary ring 
are semiprimary \cite{Scho, Th. 7.19} is based on a theorem of Bj\"ork 
\cite{Bj71, Th. 5.1} which asserts that the endomorphism rings of finitely 
presented modules over a semiprimary ring are semiprimary. We will give an 
alternative proof which makes use of basic characterizations of semiprimary 
subrings, also found in \cite{Bj71}.

An element $x$ of a ring $R$ is called \emph{strongly regular} if $Rx=Rx^2$ 
and $xR=x^2R$, and $R$ is called \emph{strongly $\pi$-regular} if for each 
element $a\in R$ there exists an integer $n>0$ such that $a^n$ is a strongly 
regular element of $R$ (see Azumaya \cite{Az54}). It was proved by Dischinger 
\cite{Di76} that the two conditions

\setitemsize(2)
\item(1)
\ for each $a\in R\,$ the descending chain of right ideals $\,a^nR$, $\,n>0$, 
stabilizes,

\item(2)
\ for each $a\in R\,$ the descending chain of left ideals $\,Ra^n$, $\,n>0$, 
stabilizes

\noindent
are equivalent to each other. Hence each of them is equivalent to the ring $R$ 
being strongly $\pi$-regular. In particular, each left or right perfect ring 
is strongly $\pi$-regular. All semiprimary rings are left and right perfect.

The next lemma reformulates part of conclusions in \cite{Az54, Lemma 1}.

\proclaim
Lemma A1.
Suppose that $x$ is a strongly regular element of a ring $R$. Then there are 
uniquely determined elements $\,e,{\mskip1mu}y\in R$ such that
$$
xe=ex=x,\qquad xy=yx=e,\qquad ey=ye=y,\qquad e^2=e.\eqno(14.1)
$$
\endproclaim

\Proof.
The equality $Rx=Rx^2$ implies that the right annihilator $\rann(x)$ of $x$ in 
$R$ coincides with the right annihilator of $x^2$, and therefore 
$xR\cap\rann(x)=0$. The equality $xR=x^2R$ implies that $R=xR+\rann(x)$. Hence
$$
R=xR\oplus\rann(x)
$$
and the left multiplication by $x$ induces a bijective endomorphism $\ph$ of 
the right ideal $xR$. Equalities (14.1) require that $e\in xR$ and 
$y\in eR=xR$. On the other hand, bijectivity of $\ph$ shows that there are 
uniquely determined $e,y\in xR$ such that $xe=x$ and $xy=e$. These elements
do satisfy (14.1). For example, since
$$
\ph(e^2-e)=xe^2-xe=x-x=0,
$$
bijectivity of $\ph$ entails $\,e^2-e=0$.
\endproof

In a strongly $\pi$-regular ring $R$ each nonnil right ideal contains a 
nonzero idempotent. Furthermore, if $R$ has no infinite sets of pairwise 
orthogonal idempotents, then each right ideal is a direct sum of a right ideal 
generated by an idempotent and a nil right ideal. Such a ring $R$ is semilocal 
with nil Jacobson radical $J$. It follows that $R$ is semiprimary when $J$ is 
nilpotent or $R$ is right perfect when $J$ is right $T$-nilpotent, as defined 
by Bass \cite{Bas60}. These additional properties of the Jacobson radical are 
easily verified in the case when $R$ is a subring, respectively, of a 
semiprimary ring \cite{Bj71, Cor. 3.6} or of a right perfect ring \cite{Bj71, 
Prop. 3.5}. We get

\proclaim
Corollary A2.
Suppose that $R$ is a strongly $\pi$-regular subring of a ring $S$. Then $R$ 
is right perfect, or left perfect, or semiprimary whenever so is $S$.
\endproclaim

The property of being strongly $\pi$-regular for a ring passes over to certain 
naturally constructed subrings. The semiprimary case is then inferred as an 
easy consequence.

\proclaim
Theorem A3.
Let $B$ be any subring of a strongly $\pi$-regular ring $S$. The dominion $R$ 
of $B$ in $S$ is a strongly $\pi$-regular ring. Hence $R$ is right perfect, or 
left perfect, or semiprimary whenever so is $S$.

\endproclaim

\Proof.
Suppose that three elements $x\in R$ and $e,y\in S$ satisfy (14.1). Since 
$\,1\ot_Bx=x\ot_B1\,$ in $S\ot_BS$, we have
$$
1\ot_Be=1\ot_Bxy=x\ot_By=ex\ot_By=e\ot_Bxy=e\ot_Be
$$
and $\,e\ot_B1=e\ot_Be\,$ by a similar calculation. Hence $\,1\ot_Be=e\ot_B1$. 
Also,
$$
1\ot_By=1\ot_Bey=e\ot_By=yx\ot_By=y\ot_Bxy=y\ot_Be=ye\ot_B1=y\ot_B1,
$$
all equalities in $S\ot_BS$. Thus $e,y\in R$. In other words, each element 
$x\in R$ which is strongly regular in $S$ has to be strongly regular in $R$. 
Given an arbitrary element $a\in R$, some its power $a^n$ is strongly regular 
in $S$, and therefore in $R$.
\endproof

One of Bj\"ork's results \cite{Bj71, Th. 5.1} asserts that the subring of 
elements fixed by a set of ring endomorphisms of a semiprimary ring is itself 
semiprimary. In \cite{Bj71, Th. 5.2} a similar conclusion is proved for the 
one-sided perfect conditions on a ring and its subring. This is generalized as 
follows:

\proclaim
Theorem A4.
Let $C$ be a coalgebra over a commutative ring $k$ and $\ep:C\to k$ its counit. 
Suppose that a measuring action of $C$ on a strongly $\pi$-regular $k$-algebra 
$A$ is given. Then for any right coideal $U$ of $C$ the set of $U$-invariants
$$
A^U=\{\,x\in A\mid\,cx=\ep(c)x\,\text{ for all $c\in U$}\}
$$
is a strongly $\pi$-regular subalgebra of $A$. Hence the algebra $A^U$ is 
right perfect, or left perfect, or semiprimary whenever so is $A$.
\endproclaim

\Proof.
Clearly, $A^U$ is a $k$-submodule of $A$. If $x\in A^U$ and $u\in U$, then
$$
u(xa)=\sum\,(u\1x)(u\2a)=\sum\,\ep(u\1)x\,(u\2a)=x(ua)
$$
for all $a\in A$. Thus the left multiplication $\mu_x:A\to A$ by any 
$x\in A^U$ commutes with the action of $U$ on $A$. It follows that $A^U$ is 
closed under multiplication in $A$.

As in Theorem A3, it remains to prove that for each triple $x,e,y\in A$ 
satisfying (14.1) the containment $x\in A^U$ implies $e,y\in A^U$. Since $x$ 
is strongly regular, the endomorphism $\ph$ of the right ideal $xA$ induced by 
$\mu_x$ is bijective. Since $\ph$ commutes with the action of $U$, so does the 
inverse endomorphism $\ph^{-1}$. Hence $e=\ph^{-1}(x)$ and 
$y=\ph^{-1}(e)$ are indeed invariant under the action of $U$.
\endproof

Historically the oldest construction of semiprimary rings is provided by 
the rings of endomorphisms.

\proclaim
Theorem A5.
Let $\calA$ be an abelian category. If $M\in\calA$ is an object of finite 
length, then its endomorphism ring\/ $\,\End_{\calA}M$ is semiprimary.
\endproclaim

In the case of module categories this theorem is stated in several textbooks 
on ring theory without any attribution. Actually it is traced to the same 
paper of Fitting one of whose results \cite{Fi33, Satz II} has been generally 
recognized as Fitting's Lemma. Fitting considered not necessarily commutative 
groups with a set of operators satisfying the ascending and descending chain 
conditions on invariant normal subgroups. Abelian groups with a set of 
operators are essentially modules over some ring. By \cite{Fi33, Satz 11} the 
endomorphism ring has a largest nilpotent ideal, called the radical, and by 
\cite{Fi33, Satz 13b} the factor ring by the radical is completely reducible 
(the terminology of \cite{Fi33} uses the word ``automorphism" for what we call 
endomorphism, while ``proper automorphisms" in \cite{Fi33} are bijective 
endomorphisms).

On the other hand, Theorem A5 can be also put in the context of $\pi$-regular 
rings. Indeed, Fitting's Lemma implies that the ring $\,\End_{\calA}M$ is 
strongly $\pi$-regular. This ring has no infinite sets of pairwise 
orthogonal idempotents just because the number of summands in any direct sum 
decomposition of $M$ is bounded by the length of the object $M$. Each nil 
one-sided ideal of $\,\End_{\calA}M$ is nilpotent by the Levitzki-Fitting 
theorem. As we have recalled earlier, these properties characterize a 
semiprimary ring.

\references
\nextref Ab-K59
\auth{E.,Abe;T.,Kanno}
\paper{Some remarks on algebraic groups}
\journal{Tohoku Math.~J.}
\Vol{11}
\Year{1959}
\Pages{376-384}

\nextref Am-M05
\auth{K.,Amano;A.,Masuoka}
\paper{Picard-Vessiot extensions of Artinian simple module algebras}
\journal{J.~Algebra}
\Vol{285}
\Year{2005}
\Pages{743-767}

\nextref An92
\auth{A.Z.,Anan'in}
\paper{Representability of noetherian finitely generated algebras}
\journal{Arch. Math.}
\Vol{59}
\Year{1992}
\Pages{1-5}

\nextref And-F
\auth{F.W.,Anderson;K.R.,Fuller}
\book{Rings and Categories of Modules}
\publisher{Springer}
\Year{1974}

\nextref Az54
\auth{G.,Azumaya}
\paper{Strongly $\pi$-regular rings}
\journal{J.~Fac. Sci. Hokkaido Univ. Ser.~I}
\Vol{12}
\Year{1954}
\Pages{34-39}

\nextref Bas60
\auth{H.,Bass}
\paper{Finitistic dimension and a homological generalization of semiprimary rings}
\journal{Trans. Amer. Math. Soc.}
\Vol{95}
\Year{1960}
\Pages{466-488}

\nextref Ber-Pl
\auth{A.,Berman;R.J.,Plemmons}
\book{Nonnegative Matrices in the Mathematical Sciences}
\publisher{SIAM}
\Year{1994}

\nextref Bia61
\auth{A.,Bia\l ynicki-Birula}
\paper{On the field of rational functions of algebraic groups}
\journal{Pacific J. Math.}
\Vol{11}
\Year{1961}
\Pages{1205-1209}

\nextref Bich23
\auth{J.,Bichon}
\paper{Faithful flatness of Hopf algebras over coideal subalgebras with a bimodule conditional expectation}
arXiv: 2301.05480.

\nextref Bj71
\auth{J.-E.,Bj\"ork}
\paper{Conditions which imply that subrings of semiprimary rings are semiprimary}
\journal{J.~Algebra}
\Vol{19}
\Year{1971}
\Pages{384-395}

\nextref Br-W
\auth{T.,Brzezi\'nski;R.,Wisbauer}
\book{Corings and Comodules}
\publisher{Cambridge Univ. Press}
\Year{2003}

\nextref Cae-MZ
\auth{S.,Caenepeel;G.,Militaru;S.,Zhu}
\book{Frobenius and Separable Functors for Generalized Module Categories and Nonlinear Equations}
\BkSer{Lecture Notes Math.}
\BkVol{1787}
\publisher{Springer}
\Year{2002}

\nextref Ch14
\auth{A.,Chirvasitu}
\paper{Cosemisimple Hopf algebras are faithfully flat over Hopf subalgebras}
\journal{Algebra Number Theory}
\Vol{8}
\Year{2014}
\Pages{1179-1199}

\nextref Ch24
\auth{A.,Chirvasitu}
\paper{Epimorphic quantum subgroups and coalgebra codominions}
\journal{Algebr. Represent. Theory}
\Vol{27}
\Year{2024}
\Pages{219-244}

\nextref Cl-PS77
\auth{E.,Cline;B.,Parshall;L.,Scott}
\paper{Induced modules and affine quotients}
\journal{Math. Ann.}
\Vol{230}
\Year{1977}
\Pages{1-14}

\nextref DG
\auth{M.,Demazure;P.,Gabriel}
\book{Groupes Alg\'ebriques I}
\publisher{Masson}
\Year{1970}

\nextref Di76
\auth{M.F.,Dischinger}
\paper{Sur les anneaux fortement $\pi$-r\'eguliers}
\journal{C.R. Acad. Sci. Paris  S\'er.~A}
\Vol{283}
\Year{1976}
\Pages{571-573}

\nextref Doi83
\auth{Y.,Doi}
\paper{On the structure of relative Hopf modules}
\journal{Comm. Algebra}
\Vol{11}
\Year{1983}
\Pages{243-255}

\nextref Doi92
\auth{Y.,Doi}
\paper{Unifying Hopf modules}
\journal{J.~Algebra}
\Vol{153}
\Year{1992}
\Pages{373-385}

\nextref Er-Sk09
\auth{M.S.,Eryashkin;S.M.,Skryabin}
\paper{The largest Hopf subalgebra of a bialgebra\inRus}
\journal{Mat. Zametki}
\Vol{86}
\Year{2009}
\Pages{942-946};
\etransl{Math. Notes}
\Vol{86}
\Year{2009}
\Pages{887-891}

\nextref Fai
\auth{C.,Faith}
\book{Algebra: Rings, Modules and Categories I}
\publisher{Springer}
\Year{1973}

\nextref Fi33
\auth{H.,Fitting}
\paper{Die Theorie der Automorphismenringe Abelscher Gruppen und ihr Analogon bei nicht kommutativen Gruppen}
\journal{Math. Ann.}
\Vol{107}
\Year{1933}
\Pages{514-542}

\nextref Gab70
\auth{P.,Gabriel}
\paper{G\'en\'eralit\'es sur les groupes alg\'ebriques}
\InBook{Propri\'et\'es G\'en\'erales des Sch\'emas en Groupes}
\BkSer{Lecture Notes Math.}
\BkVol{151}
\publisher{Springer}
\Year{1970}
\Pages{287-317}

\nextref Gant
\auth{F.R.,Gantmacher}
\book{The Theory of Matrices}
\publisher{Chelsea}
\Year{1959}

\nextref Kop93
\auth{M.,Koppinen}
\paper{Coideal subalgebras in Hopf algebras: Freeness, integrals, smash products}
\journal{Comm. Algebra}
\Vol{21}
\Year{1993}
\Pages{427-444}

\nextref Mes06
\auth{B.,Mesablishvili}
\paper{Monads of effective descent type and comomadicity}
\journal{Theory Appl. Categ.}
\Vol{16}
\Year{2006}
\Pages{1-45}

\nextref Moe-R86
\auth{C.,Moeglin;R.,Rentschler}
\paper{Sous-corps commutatifs ad-stables des anneaux de fractions des quotients des alg{\`e}bres enveloppantes; espaces homog{\`e}nes et induction de Mackey}
\journal{J.~Funct. Anal.}
\Vol{69}
\Year{1986}
\Pages{307-396}

\nextref Mo
\auth{S.,Montgomery}
\book{Hopf Algebras and their Actions on Rings}
\publisher{Amer. Math. Soc.}
\Year{1993}

\nextref Mue70
\auth{B.J.,Mueller}
\paper{On semi-perfect rings}
\journal{Ill. J. Math.}
\Vol{14}
\Year{1970}
\Pages{464-467}

\nextref Nas-BO89
\auth{C.,N\u ast\u asescu;M.,Van den Bergh;F.,Van Oystaeyen}
\paper{Separable functors applied to graded rings}
\journal{J.~Algebra}
\Vol{123}
\Year{1989}
\Pages{397-413}

\nextref Nich78
\auth{W.D.,Nichols}
\paper{Quotients of Hopf algebras}
\journal{Comm. Algebra}
\Vol{6}
\Year{1978}
\Pages{1789-1800}

\nextref Ob77
\auth{U.,Oberst}
\paper{Affine Quotientenschemata nach affinen, algebraischen Gruppen und induzierte Darstellungen}
\journal{J.~Algebra}
\Vol{44}
\Year{1977}
\Pages{503-538}

\nextref Scha94
\auth{P.,Schauenburg}
\paper{Hopf modules and Yetter-Drinfel'd modules}
\journal{J.~Algebra}
\Vol{169}
\Year{1994}
\Pages{874-890}

\nextref Schn93
\auth{H.-J.,Schneider}
\paper{Some remarks on exact sequences of quantum groups}
\journal{Comm. Algebra}
\Vol{21}
\Year{1993}
\Pages{3337-3357}

\nextref Scho
\auth{A.H.,Schofield}
\book{Representation of Rings over Skew Fields}
\publisher{Cambridge Univ. Press}
\Year{1985}

\nextref Sk06
\auth{S.,Skryabin}
\paper{New results on the bijectivity of antipode of a Hopf algebra}
\journal{J.~Algebra}
\Vol{306}
\Year{2006}
\Pages{622-633}

\nextref Sk07
\auth{S.,Skryabin}
\paper{Projectivity and freeness over comodule algebras}
\journal{Trans. Amer. Math. Soc.}
\Vol{359}
\Year{2007}
\Pages{2597-2623}

\nextref Sk10
\auth{S.,Skryabin}
\paper{Models of quasiprojective homogeneous spaces for Hopf algebras}
\journal{J.~Reine Angew. Math.}
\Vol{643}
\Year{2010}
\Pages{201-236}

\nextref Sk11
\auth{S.,Skryabin}
\paper{Structure of $H$-semiprime Artinian algebras}
\journal{Algebr. Represent. Theory}
\Vol{14}
\Year{2011}
\Pages{803-822}

\nextref Sk21
\auth{S.,Skryabin}
\paper{Flatness of Noetherian Hopf algebras over coideal subalgebras}
\journal{Algebr. Represent. Theory}
\Vol{24}
\Year{2021}
\Pages{851-875}

\nextref Sk25a
\auth{S.,Skryabin}
\paper{On Takeuchi's correspondence}
arXiv: 2501.06045.

\nextref Sk25b
\auth{S.,Skryabin}
\paper{Failure of flatness over finite-dimensional Hopf subalgebras}
\hfil\break\hbox{}\hfill
arXiv: 2506.16292.

\nextref Sk-Oy06
\auth{S.,Skryabin;F.,Van Oystaeyen}
\paper{The Goldie theorem for $H$-semiprime algebras}
\journal{J.~Algebra}
\Vol{305}
\Year{2006}
\Pages{292-320}

\nextref St
\auth{B.,Stenstr\"om}
\book{Rings of Quotients}
\publisher{Springer}
\Year{1975}

\nextref Sw
\auth{M.E.,Sweedler}
\book{Hopf Algebras}
\publisher{Benjamin}
\Year{1969}

\nextref Sw75
\auth{M.,Sweedler}
\paper{The predual theorem to the Jacobson-Bourbaki theorem}
\journal{Trans. Amer. Math. Soc.}
\Vol{213}
\Year{1975}
\Pages{391-406}

\nextref Tak72
\auth{M.,Takeuchi}
\paper{A correspondence between Hopf ideals and sub-Hopf algebras}
\journal{Manuscripta Math.}
\Vol{7}
\Year{1972}
\Pages{251-270}

\nextref Tak79
\auth{M.,Takeuchi}
\paper{Relative Hopf modules---equivalences and freeness criteria}
\journal{J.~Algebra}
\Vol{60}
\Year{1979}
\Pages{452-471}

\nextref Tak89
\auth{M.,Takeuchi}
\paper{A Hopf algebraic approach to the Picard-Vessiot theory}
\journal{J.~Algebra}
\Vol{122}
\Year{1989}
\Pages{481-509}

\nextref Tak94
\auth{M.,Takeuchi}
\paper{Quotient spaces for Hopf algebras}
\journal{Comm. Algebra}
\Vol{22}
\Year{1994}
\Pages{2503-2523}

\nextref Ulb90
\auth{K.-H.,Ulbrich}
\paper{On Hopf algebras and rigid monoidal categories}
\journal{Isr. J. Math.}
\Vol{72}
\Year{1990}
\Pages{252-256}

\nextref Von96
\auth{N.,Vonessen}
\paper{Actions of algebraic groups on the spectrum of rational ideals}
\journal{J.~Algebra}
\Vol{182}
\Year{1996}
\Pages{383-400}

\endreferences
\bye